\documentclass[11pt]{article}

\usepackage{geometry}
\usepackage{layout}
\usepackage{tikz}

\usepackage{xcolor} 
\usepackage{wrapfig}

\newcommand{\red}[1]{#1}
\newcommand{\green}[1]{#1}
\newcommand{\blue}[1]{#1}

\usepackage{soul}
\setstcolor{red}

\usepackage{url}
\usepackage{amsmath}
\usepackage{amsthm}
\usepackage{amssymb}
\usepackage{bm} %
\usepackage{bbm} %
\usepackage{mathtools}

\usepackage{comment}
\usepackage{enumitem}
\usepackage{nicefrac}

\makeatletter
\let\blx@noerroretextools\@empty
\makeatother

\usepackage[style=numeric, citestyle=authoryear, backend=biber, uniquename=false, uniquelist=false, maxnames=50, firstinits=true, maxcitenames=2]{biblatex}
\usepackage[unicode = true, pdfusetitle]{hyperref}
\usepackage{cleveref}
\hypersetup{ colorlinks = true,linkcolor = blue,anchorcolor =red,citecolor = blue,filecolor = red,urlcolor = red, pdfauthor=author}

\usepackage{subcaption}
\usepackage{diagbox}

\newtheorem{theorem}{Theorem}
\newtheorem{lemma}{Lemma}
\newtheorem{proposition}{Proposition}

\newtheorem{corollary}{Corollary}
\newtheorem{assumption}{Assumption}
\newtheorem{fact}{Fact}

\theoremstyle{definition}
\newtheorem{definition}{Definition}
\newtheorem{remark}{Remark}

\newcommand{\lp}{\left (} %
\newcommand{\rp}{\right )} %
\newcommand{\lb}{\left [}
\newcommand{\rb}{\right ]}
\newcommand{\lbr}{\left\{}
\newcommand{\rbr}{\right\}}
\newcommand{\lv}{\left \lvert}
\newcommand{\rv}{\right \rvert}

\newcommand{\mc}[1]{\mathcal{#1}}
\newcommand{\mbb}[1]{\mathbb{#1}}

\newcommand{\indi}[1]{\mathbbm{1}_{ \{#1\} }}

\newcommand{\reals}{\mathbb{R}} 

\usepackage{xspace}

\newcommand{\X}{\mathcal{X}}

\DeclareMathOperator*{\argmax}{arg\,max}

\newcommand{\tbf}[1]{\textbf{#1}}
\newcommand{\ttt}[1]{\texttt{#1}}

\usepackage{pgfplots}
\pgfplotsset{compat=newest}
\pgfplotsset{scaled y ticks=false}
\usepgfplotslibrary{groupplots}
\usepgfplotslibrary{dateplot}
\tikzstyle{every node}=[font=\small]
\pgfplotsset{
    yticklabel style={/pgf/number format/fixed},  
}

\makeatletter
\newcommand{\pushright}[1]{\ifmeasuring@#1\else\omit\hfill$\displaystyle#1$\fi\ignorespaces}
\newcommand{\pushleft}[1]{\ifmeasuring@#1\else\omit$\displaystyle#1$\hfill\fi\ignorespaces}
\makeatother

\expandafter\def\csname ver@etex.sty\endcsname{3000/12/31}

\usepackage{autonum}

\newcommand*{\wealth}{\mathcal{K}}

\newcommand*{\defined}{\coloneqq}

\newcommand{\payoff}{\mc{S}}

\newcommand{\Phat}{\widehat{P}}

\newcommand{\skeptic}{\texttt{Skeptic}\xspace} 
\newcommand{\forecaster}{\texttt{Forecaster}\xspace} 
\newcommand{\reality}{\texttt{Reality}\xspace} 

\newcommand{\iid}{i.i.d.\xspace}

\newcommand{\dks}{d_{\mathrm{KS}}}
\newcommand{\dkl}{d_{\mathrm{KL}}}
\newcommand{\dmmd}{d_{\mathrm{MMD}}} 
 
\newcommand{\dG}{d_{\mc{G}}} 
\newcommand{\gtilde}{\widetilde{g}} 
\newcommand{\Gtilde}{\widetilde{\mc{G}}}

\newcommand{\leqas}{\stackrel{a.s}{\leq }}
\newcommand{\geqas}{\stackrel{a.s}{\geq }}
\newcommand{\eqas}{\stackrel{a.s}{=}}
\newcommand{\eqdist}{\stackrel{d}{=}}

\newcommand{\prediction}{\mc{A}_{\text{pred}}}
\newcommand{\betting}{\mc{A}_{\text{bet}}}
\newcommand{\ons}{\mc{A}_{\text{ONS}}}
\newcommand{\oga}{\mc{A}_{\text{OGA}}}

\newcommand{\Aerm}{\mc{A}_{\text{ERM}}}

\newcommand{\convas}{\stackrel{a.s}{\longrightarrow}}

\newcommand{\probclasstwo}{ \mc{P}^{2}_{\epsilon, \mc{G}}}
\newcommand{\nullclass}{\mathcal{P}_{\text{null}}}
\newcommand{\altclass}{\mathcal{P}_{\text{alt}}}
\newcommand{\Tau}{\mathcal{T}}

\newcommand{\wealthU}{\wealth^{(u)}}
\newcommand{\modregret}{\mc{R}_n^{\varrho}}

\newcommand{\nzero}{n_0(\alpha, \Delta, \sigma)}
\newcommand{\sigmatilde}{\widetilde{\sigma}}

\title{Nonparametric Two-Sample Testing by Betting}
\date{}

\usepackage{authblk}
\author[1]{Shubhanshu Shekhar \thanks{shubhan2@andrew.cmu.edu}}
\author[1, 2]{Aaditya Ramdas \thanks{aramdas@stat.cmu.edu}} 
\affil[1]{Department of Statistics and Data Science, Carnegie Mellon University }
\affil[2]{Machine Learning Department, Carnegie Mellon University }

\begin{document}
\maketitle

\allowdisplaybreaks
\begin{abstract}
    We study the problem of designing consistent sequential  two-sample tests in a nonparametric setting. Guided by the principle of \emph{testing by betting}, we reframe this task  into that of selecting a sequence of payoff functions that maximize the wealth of a fictitious bettor, betting against the null in a repeated game. In this setting,  the relative increase in the bettor's wealth has a precise interpretation as the measure of evidence against the null, and thus our sequential test rejects the null when the wealth crosses an appropriate threshold. 
    We develop a general framework for setting up the betting game for two-sample testing, in which the payoffs are selected by a prediction strategy as data-driven predictable estimates of the \emph{witness function} associated with the variational representation of some statistical distance measures, such as integral probability metrics~(IPMs). We then formally relate the statistical properties of the test~(such as consistency, type-II error exponent and expected sample size) to the regret of the corresponding prediction strategy. 
    We construct a practical sequential two-sample test by instantiating our general strategy with the kernel-MMD metric, and demonstrate its ability  to adapt to the difficulty of the unknown alternative through theoretical and empirical results. 
    Our framework is versatile, and easily extends to other problems; we illustrate this by applying our approach to  construct consistent tests for the following problems: \textbf{(i)} time-varying two-sample testing with non-exchangeable observations, and \textbf{(ii)} an abstract class of ``invariant'' testing problems, including symmetry and independence testing. 
\end{abstract}

\tableofcontents

\section{Introduction}
\label{sec:introduction}
    \emph{Two-sample testing} is a fundamental problem in statistics, where  the goal is to check for the homogeneity of  samples drawn from two independent sources.   
    Prior works, with some exceptions discussed in~\Cref{subsec:related-work}, have mainly studied this problem in the \emph{batch} setting (also called the fixed sample size setting). 
    Since the sample size in batch tests is decided before collecting the observations, such tests run the risk of allocating too many observations on easier problem instances leading to wasted resources, or too few observations on harder problem instances resulting in inconclusive evidence against the null.
    To address these issues, we propose a general framework for designing consistent level-$\alpha$ sequential  nonparametric  tests for the two-sample testing problem, that automatically adapt the sample size to the unknown alternative.

\red{
    A large fraction of existing works in the sequential testing literature have focused on designing tests for simple null hypotheses, or for composite null hypotheses in parametric settings. Within these restricted (parametric) scenarios, several tests have been proposed that satisfy  strong optimality properties. These optimality properties, however, are heavily reliant on the model assumptions: for instance, see~\textcite[\S~5.2]{tartakovsky2014sequential} for an example with univariate Gaussians, where the optimality of Wald's sequential probability ratio test~(SPRT) breaks down.}

\green{    In contrast, the literature on \emph{nonparametric} sequential testing, the focus of this paper, is sparser. We work within the framework of ``sequential tests of power one'', as set out by \textcite{darling1968some}, which stop only on rejecting the (typically highly composite, multivariate and nonparametric) null. 
    Due to the generality of the composite nonparametric setting, the theoretical guarantees of these tests are slightly weaker: one typically hopes to ensure type-I error control uniformly over the null, asymptotic consistency (power one) under any alternative, and in some cases, an upper bound on expected stopping time under the alternative (hopefully implying both minimax optimality and instance optimality/adaptivity). 
    Our goal in this paper is to develop a  general framework for designing sequential two-sample tests satisfying the following properties:
    \begin{enumerate}
        \item[\textbf{(P1)}] They control  type-I error uniformly over the composite null,  and they have asymptotic power one under any alternative.
        \item[\textbf{(P2)}] Additionally,  under some conditions one can identify their  type-II error exponent and expected stopping time under the alternative (ideally minimax optimal and instance adaptive).  
        \item[\textbf{(P3)}] They are computationally efficient/feasible, and have good empirical performance. 
    \end{enumerate}
     We describe  our general strategy in~\Cref{sec:general-two-sample}, and use it  to instantiate powerful  two-sample tests satisfying \textbf{(P1)}-\textbf{(P3)}. This significantly improves upon the existing nonparametric sequential two-sample tests, none of which  satisfy all three properties. }
\green{
    Additionally, our ideas easily  generalize to (a) two-sample testing with  time-varying distributions~(\Cref{subsec:time-varying}), which is significant  since two-sample testing has long relied on exchangeability under the null  to run permutation tests, and (b) a class of abstract ``invariant'' testing problems~(e.g., testing for symmetry or independence) discussed in~\Cref{subsec:unified-approach}. 
    }
 
    Our design strategy is based on the general principle of \emph{testing by betting}, recently elucidated by~\textcite{shafer2021testing}. This principle extends the game-theoretic reformulation of the foundations of probability by~\textcite{shafer2019game} to hypothesis testing. 
    This approach establishes an equivalence between gathering evidence against the null,  and multiplying an initial wealth by a large  factor by repeatedly betting on the observations with payoff functions bought for their expected value under the null. \green{These connections between betting and probability have a long history, going back at least to the initial work of~\textcite{ville1939etude}, and we refer the reader to~\textcite[Appendix~F]{waudby2022estimating} for more details on the evolution of these ideas. }
    Overall, this approach transforms the task of (sequential) hypothesis testing into that of setting up  a betting game with  payoff functions that simultaneously ensure  \tbf{(i)} a \emph{fair} game under the null, and \tbf{(ii)} a large rate of growth of the wealth of the bettor under the alternative. 
    \blue{Note that, the main focus of~\textcite{shafer2021testing}  is on comparing the  the merits of the \emph{betting-score}~(that is, the ratio of final and initial wealth of the bettor) as a measure of evidence, against usual notions such as p-values. 
    Hence,~\textcite{shafer2021testing} does not propose constructive strategies for setting up the betting game and selecting the payoff functions, beyond some  basic tasks such as testing a simple null against a simple alternative. The two-sample testing problem considered in our paper, on the other hand, is nonparametric with composite null and composite alternative. Thus, to adapt the ideas of~\textcite{shafer2021testing} to our problem, we first need to develop a general framework for constructing appropriate betting games for composite testing problems, and design betting strategies, such that the resulting wealth process satisfies the two properties listed above. We develop such a framework in~\Cref{sec:general-two-sample} that uses statistical distance measures that admit a variational representation, and exploits certain symmetries to set up an appropriate betting game. 
    }
    
    \subsection{Testing by betting, test martingales and Ville's inequality}\label{sec:test-by-bet}
    
    To illustrate the above discussion, consider a hypothesis testing problem with a null hypothesis $H_0:P \in \nullclass$ and alternative $H_1:P \in \altclass$, and observations denoted by $Z_1, Z_2, \ldots$ lying in some space $\mc{Z}$, and drawn \iid according to $P$. To test the null $H_0$, a bettor may place repeated bets on the outcomes $\{Z_t: t \geq 1\}$ starting with an initial wealth $\wealth_0=1$. 
    A single round of betting~(say at time $t$) involves the following two steps. \textbf{(i)} First, the bettor selects a payoff function $\payoff_t: \mc{Z} \to [0, \infty)$, under the  restriction that it ensures a \emph{fair} bet if the null is true. Formally, this is imposed by requiring $\payoff_t$ to satisfy $\mbb{E}_P[\payoff_t(Z_t)|\mc{F}_{t-1}]=1$~\red{(or more generally,  $\mathbb{E}_{P}[\payoff_t(Z_t)|\mc{F}_{t-1}] \leq 1$)} for all $P \in \nullclass$, where $\mc{F}_{t-1} = \sigma(Z_1, \ldots, Z_{t-1})$. \textbf{(ii)} Then, the outcome $Z_t$ is revealed, and the bettor's wealth grows~(or possibly shrinks) by a factor of $\payoff_t(Z_t)$.
    Thus,  the bettor's wealth after $t$ rounds of betting is $\wealth_t =\wealth_0 \prod_{i=1}^t \payoff_i(Z_i)$. 
    
    The two key technical pieces that underpin the framework are \emph{test martingales}~\parencite{shafer2011test} and \emph{Ville's inequality}~\parencite{ville1939etude}.
    To elaborate, the restriction on the conditional expectation of the payoff functions implies that under the null, $\{\wealth_t: t\geq 0\}$ is a \emph{test martingale}, which is a nonnegative martingale with an initial value $1$. Due to this fact, $\wealth_t$ is unlikely to take large values for any $t \geq 1$. 
    On the other hand, when $H_1$ is true, the bettor's choice of payoff functions, $\{\payoff_t: t \geq 1\}$ should ensure that the wealth process~(or equivalently, the amount of evidence against the null) grows at a fast rate, ideally exponentially. 
    Such a  wealth process naturally leads to the following sequential test: \emph{reject the null if $\wealth_t \geq 1/\alpha$}, where $\alpha \in (0,1)$ is the desired confidence level. \red{Ville's maximal inequality~(recalled in~\Cref{fact:ville} in~\Cref{appendix:background})  ensures that this test controls the type-I error at level $\alpha$.}
    
    \blue{
    The discussion in the previous paragraph highlights some key design choices  that must be made to use these ideas for two-sample testing: in which function class should $\payoff_t$ lie; how to ensure $\mathbb{E}[\payoff_t|\mc{F}_{t-1}]=1$ uniformly over $\nullclass$; and how to ensure fast growth of $\wealth_t$ under the alternative?  
    }
    When testing simple hypotheses $H_0: Z_t \sim P$ and $H_1:Z_t \sim Q$ with $P$ and $Q$ known, an obvious choice of $\payoff_t$ is the likelihood ratio $dQ/dP$. Indeed, with this choice of payoff functions, we have $\mathbb{E}_P[\payoff_t|\mc{F}_{t-1}] =1$, meaning it is a fair bet under the null. Furthermore, it is easy to check that under $H_1$, the wealth process with this payoff grows exponentially, with an optimal (expected) growth rate of $\dkl(Q, P)$: the KL-divergence between $Q$ and $P$. 
    However, when dealing with cases where either one or both of $H_0$ and $H_1$ are composite and nonparametric~(as is the case with the two-sample testing problem considered in this paper), there is no obvious choice for the payoff functions. \emph{We propose a principled strategy of selecting  payoff functions, that result in powerful sequential tests for certain ``invariant'' testing problems like two-sample testing}.

    \subsection{Overview of Results}
        \label{subsec:contributions}
        Our goal in this paper is to design sequential level-$\alpha$ tests of \emph{power one}~\parencite{darling1968some} for the two-sample testing problem and its generalizations. 
        \begin{definition}[\ttt{sequential-test}]
        \label{def:sequential-test}
            A level-$\alpha$  sequential test can be represented by a random stopping time $\tau$ taking values in $\{1,2, \ldots \} \cup \{\infty\}$, and satisfying the condition
                $\mbb{P} \lp \tau < \infty \rp \leq \alpha$, 
            under the null $H_0$. Thus, $\tau$  denotes the random time at which the null hypothesis is rejected.
        \end{definition}
        
        In the next three sections, we focus mainly on the standard  two-sample testing problem, which can be stated as follows: given two independent streams of observations, $\{X_t: t \geq 1\}$ and $\{Y_t: t \geq 1\}$, drawn \iid from $P_X$ and $P_Y$ respectively, we wish to test the null $H_0: P_X = P_Y$ against the alternative $H_1: P_X \neq P_Y$. 
       
        \paragraph*{General two-sample testing strategy} In~\Cref{sec:general-two-sample}, we describe our framework for nonparametric two-sample testing by betting. We start by selecting a function class $\mc{G}$~(consisting of functions $g:\mc{X} \to [-1/2, 1/2]$), such that $\dG(P_X, P_Y) \defined \sup_{g \in \mc{G}} \mathbb{E}_{P_X}[g(X)] - \mathbb{E}_{P_Y}[g(Y)]>0$ whenever $P_X\neq P_Y$. An element $g^* \in \mc{G}$ achieving the supremum in the definition of~$\dG(P_X, P_Y)$, also called the \emph{witness function}, can be interpreted as the test function in $\mc{G}$ that best distinguishes $P_X$ from $P_Y$.  If $g^*$ were known, we could use it  to define the wealth process, $\wealth_t^* = \wealth_{t-1}^* \times (1 + \lambda^*(g^*(X_t) - g^*(Y_t)))$, where $\lambda^*$ is the log-optimal  fraction of wealth to bet in each round~(see~\eqref{eq:lambda-star} for formal definition). It can be verified that this process is a test martingale under the null, and grows exponentially to infinity under the alternative.  
        Using this `oracle' wealth process, we can then define the corresponding sequential test $\tau^* = \min \{n: \wealth^*_n \geq 1/\alpha\}$. However, since $g^*$ depends on the unknown $P_X$ and $P_Y$, we instead propose to use predictable empirical estimates, $\{g_t: t \geq 1\}$, of $g^*$ to design the practical sequential test $\tau$. Naturally, the quality of the estimates~(measured in terms of \emph{regret}) will govern the statistical properties of $\tau$, and we characterize this formally in~\Cref{theorem:main-result}. 
        
        \paragraph*{Instances of two-sample test} In~\Cref{subsec:kernel-mmd}, we instantiate the general approach described above for the case when $\dG$ is the kernel-MMD metric~($\dmmd$; a widely used distance metric in machine learning). We show that the resulting sequential kernel-MMD test has two key advantages over its batch counterpart: (i) adaptivity to unknown alternative, and (ii) lower computational costs~(as it does not use permutations). 
        
        Prior to introducing the kernel-MMD test, we study a related problem of testing the equality of means of two bounded random vectors in~\Cref{sec:bounded-mean-testing}. In this case, the distance $\dG$ is the $\ell_2$-distance between the means, and the resulting test can be interpreted as a kernel-MMD test with a linear kernel. The design and analysis of this test mirrors that of the kernel-MMD test, and hence it allows us to describe the key ideas involved, in a technically simpler setting. 
        
        \paragraph*{Generalizations} Our design approach easily extends to more general cases beyond standard two-sample testing. In~\Cref{subsec:time-varying}, we show how our strategy can be used to construct a consistent sequential test for two-sample testing with time-varying distributions. Then, in~\Cref{subsec:unified-approach}, we study an abstract problem of testing whether the unknown distribution is invariant to an operator $\Tau$ or not. This problem encompasses several tasks, such as two-sample testing, independence testing and symmetry testing. We describe the steps to construct a general level-$\alpha$ sequential test for this problem. Then, in~\Cref{theorem:abstract-power-result}, we characterize the detection boundary of a specific instance of this test~(based on a plug-in or empirical risk minimization method).

    \subsection{Related Work}
    \label{subsec:related-work}
        The area of sequential hypothesis testing was initiated by~\textcite{wald1945sequential}, who proposed and analyzed the Sequential Probability Ratio Test~(SPRT) for testing a simple null against a simple alternative. \textcite{wald1948optimum} established strong optimality properties of SPRT, and in particular, showed that the SPRT has the smallest expected sample size among all tests~(including fixed sample size)  that control the type-I and type-II errors below prescribed levels. Following~\textcite{wald1945sequential}, there has been a significant body of work on extending the SPRT to composite, but parametric family of hypotheses, and the reader is referred to~\textcite[Chapters~2~\&~4]{ghosh1991handbook}  for a detailed overview. 
        \red{Unlike SPRT and its generalizations, the focus of this paper is on \emph{power-one} tests~(described earlier, in \Cref{def:sequential-test}), pioneered by Robbins and collaborators. We discuss the details of some  relevant works in literature on designing power-one nonparametric two-sample tests below.}
        
        \textcite{darling1968some} considered several nonparametric one- and two-sample testing problems involving real-valued observations, and proposed sequential tests by deriving appropriate fixed sample size uniform deviation inequalities. 
        \textcite{howard2019sequential} proposed a sequential Kolmogorov-Smirnov~(KS) test by obtaining a tighter time-uniform deviation inequality for the empirical distribution functions.
        However, these sequential tests only work with real-valued observations (or more generally, observations in a totally ordered space) and cannot be applied in problems involving multivariate observations. The other tests discussed below address this issue. 
        
        \textcite{balsubramani2016sequential} derived a time-uniform empirical Bernstein inequality for random walks, and used it to design a sequential nonparametric two-sample test based on  the linear-time kernel-MMD test statistic.
        The original batch two-sample kernel-MMD test, proposed by~\textcite{gretton2012kernel}, uses a quadratic-time  empirical estimate of the squared MMD distance. 
        The reliance of the sequential test of~\textcite{balsubramani2016sequential} on the linear-time MMD statistic, while making the test computationally more efficient, also makes it less powerful than our proposed kernel-MMD test~(in~\Cref{subsec:kernel-mmd}) and the tests of~\textcite{lheritier2018sequential} and~\textcite{manole2023martingale} discussed below.
        
        \textcite{lheritier2018sequential} proposed a general approach to designing sequential nonparametric two-sample tests, by using sequentially learned probabilistic predictors of the labels indicating the population from which an observation was drawn. They identified sufficient conditions for the $\lambda$-consistency~(a weaker notion than usual consistency) of the resulting sequential test, and verified these conditions for a nearest-neighbor based predictor.
        In a subsequent paper~\parencite{lheritier2019low}, the authors proposed a new-predictor~(called \ttt{kd-switch}), which  results in  a sequential test satisfying the usual notion of consistency.  
        Compared to these tests, we construct sequential tests with stronger performance guarantees --- in addition to consistency,  we also characterize the type-II error exponent and the expected sample size of our tests.
        
        \textcite{manole2023martingale} propose a general technique for constructing confidence sequences for convex divergences between two probability distributions. Their approach relies on the key observation that the empirical divergence process is a reverse submartingale adapted to the exchangeable filtration. By instantiating the general confidence sequence for the special cases of the Kolmogorov-Smirnov metric~\parencite[\S~4.1]{manole2023martingale} and  kernel-MMD metric~\parencite[\S~4.2]{manole2023martingale}, the authors obtain consistent sequential nonparametric two-sample tests for both  univariate and multivariate   distributions. Unlike~\textcite{manole2023martingale}, our approach relies on constructing martingales~(instead of reverse submartingales) from statistical distances with a variational definition. Hence, our resulting sequential tests are expected to be less conservative than those of~\textcite{manole2023martingale} in rejecting the null. This intuition is verified in some numerical experiments in~\Cref{sec:numerical-experiments}.
  
\section{General Two-Sample Test}
\label{sec:general-two-sample}
    We describe our general strategy for constructing sequential two-sample tests based on the principle of testing by betting in~\Cref{subsec:construction-two-sample}. After describing the test, we then characterize its statistical properties in~\Cref{subsec:theoretical}. 
    We begin by  defining the two-sample testing problem.  
    \begin{definition}[Two-sample testing]
       \label{def:two-sample-testing}
           Given a stream of paired observations $\{(X_t, Y_t):t \geq 1\}$, drawn \iid according to $P_X \times P_Y$ on the observation space $\mc{X} \times \mc{X}$, our goal is to test the null, $H_0: P_X = P_Y$ against the alternative $H_1: P_X \neq P_Y$. 
    \end{definition}
    
    \begin{remark}
        \label{remark:two-sample-testing-definition}
        The above problem can be rewritten as testing $H_0: P_{X} \times P_{Y} \in \nullclass$  against $H_1: P_X \times P_Y \in \altclass$, where $\nullclass$ and $\altclass$ denote the null and alternative classes, defined as: 
    \begin{align}
        &\nullclass \defined \{P_X \times P_Y \in \mc{P}(\mc{X} \times \mc{X}): P_X, P_Y \in \mc{P}(\mc{X}), \text{ and } P_X =  P_Y \}, \quad \text{and} \label{eq:tst-nullclass}\\
        &\altclass \defined \{P_X \times P_Y \in \mc{P}(\mc{X} \times \mc{X}): P_X, P_Y \in \mc{P}(\mc{X}), \text{ and } P_X \neq  P_Y \}. \label{eq:tst-altclass}
    \end{align}
    In other words, the distributions in the null class are invariant to the action of the operator $\Tau:(\mc{X}\times \mc{X})\to (\mc{X}\times \mc{X})$ that takes elements $(x, y)  \in \mc{X}\times \mc{X}$ and flips their order; that is $\Tau (x, y) = (y,x)$. We will build upon this observation to extend our testing scheme to a significantly more general class of problems in~\Cref{subsec:unified-approach}.
    \end{remark}
    \begin{remark}
        \label{remark:one-sample}
        \blue{
            While we focus on two-sample testing in this section, the ideas we develop are also applicable to the related task of one-sample testing. In this problem, we are given a probability distribution $P_X$, and observations $\{Y_t: t \geq 1\}$ drawn \iid from an unknown distribution $P_Y$. The goal  is to test the null $H_0: P_Y=P_X$ against the alternative $H_1: P_Y \neq P_X$.  We discuss the details of using the techniques developed in this section for one-sample testing in~\Cref{appendix:one-sample}. 
        }
    \end{remark}

    \subsection{Construction of the Test}    
    \label{subsec:construction-two-sample}
    
    Our approach begins with choosing a distance measure~(often a metric)  on the space of probability distributions that admits a variational representation. Key examples of such distance measures are the integral probability metrics~or~IPMs~\parencite{muller1997integral} and $f$-divergence~\parencite{liese2006divergences} families. For concreteness, we focus on IPMs denoted by  
    $d_{\mc{G}}$,  defined as 
    \begin{align}
        d_{\mc{G}}(P, Q) = \sup_{g \in \mc{G}} |\mbb{E}_{P}[g(X)] - \mbb{E}_{Q}[g(Y)]|,  \label{eq:ipm}
    \end{align}
    where $\mc{G}$ denotes some class of functions from the observation space $\mc{X}$ to a bounded real-valued set, which we set to the interval $[-1/2,1/2]$ without loss of generality. For developing two-sample tests, we require that the function class $\mc{G}$ is rich enough to ensure $\dG(P,Q) >0$ for all $(P, Q) \in   \mc{P}(\mc{X}) \times \mc{P}(\mc{X})$, such that $P\neq Q$. \blue{In many problems, such as bounded-mean testing in~\Cref{sec:bounded-mean-testing}, we require this condition only on a subset of probability distributions. To formalize this, we introduce the notion of \emph{characteristic IPMs}, borrowing terminology from the literature on kernel mean embeddings of probability distributions~\parencite{fukumizu2004dimensionality, sriperumbudur2011universality}. 
    \begin{definition}
        \label{def:characteristic-IPM} 
        We say that the IPM $\dG$ associated with a function class $\mc{G}$ is \emph{characteristic} for a class of distribution pairs $\mc{P}_2 \subset \mc{P}(\mc{X}) \times \mc{P}(\mc{X})$, if for all $(P, Q) \in \mc{P}_2$, we have $\dG(P,Q) >0$. 
    \end{definition}
    \begin{remark}
        As mentioned above, in two-sample testing, we require $\dG$ to be characteristic for $\mc{P}_2 = \{(P,Q) \in \mc{P}(\mc{X}) \times \mc{P}(\mc{X}): P \neq Q\}$. \red{We will revisit this notion for other problems such as bounded-mean testing, symmetry testing and independence testing later in the paper.}
        For example, in symmetry testing with real-valued random variables, we will require $\dG$ to be characteristic for $\mc{P}_2$ which contains \red{all pairs of  non-symmetric distributions}. 
    \end{remark}
   } 
     Let $g^* \equiv g^*(P, Q, \mc{G})$  denote the \emph{witness function}, that is, the function in $\mc{G}$ at which the supremum in the definition of $\dG$ is achieved. If we interpret $\mc{G}$ as a class of test-functions used to distinguish two distributions $P$ and $Q$, then $g^*$ represents the element of $\mc{G}$  that provides maximum contrast between $P$ and $Q$.
     \red{If the supremum is achieved at more than one function, we can set $g^*$ to be any one of them. On the other hand, when the supremum in the definition of $\dG$ is not achieved, we can set $g^*$ to be any function that is $\epsilon$-suboptimal for an arbitrarily small $\epsilon>0$. In the cases we focus on, the supremum is achieved.}
    
    \paragraph{Oracle Test} In the two-sample testing problem, the distributions $P_X$ and $P_Y$ are unknown, and hence so is the witness function $g^* \in \mc{G}$ associated with them. Hypothetically, if  the witness function $g^*$ associated with $P_X$ and $P_Y$ were known, we could use it to define a test martingale~(introduced in~\Cref{sec:test-by-bet}), $\{\wealth_t^*: t \geq 0\}$ as follows: set $\wealth_0^*=1$ and 
    \begin{align}
        &\wealth_t^* = \wealth_{t-1}^*  \times \lp 1 + \lambda^* \big( g^*(X_t) - g^*(Y_t) \big) \rp, \quad \text{where}  \\
        &\lambda^* \in \argmax_{\lambda \in [-1,1]}\; \mbb{E}_{} [\log \lp 1 + \lambda (g^*(X)-g^*(Y))\rp], \quad \text{with } (X, Y) \sim P_X\times P_Y.  \label{eq:lambda-star}
    \end{align}
    The above expressions can be interpreted as the wealth process of a fictitious bettor, betting on the outcomes $\{(X_t, Y_t): t \geq 1\}$, with payoff function $\gtilde^*(x, y) = g^*(x) - g^*(y)$. \green{ The term $\lambda^*$ denotes the optimal constant bet that ensures the fastest growth rate of the wealth when $P_X \neq P_Y$, and it depends on $g^*$ in addition to $P_X$ and $P_Y$.
    The log-optimality criterion is motivated by older works such as~\textcite{kelly1956new, breiman1961optimal}, as well as the more recent papers of~\textcite{shafer2021testing,waudby2022estimating,grunwald2019safe}. As an additional benefit, this strategy ensures that the bettor never puts all its wealth on the line; that is, it never risks ending up with zero wealth.}
    The process $\{\wealth_t^*: t \geq 0\}$ is a test martingale when $P_X=P_Y$, and it grows to infinity at an exponential rate when $P_X\neq P_Y$. These two properties suggest the following level-$\alpha$  sequential test: 
    \begin{align}
        \label{eq:oracle-test-def}
        \tau^* \defined \min \{n \geq 1: \wealth_n^* \geq 1/\alpha \}, 
    \end{align}
    where the choice of threshold~(i.e., $1/\alpha$) is motivated by~Ville's inequality~\parencite{ville1939etude}. 
    It is easy to verify that $\tau^*$ is a consistent level-$\alpha$ test.

    \paragraph{Our proposed test} The oracle test discussed above is not practical, since it requires knowledge of the witness function $g^*$, and the optimal bet $\lambda^*$. Instead, to design a practical test, we propose to replace $g^*$ and $\lambda^*$ with  predictable estimates $\{g_t: t \geq 1\}$ and $\{\lambda_t: t \geq 1\}$. In other words, to construct a sequential test, we need a \emph{prediction strategy} $\prediction$ for selecting $\{g_t: t \geq 1\}$, and a \emph{betting strategy} $\betting$ for selecting $\{\lambda_t: t \geq 1\}$ in a data-driven manner. The strategy $\prediction$, is a collection of mappings $\{A_{t, \text{pred}}: t \geq 1\}$, where $A_{t, \text{pred}}$ maps $(X_1^{t-1}, Y_1^{t-1})$ to $g_t$, an element in $\mc{G}$. Similarly, $\betting$ is a collection of mappings, $\{A_{t, \text{bet}}:t \geq 1\}$, where $A_{t, \text{bet}}$ maps $(X_1^{t-1}, Y_1^{t-1})$ to $\lambda_t \in [-1,1]$. 
    
    For fixed prediction and betting strategies, we can define the wealth process $\{\wealth_t: t \geq 1\}$, and the corresponding sequential test $\tau$ as follows: 
    \begin{align}
        &\wealth_t = \wealth_{t-1} \times \lp 1 + \lambda_t \big(g_t(X_t) - g_t(Y_t) \big) \rp, \quad \text{with } \wealth_0 = 1, \quad \text{and} \label{eq:wealth-process-def} \\
        &\tau \defined \min \{t \geq 1: \wealth_t \geq 1/\alpha\} \label{eq:two-sample-test-def}. 
    \end{align}

    Thus, the task of designing a sequential test is equivalent to that of choosing prediction and betting strategies, $\prediction$ and $\betting$. The statistical properties of the resulting test will depend on how quickly $\prediction$ and $\betting$ approximate $g^*$ and $\lambda^*$. We will employ the Online Newton Step~(ONS) betting method of~\textcite{cutkosky2018black}, which we describe in~\Cref{def:ONS}, as the betting strategy in all subsequent tests. 
    The choice of the prediction strategy, however, will rely on the specific function class $\mc{G}$, as we show formally in~\Cref{theorem:main-result}. To simplify our presentation in this section, we will make the following assumption. 
    \begin{assumption}
        \label{assump:function-class-G} 
        For any $g$ that lies in $\mc{G}$, the function $-g$ also belongs to $\mc{G}$. This implies that we can write $d_{\mc{G}}(P,Q) = \sup_{g \in \mc{G}} \mathbb{E}_P[g(X)] - \mathbb{E}_Q[g(Y)]$, without the absolute value. 
    \end{assumption}
    
    \begin{remark}
        \label{remark:function-class-G} 
        We emphasize that~\Cref{assump:function-class-G} is made purely to simplify the presentation of the section. In cases where this condition is not satisfied, such as for Kolmogorov-Smirnov metric~($\dks$), we can construct the wealth process by hedging our bets over function classes $\mc{G}$ and  $-\mc{G} \defined \{-g: g \in \mc{G}\}$. That is, we can define the wealth process at time $t$ as $\wealth_t = \frac{1}{2} \lp \wealth_t^{+} + \wealth_t^{-} \rp$, where $\wealth_t^{+}$~(resp. $\wealth_t^{-}$) denotes the wealth process with payoffs from $\mc{G}$~(resp. $-\mc{G}$).  The details are in~\Cref{appendix:background-G}.
    \end{remark}

    The quality of a prediction strategy used for constructing the test $\tau$ is quantified in terms of its  regret, which we define below. 
    \begin{definition}[Regret of a prediction strategy] 
    \label{def:regret}
        The regret of a prediction strategy, $\prediction$,  on a sequence of observations $\{(X_t, Y_t): 1 \leq t \leq n\}$ is defined as 
        \begin{align}
            \label{eq:regret-def}
            \mc{R}_n \equiv \mc{R}_n\lp \prediction, \mc{G}, X_1^n, Y_1^n \rp \defined \lp \sup_{g \in \mc{G} }\sum_{t=1}^n  g(X_t) - g(Y_t) \rp  -  \sum_{t=1}^n g_t(X_t) - g_t(Y_t), 
        \end{align}
        where for any $t \geq 1$, the term $g_t$ is the $\mc{F}_{t-1}$-measurable element of $\mc{G}$ selected by the prediction strategy $\prediction$. Recall that $\mc{F}_{t-1} = \sigma \lp \{(X_i, Y_i): 1 \leq i \leq t-1\}\rp$. 
    \end{definition}

    Note that the above definition of regret places no probabilistic assumptions on the observations. A good strategy should ensure that the average regret, $\mc{R}_n(\prediction, \mc{G}, X_1^n, Y_1^n) / n$, converges to zero; at least on a per-sequence basis or,  better still,  uniformly over all observation sequences. In~\Cref{subsec:theoretical}, we show how the regret guarantees of a prediction scheme translates into statistical properties of the resulting sequential test. 
    
    We end this section by introducing the ONS betting strategy, $\ons$,  that we will use for selecting the bets $\{\lambda_t: t \geq 1\}$ in all the tests in this paper. 
    
    \begin{definition}[ONS betting strategy]
        \label{def:ONS}
        Let $\{v_t \in [-1,1]: t \geq 1\}$ denote a sequence of outcomes. 
        In the two-sample testing case, we have $v_t = g_t(X_t) - g_t(Y_t)$, where $\{g_t: t \geq 1\}$ are chosen using any prediction strategy $\prediction$.
        Initialize $\lambda_1 = 0$, and $a_0=1$. Then, for $t=1, 2, \ldots:$
        \begin{itemize}
            \item Observe $v_t \in [-1,1]$. 
            \item Set $z_t = {v_t / (1 + v_t \lambda_t)}$. 
            \item Update $a_t = a_{t-1} + z_t^2$. 
            \item Update $\lambda_{t+1}$ as follows: 
            \begin{align}
                \label{def:lambda-ons}
                \lambda_{t+1} = \min \left\{\frac{1}{2}, \max \left\{ - \frac{1}{2}, {\lambda_t + \frac{2}{2-\log 3} \frac{z_t}{a_t} }\right\} \right\}. 
            \end{align}
        \end{itemize}
    \end{definition}
    The ONS strategy is easy to implement, as it has a constant computational complexity at each time step $t$. Furthermore, if the observations $\{v_t: t \geq 1\}$ deviate from zero on average, the bets chosen by the ONS strategy ensure an exponential growth of the wealth process $\wealth_t$, as we discuss in~\Cref{appendix:background-ONS}.

    \subsection{Theoretical Analysis}
    \label{subsec:theoretical}
        In this section, we present the main theorem characterizing the statistical properties of the test introduced in the previous section. While the formal statements presented in~\Cref{theorem:main-result} require additional notation, the main takeaways can be summarized  as follows: 
        \begin{itemize}
            \item The test $\tau$, introduced in~\eqref{eq:two-sample-test-def}, controls the type-I error at the level $\alpha$ uniformly over the null class $\nullclass$ for any prediction strategy $\prediction$. 
            \item If $\prediction$ guarantees a limiting average regret smaller than $\dG(P_X, P_Y)$ almost surely when $P_X\neq P_Y$, the test $\tau$ is consistent. If, in addition, the limiting average payoff function is equal to $\dG(P_X, P_Y)$, then the process $\{\wealth_t: t \geq 1\}$ grows at an order optimal exponential rate. 
            \item If $\prediction$ guarantees zero limiting average regret uniformly over all sequence of observations with sufficiently high probability, then the test $\tau$ has a finite expected sample size. If this condition is guaranteed with probability one, then  we can also show that $\tau$ is \emph{exponentially consistent}; that is, its type-II error converges to zero exponentially. 
        \end{itemize}
        
        We now introduce the additional terms required to state our main result. First, we introduce 
        \begin{align}
            \Delta = \sup_{g \in \mc{G}} \mathbb{E}[g(X) - g(Y)]= \dG(P_X, P_Y), \quad \text{and} \quad 
            \sigma^2 = \sup_{g \in \mc{G}} \mathbb{V}\lp g(X) - g(Y) \rp, \label{eq:Delta-sigma}
        \end{align}
        where $\mathbb{V}(\cdot)$ denotes the variance operator. 
        For some non-increasing sequence $\{r_n \in [0,1]: n \geq 1\}$, define the  event $E_n = \{\mc{R}_n/n \leq r_n\}$ for all $n \geq 1$.  
        Introduce the term $n_0(\epsilon, \alpha)$ to denote the minimum number of observations  needed by $\prediction$ to make the average regret~(under the event $E_n$) \emph{sufficiently small}, 
        {
        \begin{align}
            \nzero \defined \min \lbr n\geq 1: r_n + \frac{\log(n/\alpha)}{n}  + \sigma \sqrt{ \frac{\log n/\alpha}{n} }< \Delta \rbr.  \label{eq:n0-def}
        \end{align}
        }
        Note that, if $\{r_n: n \geq 1\}$ converges to zero, then  $n_0(\alpha, \Delta, \sigma)$ is guaranteed to be finite for all $\Delta >0$.  
        \sloppy To simplify the notation, we suppress the $\prediction$, $\mc{G}$  and $(P_X,P_Y)$ dependence of $n_0$. 
        
        Finally, the term $\beta$, defined below, will be used to characterize the error exponent of our proposed test~(note that $\dkl$ below denotes the KL-divergence). 
            \begin{align}
               &\beta = \sup_{\epsilon>0}\; \inf_{P'\times Q' \in \mc{P}^2_{\epsilon, d}} \; \frac{ \dkl(P', P_X) + \dkl(Q', P_Y)}{2} ,   \label{eq:two-sample-exponent}\\
                \text{where} \quad &\probclasstwo = \{P'\times Q' \in \mc{P}(\X) \times \mc{P}(\X) : d_{\mc{G}}(Q', P') \leq \epsilon \}.  
            \end{align}
        To interpret this expression, first note that when $\epsilon=0$, the class $\probclasstwo$ corresponds to the set of null distributions for two-sample testing, $\nullclass$, introduced in~\eqref{eq:tst-nullclass}. 
        Then, for a fixed $\epsilon>0$,  the class $\probclasstwo$ is an $\epsilon$-expansion of $\nullclass$ in terms of the distance measure $\dG$. Hence, $\beta$ can be interpreted as the projection~(in terms of the KL-divergence, $\dkl$) of the pair $(P,Q)$ when $P\neq Q$ onto a vanishingly small expansion of the null set. In the special case, when $\dG$ metrizes weak convergence, $\beta$ admits the simpler representation: $\beta = \inf_{P' \in \mc{P}(\mc{X})}\; \frac{1}{2}\lp \dkl(P', P_X) + \dkl(P',P_Y)\rp$. 
        
        We can now state the main result of this section.        
        \begin{theorem}
            \label{theorem:main-result}
            Suppose $\dG$ is characteristic~(\Cref{def:characteristic-IPM}) for $\mc{P}_2 = \{(P,Q) \in \mc{P}(\mc{X}) \times \mc{P}(\mc{X}): P \neq Q\}$,  and~\Cref{assump:function-class-G} holds. Consider  observations $\{(X_t, Y_t): t \geq 1\}$ drawn \iid according to $P_X \times P_Y$. 
            Let $\tau \equiv \tau(\prediction, \ons)$ denote a sequential test with prediction strategy $\prediction$, and betting strategy $\ons$ introduced in~\Cref{def:ONS}. Then, the following statements are true: 
            \begin{itemize}
                \item For any prediction strategy,  the type-I error rate is controlled at the specified level $\alpha$ uniformly over the  class of null distributions~\eqref{eq:tst-nullclass}. That is, 
                \begin{align}
                    \label{eq:main-test-type-I}
                    \sup_{P_X \times P_Y \in \nullclass} \; \mbb{P}_{P_X \times P_Y}\lp \tau < \infty \rp \leq \alpha. 
                \end{align}
                \item Suppose $P_X \neq P_Y$, and  the per-sequence average regret of $\prediction$ satisfies 
                \begin{align}
                    \label{eq:no-regret}
                    \limsup_{n \to \infty} \frac{\mc{R}_n\lp \prediction, \mc{G}, X_1^n, Y_1^n\rp}{n} < \dG(P_X, P_Y) \quad \text{almost surely}. 
                \end{align}
                Then, the sequential test $\tau$ has power one under the alternative:
                \begin{align}
                    \mbb{P}_{P_X\times P_Y} \lp \tau < \infty \rp = 1, \quad \text{for each } P_X \times P_Y \in \altclass. \label{eq:main-test-consistency}
                \end{align}
                \item Suppose $P_X \neq P_Y$, and the strategy $\prediction$ for choosing $\{g_t: t \geq 1\}$ satisfies 
                \begin{align}
                    \lim_{n \to \infty} \frac{1}{n} \sum_{t=1}^n g_t(X_t) - g_t(Y_t) = \mathbb{E}_{P_X}[g^*(X)] - \mathbb{E}_{P_Y}[g^*(Y)] = \Delta. 
                    \label{eq:limiting-average-payoff}
                \end{align}
                Then, the test $\tau$ is consistent, and furthermore the process process $\wealth_n$ grows to infinity at an exponential rate:
                \begin{align}
                    \liminf_{n \to \infty} \frac{1}{n} \log \wealth_n \geq \Delta \lp \frac{ \Delta }{\mathbb{E}[(g^*(X) - g^*(Y))^2]} \bigwedge 1 \rp. 
                    \label{eq:wealth-growth-rate-1}
                \end{align}
                \item Suppose  that $P_X \neq P_Y$,  and  there exists a sequence $\{r_n: n \geq 1\}$ such that  $r_n \to 0$ and $\sum_{n \geq 1} \mathbb{P}(E_n^c) < \infty$, where $E_n = \{\mc{R}_n/n \leq r_n\}$. Then,   we have the following upper bound on the expected stopping time: 
                \begin{align}
                    \mathbb{E}[\tau]  = \mc{O}\lp n_0(\alpha, \Delta, \sigma) +\sum_{n = 1}^{\infty} \mathbb{P}(E_n^c) \rp.   \label{eq:main-test-stopping-time}
                \end{align}
                \item Suppose $P_X \neq P_Y$, and there exists a sequence $\{r_n: n \geq 1\}$ with $r_n \to 0$, and $\mathbb{P}(E_n^c)=0$ for all but finitely many $n \geq 1$. Then, we have the following under $H_1$ with $P_X \times P_Y \in \altclass$:
                \begin{align}
                    &\lim_{n \to \infty} -\frac{1}{2n} \log \lp \mbb{P}_{P_X\times P_Y} \lp \tau > n \rp \rp  \geq \beta.   \label{eq:main-test-exponent} 
                \end{align}
                Recall that the terms $n_0(\alpha, \Delta, \sigma)$, and $\beta$ were defined in~\eqref{eq:n0-def} and~\eqref{eq:two-sample-exponent} respectively. 
            \end{itemize}
        \end{theorem}
        
        \red{We will see later, that in many practical tests, both $\nzero$ and $\sum_{n \geq 1} \mathbb{P}(E_n^c) = \mc{O}(1)$  are of the order $\mc{O}(\sigma^2\log(1/\alpha \Delta)/ \Delta^2)$, and thus~\eqref{eq:main-test-stopping-time}  implies that $\mathbb{E}[\tau]$ is also $\mc{O}\lp \sigma^2 \log(1/\alpha\Delta)/ \Delta^2 \rp$. }
        A detailed proof of~\Cref{theorem:main-result} is in~\Cref{appendix:proof-main-theorem}. 
        \begin{remark}
            \label{remark:e-processes}
            As mentioned before,~\eqref{eq:main-test-type-I} is a consequence  of the fact that  $\{\wealth_t: t \geq 1\}$ is a composite test martingale for the class of distributions $\nullclass$. 
            In fact, the  process $\{\wealth_t: t \geq 0\}$ associated with any prediction strategy $\prediction$  can be shown to satisfy the property $\sup_{ P_X\times P_Y \in \nullclass } \sup_{\tau'}  \mbb{E}_P \lb \wealth_{\tau'} \rb \leq 1$, 
            where $\tau'$ is any stopping time adapted to $\{\mc{F}_t: t \geq 0\}$.  
            This inequality is the defining property of e-processes studied in recent works such as~\parencite{ramdas2021testing, grunwald2019safe}. 
            Hence, our general approach for designing the wealth process $\{\wealth_t:  t \geq 1\}$ described in this section,  can also be thought of as a method of constructing nontrivial e-processes for $\nullclass$ that grow to infinity  under the alternative if the prediction strategy suffers vanishing average regret.
        \end{remark}
        \begin{remark}
            \label{remark:wealth-growth-rate} 
            It is easy to check that the condition~\eqref{eq:limiting-average-payoff} is satisfied by the plug-in strategy for any function class $\mc{G}$ that is uniformly learnable. Furthermore, following the same arguments as~\textcite[Proposition~1]{podkopaev2022sequential}, we can also conclude that the growth rate of the oracle process, $\{\wealth_t^*: t \geq 1\}$, is of the order $\Theta \lp \Delta \lp  \frac{\Delta}{\mathbb{E}[(g^*(X)-g^*(Y))^2]} \wedge 1\rp\rp$. This in turn implies the order optimality of the rate obtained in~\eqref{eq:wealth-growth-rate-1} under the assumption~\eqref{eq:limiting-average-payoff}. 
            Additionally, when the variance $\mathbb{V}(g^*(X)-g^*(Y))$ is large~(w.r.t. $\Delta^2$), the growth rate is $\asymp \Delta^2$, while for small variance, the rate increases to $\asymp \Delta$; thus displaying empirical Bernstein-type variance adaptivity. 
        \end{remark} 
        \begin{remark}
        \label{remark:non-independent-streams} 
            We introduced the two-sample testing problem in~\Cref{def:two-sample-testing} with the assumption that the two-streams of \iid observations, $\{X_t: t \geq 1\}$ and $\{Y_t: t \geq 1\}$, are independent of each other, following the standard formulation in literature. However, we note that the results stated above in~\Cref{theorem:main-result} are valid under a much weaker assumption that the stream $\{(X_t, Y_t): t \geq 1\}$ consists of independent pairs of observations satisfying $(X_t, Y_t) \eqdist (Y_t, X_t)$ under the null, and $(X_t, Y_t) \not \eqdist (Y_t, X_t)$ under the alternative. This allows our tests to be applicable in interesting cases where standard batch methods, such as permutation tests, fail. An example is when $X_t$ and $Y_t$ are obtained from a common parent variable $Z_t$; that is, $X_t = h(Z_t,  \widetilde{X}_t)$ and $Y_t = h(Z_t, \widetilde{Y}_t)$,  for some function $h$, and  with independent $\widetilde{X}_t$ and  $\widetilde{Y}_t$. 
        \end{remark}

    \subsection{Extension to unbounded \texorpdfstring{$\mc{G}$}{G}}
    \label{subsec:unbounded-functions}
        We now discuss how to design sequential tests when the class of test functions, $\mc{G}$, contains possibly unbounded functions. \red{This allows us to significantly expand the class of test functions~($\mc{G}$) that can be employed within the framework of~\Cref{subsec:construction-two-sample}. For instance, following~\textcite{kim2021classification}, we could use as $g_t$, any classifier trained on the first $t-1$ pairs of observations to distinguish $P_X$ and $P_Y$.  
        }

        Recall that our general approach uses the fact that $\mathbb{E}[g(X)-g(Y)]=0$ for all $g \in \mc{G}$ under the null. However, note that $g(X)-g(Y)$ satisfies the stronger condition: $g(X)-g(Y) \stackrel{d}{=} g(Y)-g(X)$; that is, it is symmetric for all $g \in \mc{G}$ under the null. This suggests that we can use ideas from symmetry testing, such as those developed by~\textcite{ramdas2020admissible}, to design appropriate test martingales in the case of unbounded $\mc{G}$.
        
        The only change we need to make is to apply an odd sigmoid function to the difference $g_t(X_t) - g_t(Y_t)$, for selecting the bets and updating the wealth. We first formally state the properties required of the sigmoid function. 
        
        \begin{definition}[Sigmoid function]
        \label{def:sigmoid-func} Let $\varrho: \mathbb{R} \to [-1,1]$ denote a sigmoid function satisfying: 
        \begin{enumerate}
            \item $\varrho$ is an odd function, that is $\varrho(x) = -\varrho(-x)$ for all $x \in \reals$. 
            \item $\lim_{x \to \infty} \varrho(x) = \lim_{x \to -\infty} - \varrho(x) = 1$. 
        \end{enumerate}
        An example  of $\varrho$ suggested by \textcite{ramdas2020admissible} is $\varrho(x) = \text{tanh}(x) = \frac{e^{x} - e^{-x}}{e^{x} + e^{-x}}$.
        \end{definition}
        
        As before, let $\prediction$ denote any prediction strategy for selecting the functions $\{g_t: t \geq 1\}$. Introduce the terms $v_t \defined \varrho\big( g_t(X_t) - g_t(Y_t) \big)$, and  let  $\{\lambda_t: t \geq 1\}$ denote the bets selected by the ONS strategy, $\ons$, applied to $\{ v_t: t \geq 1\}$. We can define the following process: 
        \begin{align}
            \wealthU_0 = 1, \quad \text{and} \quad \wealthU_{t} = \wealthU_{t-1} \times \lp 1 + \lambda_t v_t \rp \; \text{for } t \geq 1, \quad  \text{where } v_t = \varrho\lp g_t(X_t) - g_t(Y_t) \rp.  \label{eq:wealth-unbounded}
        \end{align}
        Since $\varrho$ is an odd function, $\mathbb{E}[v_t|\mc{F}_{t-1}]=0$, which implies that $\{\wealthU_t: t \geq 1\}$ is a test martingale under the null. This suggests the following level-$\alpha$ test: $\tau^{\varrho} = \min \{n\geq 1: \wealthU_n \geq 1/\alpha\}$.
        
        As in the case of~\Cref{theorem:main-result}, the consistency, exponential consistency and the expected stopping time bound for the test $\tau^{\varrho}$ based on the modified process $\{\wealthU_t: t \geq 1\}$ can be characterized in terms of the modified regret, $\modregret$, defined as 
        \begin{align}
            \label{eq:modified-regret} 
            \modregret \equiv \modregret\lp \prediction, \mc{G}, X_1^n, Y_1^n \rp \defined \lp \sup_{g \in \mc{G} }\sum_{t=1}^n \varrho\big(  g(X_t) - g(Y_t) \big)\rp  -  \sum_{t=1}^n \varrho \big(g_t(X_t) - g_t(Y_t) \big). 
        \end{align}
        
        \red{ Since the function $\varrho$ is not convex, developing prediction schemes that have theoretical guarantees on the modified  regret $\modregret$ may be nontrivial. However, some experimental results in~\Cref{sec:numerical-experiments} indicate that tests based on the prediction schemes that minimize the usual regret, $\mc{R}_n$, still perform very well in practice. 
        }

\section{Bounded Mean Testing}
\label{sec:bounded-mean-testing} 
    To illustrate the steps involved in instantiating the general sequential test of~\Cref{sec:general-two-sample}, we begin by considering a related, but conceptually simpler, task of testing the equality of the means of two distributions.
    In this problem, we are given two streams of observations $\{X_t: t \geq 1\}$ and $\{Y_t: t \geq 1\}$, taking values in the unit ball in the $m$-dimensional Euclidean space,  $\X = \lbr x \in \mbb{R}^m: \|x\|_{2} \leq 1 \rbr$. We assume that all the $X_t$~(resp. $Y_t$) are drawn \iid from a distribution $P_X$~(resp. $P_Y$), and our goal is to test 
    \begin{align}
        H_0: \mbb{E}_{P_X}[X] = \mbb{E}_{P_Y}[Y], \quad \text{versus} \quad H_1: \mbb{E}_{P_X}[X] \neq \mbb{E}_{P_Y}[Y]. \label{eq:bounded-mean-def}
    \end{align}
    
    \begin{remark}
    \label{remark:bounded-mean-1}
    \blue{
        The two-sample testing problem is equivalent to checking for the equality of the characteristic functions of the two random variables --- a much more stringent requirement than~\eqref{eq:bounded-mean-def}. Nevertheless, bounded mean testing  is a nontrivial task, with  highly composite null and alternative classes. Other than being bounded, the distributions do not have any other restrictions; they can be discrete, continuous or mixed, and they do not need to have a common dominating measure. This prevents several classical techniques from being applied to this problem. This section can be seen as extending the univariate one-sample methods of~\textcite{waudby2022estimating} to the multivariate and two-sample setting.}
    \end{remark}

    For the rest of this section, we will first define a sequential test for~\eqref{eq:bounded-mean-def} using the general framework of~\Cref{sec:general-two-sample}, and then theoretically characterize its performance in~\Cref{prop:bounded-mean}. 
    
    \paragraph{Instantiating our test} To instantiate the sequential test from~\Cref{sec:general-two-sample}, we need to specify a function class $\mc{G}$ and a prediction strategy $\prediction$~(recall that the betting strategy is set to $\ons$ of~\Cref{def:ONS}). 
    
        We set $\mc{G} = \{g_u(\cdot) = \langle u, \cdot \rangle : u \in \mbb{R}^m, \, \|u\|_2 \leq 1/2 \}$, where $\langle x, y \rangle = x^Ty$ denotes the usual inner product in $\mbb{R}^m$. For the rest of this section, we will use $\mc{U}$ to denote the ball $\{u \in \mathbb{R}^m: \|u\|_2 \leq 1/2\}$.  The corresponding distance between two distributions $P$ and $Q$ taking values in $\X$ then becomes
        \begin{align}
            d_{\mc{G}} \lp P, Q\rp = \sup_{u: \|u\|_2 \leq 1/2} \; \langle \mathbb{E}_P[X] - \mathbb{E}_Q[Y], u \rangle = \frac{1}{2} \| \mathbb{E}_P[X] - \mathbb{E}_Q[Y] \|_2. 
        \end{align}
        Thus the distance $d_{\mc{G}}$ reduces to the Euclidean norm between the means of the two distributions~(scaled by $1/2$); a quantity which is zero only when the two means are equal. Hence, $\dG$ is \emph{characteristic} for the alternative class of distributions in the sense of~\Cref{def:characteristic-IPM}. 
        
        Next, we need to select a prediction strategy for selecting $\{g_t: t \geq 1\}$ from $\mc{G}$; or equivalently for selecting $\{u_t \in \mc{U}: t \geq 1\}$. The regret for this prediction problem is 
        \begin{align}
            \frac{\mc{R}_n}{n} &= \frac{1}{n} \lp  \max_{u: \|u\|_2 \leq 1/2} \sum_{t=1}^n \langle u, \, X_t - Y_t \rangle  - \sum_{t=1}^n \langle u_t, \, X_t - Y_t \rangle  \rp  \\
            & = \max_{u: \|u\|_2 \leq 1/2} \sum_{t=1}^n \langle u - u_t, \, \bar{X}_n - \bar{Y}_n \rangle.  
        \end{align}
        Since this regret corresponds to an online prediction problem with linear losses, we use an adaptive version of the  \emph{projected online gradient ascent}~(OGA) method~(see~\Cref{appendix:background-oga} for further details) as the prediction strategy. To describe this strategy, we introduce the notation $M_t = \sum_{i=1}^t \|X_i-Y_i\|_2^2$, and set $u_1 = \boldsymbol{0} \in \mc{U}$ and  $\eta_0=0$. Then, for any $t \geq 1$, we define the next payoff function, $g_{t+1}(\cdot) = \langle u_{t+1}, \cdot\rangle$ as follows: 
        \begin{align}
            \label{eq:oga-1}
            u_{t+1} =  \Pi_{\mc{U}} \lp u_{t} + \frac{1}{\sqrt{M_t}}\lp X_{t} - Y_{t} \rp \rp,
        \end{align}
        where $\Pi_{\mc{U}}(x)$ denotes the projection of $x \in \mathbb{R}^m$ onto the set $\mc{U} = \{u \in \mathbb{R}^m: \|u\|_2 \leq 1/2\}$. 
        
        We now summarize all the steps of our sequential test for bounded mean testing below. 
        \begin{definition}[Sequential Test for Bounded Means Testing]
        \label{def:bounded-means-test}
            Set $\wealth_0 = 1$, $\lambda_1=0$, and $u_1=\boldsymbol{0} \in \mathbb{R}^m$. For $t=1, 2, \ldots$:
            \begin{itemize}
                \item Observe $X_t, Y_t$. 
                \item Update the wealth: $\wealth_{t} = \wealth_{t-1} \times \lp 1 + \lambda_t \langle u_t, X_t - Y_t \rangle \rp$. Reject the null if $\wealth_t \geq 1/\alpha$. 
                \item Obtain $u_{t+1}$ from the OGA prediction strategy~\eqref{eq:oga-1}. 
                \item Obtain $\lambda_{t+1}$ from the ONS betting strategy of~\Cref{def:ONS}. 
            \end{itemize}
        \end{definition}

    \paragraph{Computational Complexity} The computational complexity of the test described in~\Cref{def:bounded-means-test} is $\mc{O}\lp \tau m\rp$, where $\tau$ denotes the (random) time at which the test rejects the null and $m$ is the dimension of $\mc{X}$. This follows  from the fact that for any $t$,  the update of $u_t$ by OGA strategy is  an $\mc{O}(m)$ operation, the update of $\wealth_t$ is also an $\mc{O}(m)$ operation, and the update of $\lambda_{t+1}$ via ONS strategy is  an $\mc{O}(1)$ operation. Hence, the complexity of one step of the test is $\mc{O}(m)$. Since the test rejects the null after $\tau$ observations, the overall complexity of the test is $\mc{O}(\tau m)$, as claimed.

    \paragraph{Statistical Properties} We can establish the statistical properties of the above test by  specializing the general results of~\Cref{theorem:main-result}. In particular,  we use the fact that the average regret of the OGA strategy in this problem converges to zero uniformly~(see~\Cref{appendix:background-oga}), which implies the following results as a consequence of~\Cref{theorem:main-result}~(details in~\Cref{proof:bounded-mean}). 
    \begin{proposition}
    \label{prop:bounded-mean}
        Introduce the terms $\Delta = \dG(P_X, P_Y) =  (1/2)\|\mathbb{E}_{P_X}[X] - \mathbb{E}_{P_Y}[Y] \|_2$,  and  $\sigma^2 = \sup_{u:\|u\|\leq 1/2} \mathbb{V}\lp \langle u, X-Y\rangle \rp = \mc{O}\lp \mathbb{E}[\|X-Y\|^2]\rp$. Then, for the sequential test described in~\Cref{def:bounded-means-test}, we have the following: 
        \begin{align}
            &\text{Under } H_0: \; \mathbb{P}(\tau < \infty) \leq \alpha. \\
            &\text{Under } H_1: \; \mathbb{P} \lp \tau <\infty \rp = 1, \quad \text{and} \quad \mathbb{E}[\tau] = \mc{O} \lp \frac{ \sigma^2 \log(1/\Delta \alpha) }{\Delta^2}  {+ \frac{\log(1/\alpha\Delta)}{\Delta}}\rp. \label{eq:tau-bounded-mean-1}
        \end{align}
    \end{proposition}
    Thus, following the general framework introduced in the previous section, we have constructed a level-$\alpha$ sequential test that is consistent against any fixed alternative. An interesting aspect of the bound in~\eqref{eq:tau-bounded-mean-1} is the presence of the second moment term, $\sigma^2$, in the numerator. In the worst case, this term can be equal to $1/\sqrt{2}$. However, for problems instances with $\sigma^2$ is small, our test has the ability to exploit this additional structure, and stop earlier.

    \paragraph{Optimality} We end this section by showing that there exist problem instances on which the expected stopping time of our test, derived in~\Cref{prop:bounded-mean}, cannot be improved (modulo polylogarithmic terms). In particular, we do this by restricting our attention to a simple class of problems, where $P_X$ and $P_Y$ are distributions with finite support. 
    \begin{proposition}
        \label{prop:lower-bound-bounded-mean}  Let $\tau'$ denote any level-$\alpha$, power-one sequential test for the problem defined in~\eqref{eq:bounded-mean-def}, such that the expected value of $\tau'$ is finite under $H_1$. Then, there exist distributions $P_X$ and $P_Y$  supported  on a finite subset of $[0,1]$, with  $\|\mathbb{E}_{P_X}[X] - \mathbb{E}_{P_Y}[Y]\|_2 = \Delta$ and $\mathbb{E}[(X-Y)^2]= \Omega(\sigma^2)$, such that we have $\mathbb{E}[\tau'] = \Omega \lp \frac{\sigma^2 \log \lp1/\alpha \rp}{\Delta^2} + {\frac{\log(1/\alpha)}{\Delta}} \rp$. 
    \end{proposition}
    The proof of this result is in~\Cref{proof:lower-bound}, and it proceeds by first obtaining a general information-theoretic lower bound on any $\mathbb{E}[\tau']$ in terms of $d_{KL}(P_X, P_Y)$, and then constructing distributions for which $d_{KL}(P_X, P_Y)$ can be approximated in terms of $\Delta$ and $\sigma$ to get the required expression. 
    
\section{Sequential Two-Sample Kernel-MMD Test}
    \label{subsec:kernel-mmd}
            We now instantiate the general strategy introduced in~\Cref{sec:general-two-sample} for the two-sample testing problem, using the kernel maximum mean discrepancy (MMD) metric. \red{The batch two-sample test based on this metric~\parencite{gretton2012kernel} is widely used  in practice, as it can be applied to observations in arbitrary spaces as long as we can define positive-definite kernels  on them. However, the existing sequential versions of this test, such as those proposed by~\textcite{balsubramani2016sequential, manole2023martingale}, often have poor empirical performance. In this section, we use our framework to design the first sequential kernel-MMD test with strong theoretical guarantees along with good empirical performance.} 
            
            Let $\mc{X}$ denote the observation space, which for simplicity, we set to $\mathbb{R}^m$ for some $m \geq 1$, and let $K: \mc{X} \times \mc{X} \to \mathbb{R}$ be a positive definite kernel on $\mc{X}$. We assume that $K$ is uniformly bounded, that is,  $\sup_{x, x' \in \mc{X}}\sqrt{K(x,x')} \leq 1$,  and  let $\mc{H}_K$ denote the reproducing kernel Hilbert space~(RKHS) associated with $K$.
            Since the unit ball in the RKHS $\mc{H}_K$ satisfies~\Cref{assump:function-class-G}, the associated IPM, called the kernel-MMD metric,  is defined as follows: 
            \begin{align}
                \label{eq:kernel-mmd}
                \dmmd(P,Q)= \sup_{\|g\|_K \leq 1} \mbb{E}_{P}[g(X)] - \mbb{E}_{Q}[g(Y)],
            \end{align}
            where $\|g\|_K$ denotes the RKHS norm of the function $g$. 
            The mean map of a distribution $P$ is a function in the RKHS  given by $\mu_P \defined \mbb{E}_P[K(X, \cdot)]$.
            When $P\neq Q$, the  ``witness'' function $h^*$ that achieves the supremum in $\dmmd(P,Q)$ (i.e. witnesses the difference between $P,Q$) is simply given by $h^* := \frac{\mu_P - \mu_Q}{\|\mu_P - \mu_Q \|_K}$,
            meaning that $\dmmd(P,Q) = \mbb{E}_{P}[h^*(X)] - \mbb{E}_{Q}[h^*(Y)]$. 
            
            \paragraph{Instantiating the test}
                The above discussion suggests the choice of $\mc{G} = \{g \in \mc{H}_K: \|g\|\leq 1/2\}$, and $g^* = \frac{1}{2}h^*$:
                \begin{align}
                    \label{eq:mmd_oracle_function} 
                    \gtilde^*(x,y) &=  g^*(x) - g^*(y)  = \lp \langle g^*, K(x, \cdot) - K(y, \cdot) \rangle \rp. \label{eq:kernel_max}
                \end{align}
                Note that the scaling $h^*$ by $1/2$ in the definition of $g^*$ ensures that $\gtilde^*$ takes values in $[-1,1]$. 
                
                Having selected $\mc{G}$, the final step in instantiating the sequential test is choosing an appropriate prediction strategy. The regret of the prediction game after $n$ observations for a prediction strategy $\prediction$ playing $\{g_t: t \geq 1\}$ is  
                \begin{align}
                    \frac{1}{n} \mc{R}_n \equiv \frac{1}{n} \mc{R}_n \lp \prediction, \mc{G}, X_1^n, Y_1^n \rp =  \max_{g \in \mc{G}}  \frac{1}{n} \sum_{t=1}^n \langle g - g_t, K(X_t, \cdot) - K(Y_t, \cdot) \rangle.  
                \end{align}
                
                A natural choice is the plug-in or the empirical risk minimization~(ERM) strategy, that simply selects $g_t =    \argmax_{g \in \mc{G}} \langle g, \mu_{\Phat_{X, t-1}} - \mu_{\Phat_{Y,t-1}} \rangle$~(see~\Cref{remark:erm-special-case}).  We can check that this choice results in a consistent sequential test. To get the exponent and bound on the expected stopping time under the alternative, however, we need to use an adaptive version of the online gradient ascent (OGA) strategy, that proceeds as follows, with $M_t \defined \sum_{i=1}^{t} \|g_i(X_i, \cdot) - g_i(Y_i, \cdot)\|_K^2$: 
                \begin{align}
                    \label{eq:oga-2}
                    g_1=\boldsymbol{0}, \quad \text{and } g_{t+1} = \Pi_{\mc{G}} \lp g_t + \frac{1}{2\sqrt{M_t}} \lp K(X_t, \cdot) - K(Y_t, \cdot) \rp\rp \; \text{for } t \geq 1.
                \end{align}
                Recall that $\Pi_{\mc{G}}$ denotes the projection operator (in terms of the RKHS norm $\|\cdot\|_K$) onto the function class $\mc{G}$, which acts as follows: $\Pi_{\mc{G}}(h) = h / (2\|h\|_K)$. 
            We  now formally describe our sequential kernel MMD test.
            
            \begin{definition}[\texttt{Sequential Kernel MMD Test}]
            \label{def:kernel_mmd_test} 
                Set $\wealth_0 = 1$, $\lambda_1=0$, and $g_1=\boldsymbol{0} \in \mc{H}_K$. For $t=1, 2, \ldots$:
                \begin{itemize}
                    \item Observe $X_t, Y_t$. 
                    \item Update the wealth: $\wealth_{t} = \wealth_{t-1} \times \lp 1 + \lambda_t \langle g_t, K(X_t,\cdot) - K(Y_t, \cdot) \rangle \rp$. 
                    \item Reject the null if $\wealth_t \geq 1/\alpha$. 
                    \item Update $g_{t+1}$ using the OGA prediction strategy~\eqref{eq:oga-2}. 
                    \item Update $\lambda_{t+1}$ using the ONS betting strategy~\Cref{def:ONS}. 
                \end{itemize}
           \end{definition}

           \paragraph{Computational Complexity} 
            The computational complexity of the test described in~\Cref{def:kernel_mmd_test} is $\mc{O}\lp \tau^2\rp$, where $\tau$ denotes the (random) time at which the test rejects the null. This differs from the test for equality of means in~\Cref{sec:bounded-mean-testing}, which has a linear dependence on $\tau$. This is because the computation of the inner product $\langle g_t, K(X_t, \cdot)-K(Y_t, \cdot)\rangle$ in the kernel-MMD test is a $\mc{O}(t)$ time operation as $g_t$ is a sum of $t-1$ terms, unlike the corresponding inner product in~\Cref{def:bounded-means-test}, which is a $\mc{O}(m)$ operator at all $t$. Adding this up over all values of $t$ from $1$ to $\tau$ implies the overall quadratic complexity.
            
           \paragraph{Statistical Properties} 
            As in the case of bounded-mean testing, we can establish the statistical properties of the kernel-MMD test by  specializing~\Cref{theorem:main-result}. Based on the regret bound for OGA strategy stated in~\Cref{appendix:background-oga}, we note that the regret satisfies $\mc{R}_n = \mc{O}(\sqrt{M_n})$. 
            Plugging this into~\Cref{theorem:main-result}, we get the following with $\Delta = \dmmd(P_X, P_Y)$ and $\sigma^2 = \mc{O}\lp \mathbb{E}[\|k(X,\cdot) - k(Y, \cdot)\|_K^2] \rp$:
            \begin{proposition}
            \label{prop:kernel-mmd}
                For the sequential test described in~\Cref{def:bounded-means-test}, we have the following: 
                \begin{align}
                    \text{Under } H_0: \; &\mathbb{P}(\tau < \infty) \leq \alpha. \\
                    \text{Under } H_1: \; &\mathbb{P} \lp \tau <\infty \rp = 1,  \quad \mathbb{E}[\tau] = \mc{O}\lp  \frac{\sigma^2  \log(1/\alpha \Delta)}{\Delta^2}  +{\frac{\log(1/\alpha \Delta)}{\Delta}}\rp,  \\ 
                    \text{and} \quad &\lim_{n \to \infty} -\frac{1}{2n} \log   \mathbb{P} \lp \tau > n \rp   \geq \inf_{P' \in \mc{P}(\X)} \; \frac{1}{2} \lp \dkl\lp P', P_X \rp + \dkl \lp P', P_Y \rp \rp. 
                \end{align}
            \end{proposition}

            {
            \begin{remark}
                \label{remark:mmd}  The upper bound on the expected stopping time in~\Cref{prop:kernel-mmd} and~\Cref{prop:bounded-mean} displays the classical Bernstein type ``two-phase behavior'' corresponding to the low and high variance regimes. In particular, when the variance is low~(i.e., when $\sigma^2 \lesssim \Delta$), then $\mathbb{E}[\tau] = {\mc{O}}(\log(1/\alpha \Delta)/\Delta)$, while   in the high-variance regime~(i.e., when $\sigma^2 \gtrsim \Delta$) the bound is ${\mathcal{O}}(\log(1/\alpha\Delta)\sigma^2/\Delta^2)$. 
            \end{remark}
            }
            This result, proved in~\Cref{proof:kernel-mmd}, implies that the expected sample size required by  our sequential kernel-MMD test under the alternative has an inverse dependence on $\dmmd(P_X, P_Y)^2$. This is of the same order as the minimum number of observations needed by the  kernel-MMD test of~\textcite[Corollary~9]{gretton2012kernel} based on uniform deviation inequality.  This highlights the key benefit of our sequential test --- the ability to adapt the sample size automatically to the unknown alternative. This is not possible with fixed sample size tests in the absence of additional information about $\dmmd(P_X, P_Y)$.
            
            \paragraph{Optimality}  As in the case of bounded-mean testing, we can again show that there exist probability distributions for which the quadratic dependence of the expected sample size on the kernel-MMD distance between $P_X$ and $P_Y$ cannot be improved.
            \begin{proposition}
                \label{prop:lower-bound-kernel-mmd}
                Let $\tau'$ denote any level-$\alpha$, power-one sequential test for the problem defined in~\eqref{eq:two-sample-test-def} with $\mc{X} = \mathbb{R}$,   such that the expected value of $\tau'$ is finite under $H_1$. Then, there exist distributions $P_X$ and $P_Y$  supported  on $\mathbb{R}$, and a kernel $K$, such that   $\dmmd(P_X, P_Y) = \Delta$ and $\mathbb{E}[\|K(X, \cdot) - K(Y, \cdot)\|_K^2] = \Omega(\sigma^2)$ for any $\Delta, \sigma^2>0$, such that we have $\mathbb{E}[\tau'] = \Omega \lp \frac{\sigma^2 \log \lp1/\alpha \rp}{\Delta^2} \rp$. 
            \end{proposition}
            Thus this result implies that {in the high variance regime~(i.e., $\sigma^2\gtrsim\Delta$)}, the performance of our sequential kernel-MMD test has the optimal dependence on $\sigma^2$ and $\alpha$, and has a near-optimal (i.e., modulo $\log$) dependence on $\Delta$.  
            The proof of this statement is in~\Cref{proof:lower-bound}.

\section{Generalizations}
\label{sec:generalizations} 
    
    We now show how the game-theoretic testing formulation can be leveraged to design powerful sequential tests in much more general settings, beyond the standard two-sample testing problem. 
    In particular,
    we discuss the following  two extensions: 
    \begin{itemize}
        
        \item In~\Cref{subsec:time-varying}, we consider a time-varying version of the two-sample testing problem, and show that the framework of~\Cref{sec:general-two-sample} can still be used to construct a consistent test for this problem. 
        
        \item Next, in~\Cref{subsec:unified-approach}, we first introduce an abstract hypothesis testing problem~(with \iid observations) that unifies several testing problems such as two-sample testing, testing for independence and symmetry testing.
        We then introduce a testing strategy for this problem using a class of test function $\mc{G}$, similar to the approach of~\Cref{sec:general-two-sample}, and  characterize the statistical properties of the test in terms of the \emph{complexity} of $\mc{G}$.
        
    \end{itemize}

    \subsection{Time-varying two-sample testing}
    \label{subsec:time-varying}
        We begin  by defining a time-varying generalization of the  two-sample testing problem. 
        \begin{definition}[Two-sample testing with time-varying distributions]
           \label{def:time-varying-problem}
          Consider two independent sequence of observations, $\{X_t \in \mc{X}: t \geq 1\}$ and $\{Y_t \in \mc{X}: t \geq 1\}$, with $X_t \sim P_t$ and $Y_t \sim Q_t$ for all $t \geq 1$, where the distributions $P_t$ and $Q_t$ are $\mc{F}_{t-1} = \sigma(X_1^{t-1}, Y_1^{t-1})$ measurable. Let  $\boldsymbol{P}$ and $\boldsymbol{Q}$ denote the joint distributions of $\{P_t\in \mc{P}(\mc{X}): t \geq 1\}$ and  $\{Q_t \in \mc{P}(\mc{X}): t\geq 1\}$ respectively.
          For some family of  test functions $\mc{G} \subset [-1/2,1/2]^{\mc{X}}$ satisfying~\Cref{assump:function-class-G},  define
           \begin{align}
               \label{eq:distance-metric-time-varying}
               d_{\mc{G}}\lp \boldsymbol{P}, \boldsymbol{Q} \rp \defined  \sup_{g \in \mc{G}}\; \liminf_{n \to \infty}\;   \frac{1}{n} \sum_{t=1}^n \mathbb{E} \lb g(X_t) - g(Y_t) | \mc{F}_{t-1} \rb. 
           \end{align}
           In this setting, we consider the following hypothesis testing problem: 
           \begin{align}
               H_0: d_{\mc{G}}\lp \boldsymbol{P}, \boldsymbol{Q} \rp \eqas 0, \quad \text{versus} \quad H_1: d_{\mc{G}}\lp \boldsymbol{P}, \boldsymbol{Q} \rp \stackrel{a.s.}{>}0. 
           \end{align}
        \end{definition}
        
        \begin{remark}
        \label{remark:time-varying-example}
            The above problem can be considered as a game between  the {statistician} and an {adversary}. 
            The  adversary adaptively selects a sequence of distributions to induce the statistician to reject the null~(if null is true) or to continue sampling~(if the alternative is true). 
            In the special case of $P_t = P_X$, and $Q_t=P_Y$ for all $t\geq 1$, this problem reduces exactly to the usual two-sample testing problem considered in~\Cref{sec:general-two-sample}. Beyond this special case, this formulation covers a significantly larger set of problems. For example,  when $P_t = Q_t$ for all even $t$ and $\mathbb{E}_{P_t}[g^*(X_t)] - \mathbb{E}_{Q_t}[g^*(Y_t)] \geq  2\Delta$ for all odd $t$ for some $g^* \in \mc{G}$, the distance $d_{\mc{G}}(\boldsymbol{P}, \boldsymbol{Q})\geq \Delta$, even though $P_t = Q_t$ infinitely often. 
        \end{remark}
        
        \begin{remark}
            \label{remark:non-exchangeable} 
            \blue{
                Note that, for the problem stated in~\Cref{def:time-varying-problem},  the joint distributions of $(X_t)_{t=1}^n$ and $(Y_t)_{t=1}^n$, for any fixed $n$, are non-exchangeable in general. This renders permutation tests, a common approach for nonparametric testing in the batch setting, inapplicable to this problem.
            } 
        \end{remark}

        We now show that, similar to the standard two-sample testing problem, the existence of a  prediction strategy with small limiting average regret on a per-sequence basis implies the existence of a consistent sequential test. 
        \begin{proposition}
            \label{prop:time-varying} 
            Let $\tau \equiv  \tau(\ons, \prediction)$ denote a sequential test for the problem defined in~\Cref{def:time-varying-problem}. Then, for any prediction strategy $\prediction$, we have $\mathbb{P}\lp \tau < \infty \rp  \leq \alpha$, under the null  $H_0$. Furthermore, under $H_1$,   let $\prediction$ denote a  prediction strategy satisfying $\limsup_{n \to \infty} \frac{1}{n} \mc{R}_n \lp \prediction, \mc{G}, X_1^n, Y_1^n \rp < \dG(\boldsymbol{P}, \boldsymbol{Q})$ almost surely. Then, if~\Cref{assump:function-class-G} holds, the resulting test $\tau$ satisfies $\mathbb{P}\lp \tau < \infty \rp  = 1$.
        \end{proposition}
        The proof of this statement is given in~\Cref{proof:time-varying}. 
        
    \subsection{Unified approach to several testing problems}
    \label{subsec:unified-approach}
        We noted in~\Cref{remark:two-sample-testing-definition} that our strategy exploits the invariance of  the null distributions in two-sample testing  to an operator $\Tau$ that flips the order of the paired observations; that is $\Tau(X, Y) = (Y, X)$. 
        Building upon this observation, and following~\textcite{romano1989bootstrap}, we now state an abstract testing problem in which the null class contains precisely those distributions that are invariant to some given operator $\Tau$. 
        \begin{definition}[Abstract Hypothesis Testing Problem]
            \label{def:abstract-hypothesis-test} 
            Let $\nullclass$ and $\altclass$ denote disjoint classes of distributions on the observations space $\mc{Z}$. Let $\Tau: \mc{Z} \to \mc{Z}$ denote an operator on the observation space satisfying  $P \circ \Tau^{-1} = P$ for all $P \in \nullclass$, and $P \circ \Tau^{-1} \neq P$ for all $P \in \altclass$. 
            Given an \iid sequence of observations $Z_1, Z_2, \ldots, $  drawn according to $P \in \nullclass \cup \altclass$, we want to test 
                $H_0: P \in \nullclass$, versus $H_1: P \in \altclass$. 
        \end{definition}
        
        To proceed as in~\Cref{sec:general-two-sample}, we need an appropriate integral probability metric~(IPM). With $\mc{G}= \{g: \mc{Z} \to [-1/2,1/2]\}$ denoting a collection of test functions, define 
        \begin{align} 
            \label{eq:general-distinguishability} 
            \dG(P, P \circ \Tau^{-1}) \defined \sup_{g \in \mc{G}}\; \lv \mathbb{E}_{P}[g(Z)] - \mathbb{E}_{P}[g(\Tau Z)] \rv. 
        \end{align}
        To define a sequential test, we require that the IPM $\dG$ defined above, associated with the function class $\mc{G}$, is \emph{characteristic} for the family of distributions $\mc{P}_2 \subset \mc{P}(\mc{Z}) \times \mc{P}(\mc{Z})$, defined as 
        \begin{align}
        \label{eq:general-pair-of-distributions}
            \mc{P}_2 = \{(P, P\circ \Tau^{-1}): P \in \altclass\}. 
        \end{align}

        Under the above assumption, we can define a sequential test similar to the two-sample case in~\Cref{sec:general-two-sample}. In particular, assuming that $P_Z \in \altclass$ is the true distribution, we can find a function $g^*$ in $\mc{G}$, that maximizes the distance $\dG$ in~\eqref{eq:general-distinguishability}. With $\lambda^*$ representing $\argmax_{\lambda \in (-1,1)} \mathbb{E}_{P_Z}[\log\lp 1+\lambda\lp g(Z) - g( \Tau Z) \rp \rp$, we can define the \emph{oracle} test, $\tau^*$, as 
        \begin{align}
            \tau^* &= \min \{n \geq 1: \wealth^*_n \geq 1/\alpha \},  \\
            \text{where} \quad \wealth^*_n & = \wealth_{n-1}^* \times \lp 1 + \lambda^* \lp g^*(Z_n) - g^*( \Tau Z_n) \rp \rp, \quad \text{and} \quad \wealth^*_0 = 1. 
        \end{align}
        
        Following the same steps as in~\Cref{sec:general-two-sample}, we need to select a prediction strategy $\prediction$ for selecting  $\{g_t: t \geq 1\}$, estimates of the witness function $g^*$, and a betting strategy $\betting$, for selecting $\{\lambda_t: t \geq 1\}$; estimates of the optimal betting fraction $\lambda^*$.
        
        \begin{definition}[Abstract Sequential Test]
           \label{def:abstract-sequential-test}
           Let $Z_1, Z_2, \ldots$ denote an \iid sequence drawn from a distribution $P \in \nullclass \cup \altclass$, and let $\prediction$ denote any prediction strategy that generates a sequence $\{g_t \in \mc{G}: t \geq 1\}$, with $g_t$ being $\mc{F}_{t-1} = \sigma(Z_1, \ldots, Z_{t-1})$ measurable. Let $\{\lambda_t \in (-1,1): t \geq 1\}$ denote the sequence of bets generated by the ONS strategy introduced in~\Cref{def:ONS}. Then, for a given level $\alpha \in (0,1)$, we define the following sequential test: 
           \begin{align}
               \label{eq:abstract-sequential-test-def}
               \tau &= \min \{n \geq 1: \wealth_t \geq 1/\alpha\},  \\
                \text{where} \quad \wealth_t & = \wealth_{t-1} \times \lp 1 + \lambda_t \lp g_t(Z_t) - g_t( \Tau Z_t) \rp \rp\; \text{for } t \geq 1, \quad \text{and} \quad \wealth_0 = 1. 
           \end{align}
        \end{definition}
        It is easy to check that the wealth process $\{\wealth_t: t \geq 1\}$ is a test martingale under the null for any prediction strategy $\prediction$. Hence, $\tau$  controls type-I error at level $\alpha$ uniformly over the possibly composite null.
        Furthermore, as in~\Cref{sec:general-two-sample}, the power of our sequential test can again be characterized in terms of the regret of the prediction strategy $\prediction$, defined as 
        \begin{align}
            \label{eq:abstract-regret} 
            \mc{R}_n \equiv \mc{R}_n \lp \prediction, \mc{G}, \Tau, Z_1^n \rp  \defined \sup_{g \in \mc{G}} \sum_{t=1}^n \bigg( \big(g(Z_t) - g(\Tau Z_t) \big) - \big(g_t(Z_t) - g_t(\Tau Z_t) \big) \bigg). 
        \end{align}
        \begin{corollary}
            \label{corollary:abstract-test-1} 
            For the hypothesis testing problem defined in~\Cref{def:abstract-hypothesis-test}, let $\tau \equiv \tau(\ons, \prediction)$ denote the sequential test introduced in~\Cref{def:abstract-sequential-test} with function class $\mc{G}$.
            If $\mc{G}$ satisfies~\Cref{assump:function-class-G}, and the associated IPM $\dG$, defined in~\Cref{eq:general-distinguishability}, is characteristic~(\Cref{def:characteristic-IPM}) for the class of distributions defined in~\eqref{eq:general-pair-of-distributions},  we have the following: 
            \begin{itemize}
                \item  For any prediction strategy $\prediction$,  the test $\tau(\ons, \prediction)$ controls the type-I error uniformly over the null; that is, $\sup_{P \in \nullclass}  \mathbb{P}_{P}(\tau<\infty) \leq \alpha$. 
                
                \item If  $\prediction$ ensures that $\limsup_{n \to \infty} \frac{\mc{R}_n \lp \prediction, \mc{G}, \Tau, Z_1^n \rp}{n} < \dG\lp P, P \circ \Tau^{-1}\rp$ almost surely under the alternative, then the test $\tau$ is consistent. That is, $\mathbb{P}_P\lp \tau < \infty \rp = 1$ for all $P \in \altclass$. 
                
                \item If  $\prediction$ ensures that there exists a sequence $\{r_n: n\geq 1\}$ with $r_n \to 0$, and events $E_n = \{ \mc{R}_n/n\leq r_n \}$ with $\sum_{n=1}^{\infty} \mathbb{P}(E_n^c) < \infty$, then the expected stopping time satisfies the upper bound~\eqref{eq:main-test-stopping-time}. If we have the stronger conditions that $\mathbb{P}(E_n^c)=0$ for all $n \geq 1$, then the test $\tau$ also satisfies~\eqref{eq:main-test-exponent}. 
            \end{itemize}
        \end{corollary}
        The proof  follows the same steps as the proof of~\Cref{theorem:main-result}, and we omit the details.

       While~\Cref{corollary:abstract-test-1} identifies sufficient conditions for the consistency of the test $\tau$, it is non-constructive in nature. We now analyze the properties of our test $\tau$ initialized with  a natural prediction strategy,  called the empirical risk minimization~(ERM) strategy. 
       
       \begin{definition}[ERM strategy]
          \label{def:erm-strategy} 
          For a stream of observations $\{Z_t: t \geq 1\}$, the ERM prediction strategy, $\Aerm$,  selects $\{g_t\equiv g_t(Z_1^{t-1}): t \geq 1\}$ as follows:
          \begin{align}
                 \label{eq:erm-prediction-0}
                g_t \in \argmax_{g \in \mc{G}} \frac{1}{t-1} \sum_{i=1}^{t-1} g(Z_i) - g(\Tau Z_i), \quad \text{for all } t\geq 2,
          \end{align}
          and at $t=1$, $\Aerm$ sets $g_1$ to be an arbitrary element of $\mc{G}$. 
       \end{definition}
       We will analyze the performance of our test $\tau(\ons, \Aerm)$ under certain assumptions on the richness of the function class $\mc{G}$. A suitable measure of complexity is the Rademacher complexity, whose definition we recall next. 
       \begin{definition}
          \label{def:rademacher-complexity}
          Consider a  function class $\mc{H}$ containing mappings from some observations space $\mc{Z}$ to $\mathbb{R}$, and let $P \in \mc{P}(\mc{Z})$ denote a probability distribution on $\mc{Z}$. For a natural number $n \geq 1$, let $\boldsymbol{\sigma_n} = (\sigma_1, \ldots, \sigma_n)$ denote a random vector distributed uniformly over $\{-1,+1\}^n$.  Then, given $Z_1, Z_2, \ldots, Z_n$ drawn \iid from $P$, introduce the the following complexity terms: 
        \begin{align}
            \label{eq:rademacher-complexity-0} 
            C_n(\mc{H}, P) \defined \frac{1}{n} \mathbb{E}\lb  \sup_{h \in \mc{H}} \sum_{t=1}^n h(Z_t) \sigma_t \rb, \quad \text{and} \quad  C_n(\mc{H}) \defined \sup_{P \in \mc{P}(\mc{Z})} C_n(\mc{H}, P). 
          \end{align}
        \end{definition}%
        Before stating \Cref{theorem:abstract-power-result}, we need to introduce two more terms: the function class $\Gtilde$, and the notion of $\Delta$-separated alternatives, $\altclass(\Delta)$. 
        \begin{align}
            \Gtilde \defined \{\gtilde(\cdot) = g(\cdot) - g(\Tau \cdot), \; g \in \mc{G}\}, \quad \text{and} \quad \altclass(\Delta)=\{P \in \altclass: \dG(P, P\circ\Tau^{-1})>\Delta)\}. \label{eq:G-tilde}
        \end{align}
       We now present the main result of this section, that relates the consistency and detection boundary of $\tau(\ons, \Aerm)$ to the complexity of the function class $\Gtilde$.
        \begin{theorem}
            \label{theorem:abstract-power-result}
            For the sequential test $\tau \equiv \tau(\ons, \Aerm)$ for the  testing problem of~\Cref{def:abstract-hypothesis-test}, with prediction strategy  $\Aerm$ introduced in~\Cref{def:erm-strategy},  we have the following: 
            \begin{itemize}
                \item $\tau$ is consistent against any $P \in \altclass$, for which  $C_n(\Gtilde, P)$ converges to $0$, that is,
                \begin{align}
                     \lim_{n \to \infty} C_n(\Gtilde, P) \to 0  \quad \Rightarrow \quad \mathbb{P}_P(\tau < \infty) = 1. 
                \end{align}
                \item Suppose $C_n(\Gtilde)$ converges to $0$ with $n$, and for  a small $\gamma \in (0,1)$, introduce the term 
                \begin{align}
                    \Delta_n^* =  \sqrt{ \frac{ 8 \log{n/\alpha}}{n}} + \frac{2}{n} \lp 2 + \sum_{t=1}^{n-1} \lp C_t(\Gtilde) + 5\sqrt{ \frac{\log(16n/\gamma)}{2t}} \rp \rp + \sqrt{\frac{8 \log(4/\gamma)}{n}}.
                \end{align}
                
                Then, for any $n \geq 1$, and $\Delta_n>\Delta_n^*$, 
                \begin{align}
                     \sup_{P \in \altclass(\Delta_n)} \mathbb{P}_P(\tau >n ) \leq \gamma. 
                \end{align}
                In other words, $\Delta_n^*$ denotes the minimum separation that can be detected with power greater than $1-\gamma$ by our sequential test within the first $n$ observations.
            \end{itemize}
        \end{theorem}
        The proof of this statement is given in~\Cref{proof:abstract-power-result}. 
        \begin{remark}
            \label{remark:erm-special-case}
            The above result implies that the detection boundary for our test in terms of the $\dG$ distance measure is given by $\Delta_n^* =  \Omega\lp \frac{1}{n} \sum_{t=1}^{n-1} C_t(\Gtilde) +  \sqrt{\frac{\log(n/\alpha)}{n}} + \sqrt{\frac{\log(n/\gamma)}{n}}\rp$, where $\alpha$ and $\gamma$ correspond to the type-I and type-II errors. For the bounded-mean test~(\Cref{sec:bounded-mean-testing}) and the kernel-MMD test~(\Cref{subsec:kernel-mmd}) introduced earlier, it is known that $C_t(\Gtilde)$ decays to zero at a $1/\sqrt{t}$ rate. Hence, for both these tests, we have $\Delta_n^* = \Omega \lp   \sqrt{\frac{\log(n/\alpha)}{n}} + \sqrt{\frac{\log(n/\gamma)}{n}}\rp$. 
        \end{remark}
        
        As mentioned earlier, the abstract test of~\Cref{def:abstract-hypothesis-test} with $\Tau$ such that $\Tau(X,Y)= (Y,X)$ reduces to the two-sample testing problem. We now show that two other important testing problems are also covered by this definition. 
        
        \subsubsection{Testing for symmetry} We state the simplest version, in which $\mc{Z} = \mathbb{R}$, and we assume that the null distributions are symmetric about the origin. That is, $\Tau: \mc{Z} \to \mc{Z}$, such that $\Tau z = -z$. Hence, for any continuous $P \in \mc{P}(Z)$, we have $P\circ \Tau^{-1} =Q$ such that $F_P(z) = 1 - F_Q(-z)$ for all $z \in \mc{Z}$. 
        To define the sequential test, we can use the function class $\mc{G} = \{g(x)= \indi{x \leq u}: u \in \mc{Z}\}$. Then the distance $\dG$ in~\eqref{eq:general-distinguishability} reduces to the KS distance between $P$ and $P\circ \Tau^{-1}$:
        \begin{align}
            \sup_{g \in \mc{G}} |\mathbb{E}_{P}[g(Z)] - \mathbb{E}_{P}[g(\Tau Z)]| = \sup_{u \in \mathbb{R}} |F_{P}(u) - 1 + F_P(-u)| > 0.
        \end{align}
        \Cref{theorem:abstract-power-result} implies that the test with ERM prediction strategy is consistent, and has a detection boundary of the order $\mc{O}(\sqrt{\log n/n})$. 
        
        \subsubsection{Testing for independence} In this case, we have two observation spaces $\mc{X}$ and $\mc{Y}$, not necessarily the same, and define $\mc{Z} = \lp \mc{X} \times \mc{Y} \rp \times \lp \mc{X} \times \mc{Y} \rp$. Let $P_{XY}$ denote a distribution in $\mc{P}(\mc{X} \times \mc{Y})$, and let $P_X \in \mc{P}(\mc{X})$  and $P_Y \in \mc{P}(\mc{Y})$ denote its marginals. Under the null hypothesis, we have $P_{XY} = P_X \times P_Y$, which can be encoded via the operator $\Tau: \mc{Z} \to \mc{Z}$, with $\Tau((x, y), (x', y')) = ((x,y'), (x', y))$. 
       
         When $\mc{X}=\mc{Y}=\mathbb{R}$, we can again select $\mc{G} = \{g(x) = \indi{x \leq u}: u \in \mbb{R}\}$, which leads to  $\dG$ being the KS distance between $P_{XY}$ and the product of its marginals $P_X \times P_Y$. For general $\mc{X} \neq \mc{Y}$, a suitable choice of $\mc{G}$ is a norm ball in the RKHS of the product kernel $K((x, y), (x', y')) \defined K_X(x, x') K_Y(y, y')$ for positive definite kernels $K_X: \mc{X}\times \mc{X} \to \mbb{R}$ and $K_Y: \mc{Y}\times \mc{Y} \to \mbb{R}$. In this case, the distance $\dG$ is the kernel-MMD distance between $P_{XY}$ and $P_{X}\times P_Y$; also called the HSIC criterion \parencite{gretton2005measuring}. In both cases,~\Cref{theorem:abstract-power-result} implies that the test is ERM strategy is consistent, and furthermore has a detection boundary of the order $\mc{O}(\sqrt{\log n /n})$ in their respective distance metrics.

\section{Numerical Simulations}
\label{sec:numerical-experiments}
  
    In this section, we demonstrate, through experiments,  the key advantages of our betting-based sequential kernel-MMD test over existing batch and sequential tests. 
     \begin{figure}[htb!]
        \def\figwidth{0.80\linewidth}
        \def\figheight{0.30\textheight} %
        \centering
           \input{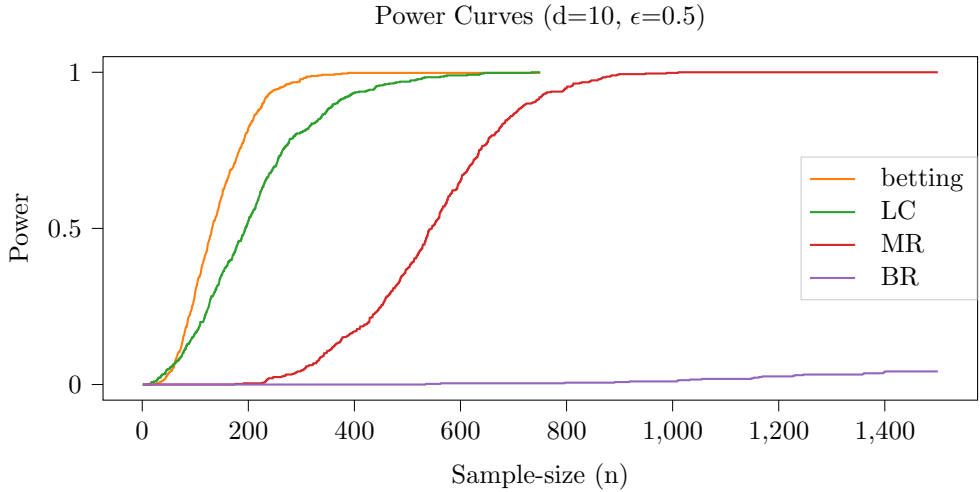}
        \caption{Comparison of the power of our sequential kernel-MMD test~(betting) with  other sequential tests~(LC, MR, BR). 
        }
        \label{fig:experiments-1}
   \end{figure}  
    \begin{figure}[htb!]
        \def\figwidth{0.80\linewidth}
        \def\figheight{0.30\textheight} %
        \centering
           \input{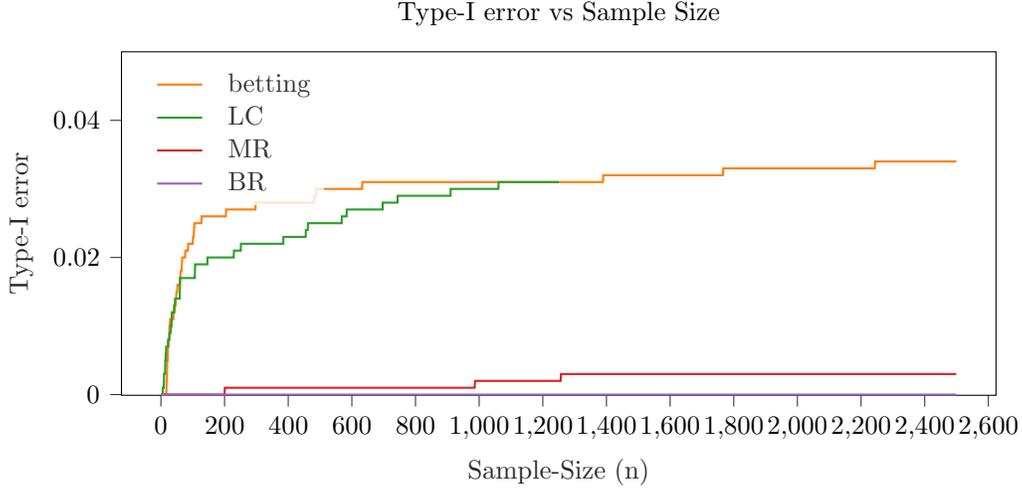}
        \caption{Comparison of the type-I error of our sequential kernel-MMD test~(betting) with  other sequential tests~(LC, MR, BR). 
        }
        \label{fig:experiments-type-I}
   \end{figure}  
   
   \begin{figure}[htb!]
       \centering
        \def\figwidth{0.80\linewidth}
        \def\figheight{0.30\textheight} %
       \begin{tikzpicture}[baseline]

\definecolor{crimson2143940}{RGB}{214,39,40}
\definecolor{darkgray176}{RGB}{176,176,176}
\definecolor{darkorange25512714}{RGB}{255,127,14}
\definecolor{forestgreen4416044}{RGB}{44,160,44}
\definecolor{lightgray204}{RGB}{204,204,204}
\definecolor{mediumpurple148103189}{RGB}{148,103,189}
\definecolor{steelblue31119180}{RGB}{31,119,180}
\pgfplotsset{
    yticklabel style={/pgf/number format/fixed},  
}
\begin{axis}[
height=\figheight,
legend cell align={left},
legend style={
  fill opacity=0.8,
  draw opacity=1,
  text opacity=1,
  at={(0.97,0.03)},
  anchor=south east,
  draw=lightgray204
},
tick align=outside,
tick pos=left,
title={Adaptivity of sequential kernel-MMD test with Gaussian kernel},
width=\figwidth,
x grid style={darkgray176},
xlabel={Sample-size (n)},
xmin=-49.5, xmax= 950, %
xtick style={color=black},
y grid style={darkgray176},
ylabel={Power},
ymin=0.0165, ymax=1.04683333333333,
ytick style={color=black}
]
\addplot [semithick, steelblue31119180]
table {%
10 0.116666666666667
22 0.243333333333333
35 0.4
47 0.513333333333333
60 0.596666666666667
73 0.74
85 0.796666666666667
98 0.88
111 0.93
123 0.973333333333333
136 0.966666666666667
148 0.993333333333333
161 1
174 0.993333333333333
186 0.99
199 0.996666666666667
212 1
224 1
237 1
250 1
};
\addlegendentry{$\epsilon$=0.70}
\addplot [semithick, steelblue31119180, dashed, forget plot]
table {%
149.73 0.0165
149.73 1.04683333333333
};
\addplot [semithick, darkorange25512714]
table {%
10 0.0866666666666667
33 0.176666666666667
56 0.323333333333333
79 0.43
102 0.573333333333333
125 0.65
148 0.756666666666667
172 0.85
195 0.876666666666667
218 0.906666666666667
241 0.956666666666667
264 0.966666666666667
287 0.98
311 0.973333333333333
334 0.996666666666667
357 0.996666666666667
380 1
403 0.996666666666667
426 1
450 1
};
\addlegendentry{$\epsilon$=0.50}
\addplot [semithick, darkorange25512714, dashed, forget plot]
table {%
288.586666666667 0.0165
288.586666666667 1.04683333333333
};
\addplot [semithick, forestgreen4416044]
table {%
10 0.0666666666666667
41 0.166666666666667
72 0.263333333333333
103 0.376666666666667
134 0.496666666666667
165 0.58
196 0.696666666666667
227 0.76
258 0.826666666666667
289 0.866666666666667
320 0.886666666666667
351 0.943333333333333
382 0.976666666666667
413 0.966666666666667
444 0.98
475 0.993333333333333
506 1
537 0.993333333333333
568 0.996666666666667
600 1
};
\addlegendentry{$\epsilon$=0.40}
\addplot [semithick, forestgreen4416044, dashed, forget plot]
table {%
437.37 0.0165
437.37 1.04683333333333
};
\addplot [semithick, crimson2143940]
table {%
10 0.0633333333333333
62 0.196666666666667
114 0.29
166 0.433333333333333
218 0.6
270 0.716666666666667
322 0.773333333333333
374 0.876666666666667
426 0.923333333333333
478 0.963333333333333
531 0.97
583 0.99
635 0.983333333333333
687 0.996666666666667
739 0.996666666666667
791 1
843 1
895 1
947 1
1000 1
};
\addlegendentry{$\epsilon$=0.35}
\addplot [semithick, crimson2143940, dashed, forget plot]
table {%
583.72 0.0165
583.72 1.04683333333333
};
\addplot [semithick, mediumpurple148103189]
table {%
10 0.0633333333333333
72 0.173333333333333
135 0.286666666666667
197 0.433333333333333
260 0.553333333333333
323 0.686666666666667
385 0.746666666666667
448 0.786666666666667
511 0.903333333333333
573 0.9
636 0.943333333333333
698 0.96
761 0.99
824 0.97
886 0.996666666666667
949 0.993333333333333
1012 0.996666666666667
1074 1
1137 1
1200 1
};
\addlegendentry{$\epsilon$=0.30}
\addplot [semithick, mediumpurple148103189, dashed, forget plot]
table {%
870.17 0.0165
870.17 1.04683333333333
};
\end{axis}

\end{tikzpicture}
       \caption{This figure demonstrates the ability of our sequential test to adapt to the unknown hardness of the problem: we set $P_X \sim N(\boldsymbol{O}, I_m)$ and $P_Y \sim N(a_{\epsilon, 1}, I_m)$, where $a_{\epsilon, 1}$ is obtained by setting the first coordinate of $\boldsymbol{0} \in \mathbb{R}^m$ to $\epsilon$. Each solid curve is obtained by running~($200$ trials of) the kernel-MMD permutation test for $20$ different values of the sample size, while the dashed vertical line shows the corresponding average stopping time of our sequential kernel-MMD test.} 
       \label{fig:adaptivity}
   \end{figure}
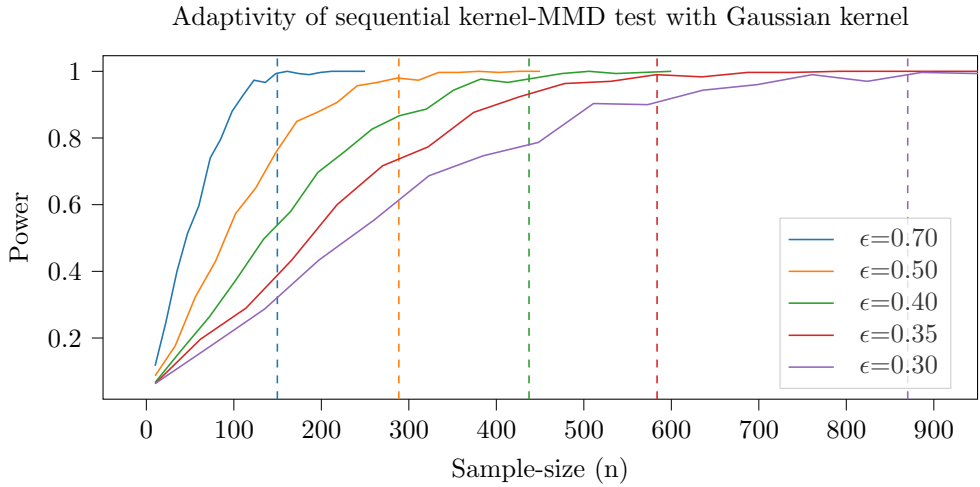
   
    \begin{figure}[htb!]
       \centering
        \def\figwidth{0.80\linewidth}
        \def\figheight{0.30\textheight} %
       \begin{tikzpicture}

\definecolor{crimson2143940}{RGB}{214,39,40}
\definecolor{darkgray176}{RGB}{176,176,176}
\definecolor{darkorange25512714}{RGB}{255,127,14}
\definecolor{forestgreen4416044}{RGB}{44,160,44}
\definecolor{lightgray204}{RGB}{204,204,204}
\definecolor{mediumpurple148103189}{RGB}{148,103,189}
\definecolor{steelblue31119180}{RGB}{31,119,180}

\begin{axis}[
height=\figheight,
legend cell align={left},
legend style={
  fill opacity=0.8,
  draw opacity=1,
  text opacity=1,
  at={(0.97,0.03)},
  anchor=south east,
  draw=lightgray204
},
tick align=outside,
tick pos=left,
title={Adaptivity of sequential kernel-MMD test with unbounded kernel},
width=\figwidth,
x grid style={darkgray176},
xlabel={Sample-size (n)},
xmin=-29.5, xmax=839.5,
xtick style={color=black},
y grid style={darkgray176},
ylabel={Power},
ymin=0.0025, ymax=1.0475,
ytick style={color=black}
]
\addplot [semithick, steelblue31119180]
table {%
10 0.21
27 0.583333333333333
44 0.84
61 0.97
79 0.99
96 1
113 1
130 1
148 1
165 1
182 1
200 1
};
\addlegendentry{2.0}
\addplot [semithick, steelblue31119180, dashed, forget plot]
table {%
66.44 0.0025
66.44 1.0475
};
\addplot [semithick, darkorange25512714]
table {%
10 0.0966666666666667
36 0.303333333333333
62 0.503333333333333
89 0.75
115 0.79
141 0.896666666666667
168 0.94
194 0.953333333333333
220 0.996666666666667
247 1
273 1
300 0.996666666666667
};
\addlegendentry{1.2}
\addplot [semithick, darkorange25512714, dashed, forget plot]
table {%
179.836666666667 0.0025
179.836666666667 1.0475
};
\addplot [semithick, forestgreen4416044]
table {%
10 0.0666666666666667
45 0.256666666666667
80 0.49
116 0.6
151 0.753333333333333
187 0.87
222 0.923333333333333
258 0.97
293 0.97
329 0.986666666666667
364 0.993333333333333
400 0.996666666666667
};
\addlegendentry{1.0}
\addplot [semithick, forestgreen4416044, dashed, forget plot]
table {%
262.16 0.0025
262.16 1.0475
};
\addplot [semithick, crimson2143940]
table {%
10 0.07
63 0.23
117 0.426666666666667
170 0.566666666666667
224 0.71
278 0.823333333333333
331 0.916666666666667
385 0.956666666666667
439 0.976666666666667
492 0.976666666666667
546 0.996666666666667
600 0.996666666666667
};
\addlegendentry{0.8}
\addplot [semithick, crimson2143940, dashed, forget plot]
table {%
419.7 0.0025
419.7 1.0475
};
\addplot [semithick, mediumpurple148103189]
table {%
10 0.05
81 0.166666666666667
153 0.326666666666667
225 0.433333333333333
297 0.566666666666667
369 0.71
440 0.76
512 0.84
584 0.913333333333333
656 0.953333333333333
728 0.96
800 0.99
};
\addlegendentry{0.6}
\addplot [semithick, mediumpurple148103189, dashed, forget plot]
table {%
780.653333333333 0.0025
780.653333333333 1.0475
};
\end{axis}

\end{tikzpicture}
       \caption{This figure demonstrates the ability of our sequential test to adapt to the unknown hardness of the problem: we set $P_X \sim N(\boldsymbol{O}, I_m)$ and $P_Y \sim N(a_{\epsilon, 1}, I_m)$, where $a_{\epsilon, 1}$ is obtained by setting the first coordinate of $\boldsymbol{0} \in \mathbb{R}^m$ to $\epsilon$. Each solid curve is obtained by running~($200$ trials of) the kernel-MMD permutation test with linear kernel for $20$ different values of the sample size, while the dashed vertical line shows the corresponding average stopping time of our sequential kernel-MMD test.} 
       \label{fig:adaptivity-unbounded}
   \end{figure}
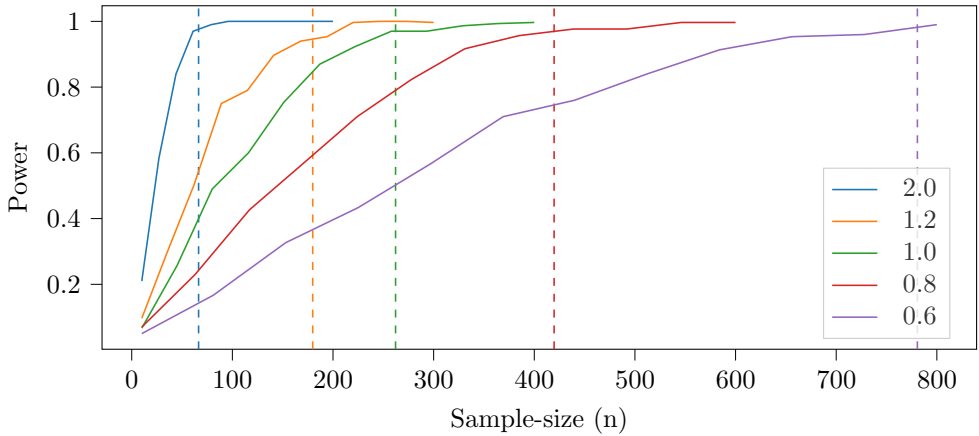
 
    \paragraph*{Experiment Setup} We will consider the two-sample testing problem with observations taking values in the set $\mc{X} = \mathbb{R}^m$ for some $m \geq 1$. In all cases, we will fix the distribution $P_X$ to $N(\boldsymbol{0}, I_m)$, where $I_m$ is the $m\times m$ identity matrix. For any integer $1\leq j \leq m$, and any $\epsilon \in \reals$, let $a_{\epsilon, j} \in \mathbb{R}^m$ denote the element obtained by changing the first $j$ coordinates of $\boldsymbol{0} \in \reals$ to $\epsilon$. We  set the distribution $P_{Y}$ to $N(a_{\epsilon, j}, I_m)$ for different choices of $j, \epsilon$. For the first two experiments we use the Gaussian kernel $K(x, y) = \exp \lp - \frac{1}{2 b^2} \|x-y\|_2^2 \rp$, with the bandwidth $b$ set to $\sqrt{m}$. In the third experiment, we use a linear kernel $K(x,y) = x^{T} y$, to study the performance of our modified sequential test with unbounded kernels, discussed in~\Cref{subsec:unbounded-functions}. 
    
    \paragraph*{Experiment 1: Comparison of power and type-I error with existing sequential tests} We  compare the power  of our kernel-MMD test with the following baselines:  \textbf{(i)} the sequential test of~\textcite{lheritier2019low}, denoted by LC, \textbf{(ii)} the sequential test proposed by~\textcite{manole2023martingale} based on reverse submartingales~(MR), and \textbf{(iii)} the linear-time sequential test of~\textcite{balsubramani2016sequential}, denoted by BR. 
    The power curves of the above tests are shown in~\Cref{fig:experiments-1}. 
    For each test, we run $400$ trials, and the power curves are obtained by plotting the empirical CDF of the stopping times. 

    Similarly, we plot the type-I error curves of the four sequential tests under the null in~\Cref{fig:experiments-type-I}, averaged over $1000$ trials. As expected, the MR and BR tests are significantly more conservative as compared to the martingale based tests~(betting and LC).

    \paragraph*{Experiment 2: Verification of adaptivity}~\Cref{prop:kernel-mmd} shows that the expected stopping time of our kernel-MMD test adapts to the hardness of the problem, without any additional information. To verify this claim, we again fix $m=10$, $j=1$ and vary $\epsilon$ in the range $\{0.3, 0.35, 0.4, 0.45, 0.5, 0.7\}$. For every $\epsilon$ value, we run $200$ trials of our kernel-MMD test and obtain the expected stopping time. A a baseline, we run $200$ trials each of the kernel-MMD permutation test with $150$ permutation, at $20$ different sample size values to obtain the power curve. The results, shown in~\Cref{fig:adaptivity}, demonstrate the ability of our sequential test to adapt the expected stopping time to the alternative.

    \paragraph*{Experiment 3: Sequential test with unbounded~(linear) kernel} The sequential kernel-MMD test, defined in~\Cref{subsec:kernel-mmd}, requires the kernel $K$ to be uniformly bounded; that is, $\sup_{x \in \mc{X}}K(x,x) \leq B<\infty$. However, as we discussed in~\Cref{subsec:unbounded-functions}, this restriction can be easily addressed within our framework by using an anti-symmetric sigmoidal function $\varrho$ taking values in $[-1,1]$. To empirically verify the performance of our modified sequential test with $\varrho(\cdot)=\text{tanh}(\cdot)$, we repeat the previous experiment with an unbounded linear kernel $K(x, y) = x^{T}y$. The result is plotted in~\Cref{fig:adaptivity-unbounded}, and shows the adaptivity of our sequential test even in the case of unbounded kernels.

\section{Discussion}
\label{sec:discussion}

    \subsection{Advantage of sequential tests} Our sequential tests have the following advantages compared to existing batch tests: 
    \begin{itemize}
        \item \emph{Adaptivity to alternative:} As shown analytically in~\Cref{prop:bounded-mean} and~\Cref{prop:kernel-mmd}, as well as empirically in~\Cref{sec:numerical-experiments}, the expected stopping times of our sequential tests automatically adapt to the unknown \emph{hardness} of the problem, measured by~$\dG(P_X, P_Y)$. This is in contrast to the fixed-sample size tests, where to choose the `right' sample size, we require additional prior information about the problem~(i.e., a lower bound on $\dG(P_X, P_Y)$). 
        
        \item \emph{Lower computational complexity:} Selecting the rejection threshold in fixed-sample size tests is often a nontrivial and computationally expensive task. For example, kernel-MMD permutation test requires recomputing the quadratic-time kernel-MMD statistic $b$ times, with $b$ usually in $[200, 1000]$.  In contrast, the rejection threshold for our level-$\alpha$ sequential tests is $1/\alpha$ --- a direct consequence of using \emph{test martingales} in our design. This leads to  significantly lower running times of our sequential test in comparison to its fixed sample size counterpart. More importantly, this threshold does not lead to an overly conservative test: Ville's inequality  holds with equality for continuous-time nonnegative martingales, and often holds almost with equality for our discrete-time nonnegative martingales.
        
    \end{itemize}
    
    \subsection{Benefits of game-theoretic formulation} Our design strategy is based on a game-theoretic view of sequential testing~(\Cref{subsec:testing_by_betting}), recently popularized by~\textcite{shafer2021testing}, where the gain in the wealth of a (fictitious) bettor  has a precise interpretation as the strength of evidence against the null. Working in this framework has two main benefits:
    \begin{itemize}
        \item \emph{Connections to online learning:} As we described in~\Cref{theorem:main-result}, the game-theoretic approach allows us to connect the statistical properties of the sequential test to the regret achieved by the prediction strategy in an associated online learning problem. Since there exists a well-developed theory of online learning algorithms ~\parencite{cesa2006prediction} for a wide range of  function classes $\mc{G}$, our work provides a simple method for  using these algorithms to design new sequential tests with strong performance guarantees. 
        
        \item \emph{Extension to non-\iid observations:}
        Many practical two-sample tests in the batch  setting rely strongly on the observations being  \iid, or at least exchangeable. For example, the kernel-MMD permutation test considered in our experiments in~\Cref{sec:numerical-experiments} uses the exchangeability of the observations to obtain the rejection threshold. 
        However, these conditions are often not satisfied in  many  applications, preventing the use of such  tests. In contrast, our techniques using the game-theoretic framework easily extend to a variant of the two-sample testing with time-varying probability distributions, as we described in~\Cref{subsec:time-varying}. 
    \end{itemize}
    
    \subsection{Working with unpaired observations}   
    \label{subsec:unpaired-observations} 
        In~\Cref{sec:general-two-sample}, we developed our two-sample testing framework under the assumption that, in each round $t$, we observe the pair $(X_t, Y_t)$ drawn from $P_X \times P_Y$. However, in many applications, the observations arrive in batches, often of unequal sizes. That is, in round $t$, we observe $\{X_{t,i}: 1 \leq i \leq n_t\}$ and $\{Y_{t,j}: 1 \leq j \leq m_t\}$ drawn \iid from $P_X$ and $P_Y$ respectively,  with $n_t$ possibly different from $m_t$. Such observation models can be easily handled in our framework by averaging the payoffs in the update rule for the test martingale, as follows: 
        \begin{align}
            \wealth_{t} &= \wealth_{t-1} \times \lp 1 + \lambda_t\lp \frac{1}{n_t} \sum_{i=1}^{n_t} g_t(X_{t,i}) - \frac{1}{m_t} \sum_{j=1}^{m_t} g_t(Y_{t,j}) \rp \rp \label{eq:batch-update-rule} \\
            &= \wealth_{t-1} \times \lp 1 + \lambda_t\lp \frac{1}{n_t m_t} \sum_{i=1}^{n_t} \sum_{j=1}^{m_t} g_t(X_{t,i}) -  g_t(Y_{t,j}) \rp \rp. \label{eq:batch-update-rule-2}
        \end{align}
        In some applications, instead of batches,  we only have access to a single stream  of observations $\{(Z_t, L_t): t \geq 1\}$, with $L_t \in \{0,1\}$ and $L_t=0$ indicates $Z_t \sim P_X$, while $L_t=1$ implies that $Z_t \sim P_Y$. The averaging idea, described above in~\eqref{eq:batch-update-rule}, can be used to address this case as well, with the modification that we only update the test martingale when we have at least one observation from both distributions. More formally, we update the test martingale using~\eqref{eq:batch-update-rule} at random stopping times $\{T_s: s \geq 1\}$, defined as follows: 
        \begin{align}
        T_0=0, \quad \text{and} \quad T_s = \min \big\{ t>T_{s-1}: 0 < \sum_{t'=T_{s-1}+1}^t L_{t'} < t-T_{s-1}\big\}.  
        \end{align}
        Together with the results of~\Cref{sec:general-two-sample} and~\Cref{subsec:time-varying}, the above discussion implies that the framework developed in this paper can be used to design consistent sequential tests that work under significantly weaker assumptions than existing tests in literature. In particular,  our strategy works with arbitrary observation models, with possibly dependent streams of observations~(\Cref{remark:non-independent-streams}), that are drawn from  time-varying distributions~(\Cref{subsec:time-varying}). 
           
\section{Conclusion}
\label{sec:conclusion}
    In this paper, we described a general strategy of constructing  sequential tests for the two-sample testing problem and its generalizations. The fundamental idea underlying our approach is the principle of testing by betting, which motivates a game-theoretic formulation of the problem. We presented a general strategy of constructing sequential tests within this framework based on a class of integral probability metrics~(IPMs), and instantiated this strategy for the  kernel-MMD metric. Both theoretical and empirical results demonstrate the computational efficiency and the ability of the test to adapt to unknown alternatives.
    
    \red{
    Since the framework developed in our paper is quite general, adapting it to new testing problems is an interesting direction for future work. For instance, inspired by this work,~\textcite{shaer2022model} very recently proposed a conditional independence test based on symmetry testing, concurrently developing ideas similar to those discussed in~\Cref{subsec:unbounded-functions}. Another important direction is to instantiate our framework using other distance measures with variational representations, such as Wasserstein metric and $f$-divergences, and examine the stopping times and exponents that result. Finally, a rigorous empirical evaluation of the performance of the ``invariant" tests proposed in~\Cref{subsec:unified-approach}  is also an important question. 
    }

\subsection*{Acknowledgements}
{
We thank Amaury Durand and Olivier Wintenberger for informing us about an error in the justification of~\eqref{eq:main-test-stopping-time} in the previous version of the manuscript. A standalone note discussing the error and its resolution can be found at \href{https://github.com/sshekhar17/nonparametric-testing-by-betting/blob/main/Correction.pdf}{https://github.com/sshekhar17/nonparametric-testing-by-betting/blob/main/Correction.pdf}.
}
\newpage 
\printbibliography
\newpage 
\begin{appendix}
\addtocontents{toc}{\protect\setcounter{tocdepth}{0}}
    \section{Additional Background}
        \label{appendix:background}
        
    \subsection{Testing by Betting}
    \label{subsec:testing_by_betting}
       
       The principle of \emph{testing by betting}~\parencite{shafer2021testing} provides the conceptual foundation for our testing strategy.  
        The basic idea can be stated as follows~\parencite[\S~2]{shafer2021testing}: \emph{the claim that a random variable $Z$ is distributed according to $P_Z$ can be interpreted equivalently as the offer of a bet with any payoff function  sold for its expected value under the distribution $P_Z$}. 
        This principle implies that  a person claiming that $Z \sim P_Z$ (a ``\forecaster'') should be willing to play any betting game that is fair under this claim, since he is not expected to lose money in such a game. Following the \forecaster's claim about $Z$, a \skeptic (who does not believe the \forecaster) can announce a betting score $S: \mathcal Z \to \mathbb{R}_+$ which represents the amount that  \skeptic's wealth is multiplied by: if the outcome is $z$, then \skeptic gets back $S(z)$ dollars for each dollar that was bet. This bet is ``fair'' from the point of view of the \forecaster if $\mathbb{E}_{P_Z}[S(Z)] \leq 1$, because it implies that the \skeptic cannot make money if the \forecaster is correct. So if a \skeptic chooses to bet against the \forecaster with any such pre-announced betting score, the \forecaster will be happy to play this game.

        We can formally describe this as a repeated  game involving three players, \forecaster, \skeptic and \reality following \textcite{shafer2019game, shafer2021testing}. 
        
        \begin{definition}[\ttt{Betting Protocol}] 
        \label{def:betting_protocol}
            Before the start of the game, \forecaster declares that the observations $\{Z_t: t \geq 1\}$ taking values in $\mc{Z}$ are distributed \iid according to some $P_Z$ in the class of distributions  $\mc{P}_0$. 
            \skeptic begins with an initial wealth  $\wealth_0 = 1$, and the game proceeds as follows for $t=1, 2, \ldots :$
            \begin{itemize}[leftmargin=3em]
                \item \skeptic selects a  function $\gtilde_t:\mc{Z} \to [-1, 1]$ such that $\mbb{E}_{P'}\lb \gtilde_t(Z) \vert \mc{F}_{t-1} \rb = 0$ for all $P' \in \mc{P}_0$, and $\mc{F}_{t-1} = \sigma \lp Z_1, \ldots, Z_{t-1} \rp$.
                
                \item \skeptic bets an amount $\lambda_t \wealth_{t-1}$ for an $\mc{F}_{t-1}$-measurable $\lambda_t \in [0,1]$ on the next realization of $Z$.
                
                \item \reality reveals the next realization $Z_t$.
                
                \item The wealth of the \skeptic is updated as $\wealth_t = \wealth_{t-1} +  \lambda_t \wealth_{t-1} \gtilde_t(Z_t) =  \wealth_{t-1} \lp 1 + \lambda_t \gtilde_t(Z_t) \rp$. 
            \end{itemize}
        \end{definition}
        
        \begin{remark}
            \label{remark:testing_by_betting}
            In the two-sample testing problem, the role of \forecaster is played by the null hypothesis, the role of \skeptic is played by the statistician and the role of \reality is played by the independent and identically distributed~(\iid) source, generating the observations. 
        \end{remark}

        \subsection{Ville's inequality}  We now recall a time-uniform analog of Markov's inequality,  derived by~\textcite{ville1939etude}. 
            \begin{fact}[Ville's Inequality]
            \label{fact:ville}
                Suppose $\{\wealth_t: t \geq 0\}$ is a nonnegative supermartingale process adapted to a filtration $\{\mc{F}_t: t \geq 0\}$. Then, we have,   for any $a>0$, $\mbb{P}\lp \exists t \geq 1 :  \wealth_t \geq a \rp  \leq \frac{\mbb{E}[\wealth_0]}{a}$.  
            \end{fact}

        \subsection{Details of ONS betting strategy}
        \label{appendix:background-ONS}
            We used the ONS betting strategy proposed by~\textcite[Algorithm 1]{cutkosky2018black}. The key property of this betting strategy that we repeatedly  use in designing our tests is that for any sequence of outcomes $\{v_t \in [-1,1]: t \geq 1\}$, the wealth process satisfies the following for any $n \geq 1$.  
            \begin{align}
                \wealth_n & \geq \frac{1}{\sum_{t=1}^n v_t^2} \exp \lp \frac{ (\sum_{t=1}^n v_t)^2}{ 4 \lp \sum_{t=1}^n v_t^2 + \sum_{t=1}^n v_t \rp}  \rp
                \geq \exp \lp \frac{n}{8} \lp \frac{1}{n} \sum_{t=1}^n v_t \rp^2 - \log n \rp.  \label{eq:ons-wealth}
            \end{align}
            This statement can be extracted from a more general result obtained by~\textcite{cutkosky2018black} while proving their Theorem~1. 
            
        \subsection{Regret Bound for OGA strategy}
        \label{appendix:background-oga}
            Let $\mc{U}$ denote an inner-product space, and let $\mc{U}_1$ and $\mc{U}_2$ denote two bounded subsets of $\mc{U}$. Consider the following prediction problem: For $t=1, 2, \ldots$: 
            \begin{itemize}
                \item Player plays $u_t \in \mc{U}_1$ 
                \item Adversary selects $v_t \in \mc{U}_2$ 
                \item Player gains reward $\langle u_t, v_t \rangle$. 
            \end{itemize}
            For any prediction strategy, $\mc{A}$, that selects $u_1, u_2, \ldots$, the regret after $n$ rounds is defined as 
            \begin{align}
                \mc{R}_n \equiv \mc{R}_n(\mc{A}, \mc{U}_1,  v_1^n) = \sup_{u \in \mc{U}_1} \sum_{t=1}^n \langle u - u_t, v_t \rangle. 
            \end{align}
            A useful strategy for this problem is the Online Gradient Ascent strategy~($\oga$), that sets $u_0$ as an arbitrary element of $\mc{U}_1$ and proceeds as follows for $t \geq 1$
            \begin{align}
                u_{t} = u_{t-1} + \eta_t v_{t-1}, 
            \end{align}
            where $\{\eta_t: t \geq 1\}$ are a non-negative decreasing sequence of step-sizes. 
            
            A standard result in online-learning~\parencite[Theorem 2.13]{orabona2019modern} implies the following regret bound for the OGA strategy: 
            \begin{align}
                \mc{R}_n \lp \oga \rp \equiv \mc{R}_n \lp \oga, \mc{U}_1, v_1^n \rp  &\leq \frac{D^2}{2\eta_n} + \frac{1}{2} \sum_{t=1}^n \eta_t \|v_t\|_2^2  \label{eq:oga-regret-adaptive}\\
                & \leq 
                \frac{D^2}{2\eta_n} + \frac{G^2}{2} \sum_{t=1}^n \eta_t, \label{eq:oga-regret}
            \end{align}
            where the terms $D$ and $G$ are defined as 
            \begin{align}
                D \defined \sup_{u, u' \in \mc{U}_1} \|u - u'\|, \quad \text{and} \quad G \defined \sup_{v \in \mc{U}_2} \|v\|. \label{eq:oga-D-G}
            \end{align}
            
            Setting $\eta_t = D/(G\sqrt{t})$ in~\eqref{eq:oga-regret} implies a regret upper bound of $ \mc{R}_n \leq \sqrt{n}DG$. Alternatively, with $\eta_t = \frac{D}{\sqrt{M_t}}$ with $M_t = \sum_{i=1}^t \|v_t\|_2^2$, the bound in~\eqref{eq:oga-regret-adaptive} implies a more refined, observation-dependent, regret bound of $\mc{R}_n = \frac{3D}{2} \sqrt{M_n}$. 
    
    \subsection{Discussion of~\Cref{assump:function-class-G}}
    \label{appendix:background-G}
        Recall that the IPM $\dG$ between two probability distributions is defined as 
        \begin{align}
            \dG(P,Q) &= \max_{g \in \mc{G}} \; |\mathbb{E}_P[g(X)] - \mathbb{E}_Q[g(Y)]| \\
            &= \max \lp \max_{g \in \mc{G}} \; \mathbb{E}_P[g(X)] - \mathbb{E}_Q[g(Y)], \;  \max_{g \in \mc{G}} \; \mathbb{E}_P[-g(X)] - \mathbb{E}_Q[-g(Y)]  \rp \\
            &= \max \lp \max_{g \in \mc{G}} \; \mathbb{E}_P[g(X)] - \mathbb{E}_Q[g(Y)], \;  \max_{g^{-} \in \mc{G}^{-}} \; \mathbb{E}_P[g^{-}(X)] - \mathbb{E}_Q[g^{-}(Y)]  \rp, 
        \end{align}
        where $\mc{G}^{-} \defined \{-g: g \in \mc{G}\}$. When,~\Cref{assump:function-class-G} is satisfied, the function classes $\mc{G}$ and $\mc{G}^-$ coincide. However, in general, when $\mc{G} \neq \mc{G}^-$, we need to consider the two possibilities that the witness function $g^*$ lies in $\mc{G}$ or $\mc{G}^-$. We can do this easily by defining the wealth process as the average of two wealth processes with payoff functions lying in $\mc{G}$ and $\mc{G}^-$ respectively. In the betting language, we divide our initial capital of $\$1$, equally among two bettors, who play separate betting games on the same outcomes $\{(X_t, Y_t): t \geq 1\}$, but with payoffs chosen from $\mc{G}$ and $\mc{G}^-$ respectively, using  (possibly) different betting strategies $\prediction^+$  and $\prediction^-$. 
        
        More formally, let $\prediction^+$ and $\prediction^-$ denote two prediction strategies for selecting functions  $\{g_t^+ \in \mc{G}: t \geq 1\}$ and $\{g_t^- \in \mc{G}^-: t \geq 1\}$ respectively, based on the same sequence of observations $\{(X_t, Y_t): t \geq 1\}$. We can define the corresponding wealth processes as follows, with bets $\lambda_t^+$ and $\lambda_t^-$ chosen via the ONS strategy : 
        \begin{align}
            &\wealth_0^{+} = \wealth_0^{-} = \frac{1}{2}, \quad \text{and} \\
            &\wealth_t^+ = \wealth_{t-1}^+ \times \lp 1 + \lambda_t^{+} \lp g_t^+(X_t) - g_t^+(Y_t) \rp \rp, \quad \text{and} \\
            &\wealth_t^- = \wealth_{t-1}^- \times \lp 1 + \lambda_t^{-} \lp g_t^-(X_t) - g_t^-(Y_t) \rp \rp, \quad \text{for } t \geq 1. 
        \end{align}
        Under the null, when $P_X=P_Y$, both these processes are nonnegative martingales with initial value equal to $1/2$. Hence, by adding them, we can get a new process, $\{\wealth_t: t \geq 1\}$ with an initial value equal to $1$. This can be used to define our test, $\tau$, as follows: 
        \begin{align}
            \tau = \min\{n: \wealth_n \geq 1/\alpha\}, \quad \text{where} \quad  \wealth_t = \wealth_t^+ + \wealth_t^-. 
        \end{align}
        Ville's inequality implies that the type-I error of this test is controlled at level $\alpha$. Furthermore, under the alternative when $P\neq Q$, at least one of $\wealth_t^+$ or $\wealth_t^-$ grows to infinity, and we can obtain analogs of~\eqref{eq:main-test-consistency},~\eqref{eq:main-test-exponent} and~\eqref{eq:main-test-stopping-time} of~\Cref{theorem:main-result} depending on the regret behavior of the strategies $\prediction^+$ or $\prediction^-$.

    \section{One-Sample Testing}
    \label{appendix:one-sample}
        We now consider the one-sample testing problem, where we are 
        given a sequence of observations $\{Y_t: t \geq 1\}$, drawn \iid from some unknown distribution $P_X \in \mc{P}(\mc{X})$, and a probability distribution $P_X \in \mc{P}(\mc{X})$; and  the goal is to test 
        \begin{align}
            H_0: P_Y = P_X, \quad \text{versus} \quad H_1: P_Y \neq P_X.   \label{eq:one-sample-def}
        \end{align}
        
        \paragraph*{Oracle One-Sample Test} The general strategy developed for two-sample testing, in~\Cref{sec:general-two-sample}, are also applicable to the problem stated above. In particular, we again select a class of test functions $\mc{G} \subset [-1/2,1/2]^\mc{X}$, which defines an IPM $\dG$ on $\mc{P}(\mc{X})$ as in~\eqref{eq:ipm}.  Let $g^* \equiv g^*(P_X, P_Y, \mc{G})$ represent the witness function in $\mc{G}$ associated with  a pair of distributions, $P_X$ and $P_Y$. Then, we can define the \emph{oracle one-sample} test as 
        \begin{align}
            &\tau^* = \min \{n \geq 1: \wealth^*_n \geq 1/\alpha\} \quad \text{where} \\
            &\wealth^*_n = \wealth^*_{n-1} \times \lp 1 + \lambda^* \lp g^*(Y_n) - \mathbb{E}_{P_X}[g^*(X)] \rp \rp, \; \text{for } n \geq 1, \quad \text{and } \wealth^*_0 = 1. 
        \end{align}
        In the above display, $\lambda^*$ denotes the log-optimal bet value as defined in~\eqref{eq:lambda-star}.  Since both $g^*$ and $\lambda^*$ depend on the unknown distribution $P_Y$, to instantiate a practical sequential test using this approach, we need a prediction strategy and a betting  strategy. 
        
        \paragraph*{Practical One-Sample Test}
            As before, we will fix the betting strategy to the ONS strategy~($\ons$). Let $\prediction$ denote any feasible prediction strategy for selecting the sequence of functions $\{g_t \in \mc{G}: t \geq 1\}$. Then, we can define a practical one-sample test as 
            \begin{align}
                &\tau = \inf \{n \geq 1: \wealth_n \geq 1/\alpha \}, \quad \text{where} \label{eq:one-sample-test-def}\\
                & \wealth_n = \wealth_{n-1} \times \lp 1 + \lambda_n \lp g_n(Y_n) - \mathbb{E}_{P_X}[g_n(X)]\rp \rp, \; \text{for } n \geq 1, \quad \text{and } \wealth_0 = 1.  
            \end{align}
        
            As in the two-sample case, the performance of the test described in~\eqref{eq:one-sample-test-def} depends on the quality of the prediction strategy, $\prediction$, measured via its regret, defined as 
            \begin{align}
                \label{eq:one-sample-regret}
                \mc{R}_n \equiv \mc{R}_n\lp \prediction, \mc{G}, P_X, Y_1^n \rp \defined \lp \sup_{g \in \mc{G} }\sum_{t=1}^n  g(Y_t) - \mathbb{E}_{P_X}[g(X)] \rp  -  \sum_{t=1}^n g_t(Y_t) - \mathbb{E}_{P_X}[g_t(X)]. 
            \end{align}
            
            Depending on the behavior of the regret defined above, we can characterize the statistical properties of the sequential one-sample test of~\eqref{eq:one-sample-test-def}, similar to~\Cref{theorem:main-result}. To state the analogous result for our one-sample test, we need to introduce the term, $\beta_1$, which represents the type-II error exponent for the one-sample test. 
            \begin{align}
            \label{eq:one-sample-exponent}
                \beta_1 = \sup_{\epsilon>0} \inf_{P' \in \mc{P}_{\epsilon, \dG}} \dkl(P', P_Y), \quad \text{where} \quad \mc{P}_{\epsilon, \dG} \defined \{P' \in \mc{P}(\mc{X}):  \dG(P', P_X) \leq \epsilon \}. 
            \end{align}
            We now present the main result of this section. 
            \begin{proposition}
            \label{prop:one-sample}
            Suppose $\dG$ is characteristic~(\Cref{def:characteristic-IPM}) for $\mc{P}_2 = \{ P_Y\in \mc{P}(\mc{X}): P_Y \neq P_X\}$,  and~\Cref{assump:function-class-G} holds. Consider  observations $\{Y_t: t \geq 1\}$ drawn \iid according to $P_Y$. 
            Let $\tau \equiv \tau(\prediction, \ons)$ denote a sequential one-sample test with prediction strategy $\prediction$, and betting strategy $\ons$ introduced in~\Cref{def:ONS}. Then, the following statements are true: 
            \begin{itemize}
                \item For any $\prediction$,  the type-I error rate is controlled at the specified level $\alpha$. That is, 
                \begin{align}
                    \label{eq:one-sample-test-type-I}
                    \mbb{P}_{P_X}\lp \tau < \infty \rp \leq \alpha. 
                \end{align}
                \item Suppose  the per-sequence average regret of $\prediction$ satisfies
                \begin{align}
                    \label{eq:one-sample-no-regret}
                    \limsup_{n \to \infty} \frac{\mc{R}_n\lp \prediction, \mc{G}, P_X, Y_1^n\rp}{n} < \dG(P_X, P_Y), \; \text{almost surely}. 
                \end{align}
                Then, the sequential one-sample test $\tau$ has power one under the alternative, 
                \begin{align}
                    \mbb{P}_{P_Y} \lp \tau < \infty \rp = 1, \quad \text{for each }  P_Y \neq P_X. \label{eq:one-sample-test-consistency}
                \end{align}
                
                \item If there exists a sequence $\{r_n: n \geq 1\}$ such that  $r_n \to 0$ and $\sum_{n \geq 1} \mathbb{P}(E_n^c) < \infty$, then   we have the following upper bound on the expected stopping time, under $H_1$: 
                \begin{align}
                    \mathbb{E}[\tau]  = \mc{O} \lp \nzero +  n_0(\epsilon, \alpha) + \sum_{n \geq 1}  \mathbb{P}(E_n^c) \rp,   
                \end{align}
                where $\Delta$ and~$\sigma$ were defined in~\eqref{eq:Delta-sigma}.  
                \item Suppose  there exists a sequence $\{r_n: n \geq 1\}$ with $r_n \to 0$, and $\mathbb{P}(E_n^c)=0$ for all but finitely many $n \geq 1$. Then, we have the following under $H_1$:
                \begin{align}
                    &\lim_{n \to \infty} -\frac{1}{n} \log \lp \mbb{P}_{P_Y} \lp \tau > n \rp \rp  \geq \beta_1.  
                \end{align}
                Recall that the terms $\nzero$, and $\beta_1$ were defined in~\eqref{eq:n0-def} and~\eqref{eq:one-sample-exponent} respectively. 
            \end{itemize}
        \end{proposition}
        The proof of this result follows exactly along the lines of the proof of~\Cref{theorem:main-result}, and we omit the details.

    \section{Proof of \Cref{theorem:main-result}}
    \label{appendix:proof-main-theorem}   
        Before presenting the details, we first fix some notation. 
        As before, we will use $\Delta$ and $\sigma^2$ to denote $\dG(P_X, P_Y)$ and $\sup_{g \in \mc{G}} \mathbb{V}[g(X) - g(Y)]$ respectively. Furthermore, we also introduce a variant of $\sigma^2$, that we denote by $\sigmatilde^2$, defined as 
        \begin{align}
            \sigmatilde^2 \defined \sup_{g \in \mc{G}} \mathbb{E}\lb \big( g(X)- g(Y) \big)^2 \rb. \label{eq:sigmatilde}
        \end{align}
        It is easy to verify that $\sigmatilde^2 \leq \sigma^2 + \Delta^2$. Another term to be used in proving~\eqref{eq:main-test-stopping-time} is $\gamma^2$, defined as 
        \begin{align}
            \gamma^2 &= \sup_{g \in \mc{G}} \mathbb{E}\lb \lp \lp  g(X) - g(Y) \rp^2 - \mathbb{E}[(g(X')-g(Y'))^2] \rp^2 \rb 
            = \sup_{g \in \mc{G}} \mathbb{V}\lp (g(X)-g(Y))^2\rp. \label{eq:gamma-2}  
        \end{align}
        For any $t \geq 1$, we use $v_t$ to denote $g_t(X_t) - g_t(Y_t)$, and introduce the running sums $S_n = \sum_{t=1}^n v_t$, and $M_n = \sum_{t=1}^n v_t^2$. 
        
        \subsection{Proof of~\eqref{eq:main-test-type-I}}
            This follows as a direct consequence of the fact that, by construction, the process $\{\wealth_t: t \geq 0\}$ is a non-negative martingale with an initial value equal to $1$ under the null. This is because
            \begin{align}
                \mbb{E}[\wealth_t \vert \mc{F}_{t-1}] &= \mbb{E}[\wealth_{t-1} \times \big(1+\lambda_t(g_t(X_t)-g_t(Y_t))\big)\vert \mc{F}_{t-1}]     \\
                & = \wealth_{t-1}\lp 1 + \lambda_t \mbb{E}[g_t(X_t) - g_t(Y_t) | \mc{F}_{t-1}] \rp= \wealth_{t-1}, 
            \end{align}
            where we used the fact that $\lambda_t$ and $g_t$ are $\mc{F}_{t-1}$ measurable, and that under the null, both $X_t$ and $Y_t$ have the same distribution. 
            The result stated in~\eqref{eq:main-test-type-I} then follows by an application of Ville's inequality. 
    
        \subsection{Proof of~\eqref{eq:main-test-consistency}}  
            We note that 
            \begin{align}
                \{\tau = \infty\} = \cap_{t=1}^{\infty} \{\tau > t \} \subset \{\tau > n\} \quad\text{for any } n \geq 1, 
            \end{align}
            which implies
            \begin{align}
                 \mbb{P}\lp \tau = \infty \rp \leq  \mbb{P} \lp \tau > n \rp \; \Rightarrow\;      \mbb{P}\lp \tau = \infty \rp \leq \liminf_{n \to \infty} \mbb{P} \lp \tau > n \rp. \label{eq:proof-main-0}
            \end{align}
            Thus, to show~\eqref{eq:main-test-consistency}, it suffices to show that the $\liminf$ in above display is equal to $0$. 
            To prove that, we first recall that due to the use of ONS betting strategy, the wealth process at time $t$, for any prediction strategy, satisfies 
            \begin{align}
                \wealth_t &\geq \exp \lp \frac{n}{8} \lp \frac{1}{n} \sum_{t=1}^n g_t(X_t) - g_t(Y_t) \rp^2 - \log n\rp \\
                & \geq  \lp \frac{n}{8} \lp \max \left\{ \frac{1}{n} \sum_{t=1}^n g_t(X_t) - g_t(Y_t), \; 0 \right\} \rp^2 - \log n\rp \label{eq:proof-main-1}
            \end{align}
            By definition, the event $\{\tau>n\}$ is contained in the event  $\{\wealth_n<1/\alpha\}$, which, due to~\eqref{eq:proof-main-1}, implies 
            \begin{align}
                \{\tau>n\} &\subset \left \{ \frac{1}{n} \sum_{t=1}^n g_t(X_t) - g_t(Y_t) < \sqrt{ \frac{8 \log(1/\alpha) + \log n}{n}}  \right\} \\
                & \subset \left\{ d_{\mc{G}}(\Phat_{X,n}, \Phat_{Y,n}) - \frac{\mc{R}_n(\prediction,X_1^n, Y_1^n)} {n} < \sqrt{\frac{8 \log(1/\alpha) + \log n}{n}} \right\} \label{eq:proof-main-2}\\
                & \defined F_n. 
            \end{align}
            To complete the proof, it suffices to show that $\indi{F_n} \convas 0$, since it implies 
            \begin{align}
                0 \leq \liminf_{n \to \infty} \mbb{P}\lp \tau > n \rp \leq \lim_{n \to \infty}\mbb{E}[\indi{F_n}] = 0, 
            \end{align}
            where the second inequality uses Fatou's lemma, and the equality follows from an application of the Bounded Convergence Theorem.
            
            We now show that $\indi{F_n} \convas 0$. To see this, first note that by definition, $d_{\mc{G}}(\Phat_{X,n}, \Phat_{Y,n}) = \sup_{g \in \mc{G}} \frac{1}{n} \sum_{t=1}^n g(X_t) - g(Y_t) \geq \frac{1}{n} \sum_{t=1}^n g^*(X_t) - g^*(Y_t) \convas d_{\mc{G}}(P_X, P_Y)>0$.  Here $g^* \equiv g^*(P_X, P_Y, \mc{G})$ denotes the witness function associated with $P_X$ and $P_Y$. Thus, we have 
            \begin{align}
                \dG\lp \Phat_{X, n}, \Phat_{Y, n} \rp - \frac{\mc{R}_n}{n} &\geqas \frac{1}{n} \sum_{t=1}^n g^*(X_t) - g^*(Y_t) \; - \frac{\mc{R}_n}{n}, 
            \end{align}
            which implies that 
            \begin{align}
                \liminf_{n \to \infty} \lp  \dG\lp \Phat_{X, n}, \Phat_{Y, n} \rp - \frac{\mc{R}_n}{n} \rp  & ~\geqas~ \liminf_{n \to \infty} \lp \frac{1}{n} \sum_{t=1}^n g^*(X_t) - g^*(Y_t) \; - \frac{\mc{R}_n}{n} \rp \\
                & ~=~ \Delta - \limsup_{n \to \infty} \frac{\mc{R}_n}{n} > 0. 
            \end{align}
            The last inequality above follows from the no-regret assumption on the prediction strategy. Hence,  the $\liminf$ of the term on the left of $<$ in the definition of $F_n$ in~\eqref{eq:proof-main-2} is positive almost surely. On the other hand, the term on the right converges to $0$, implying that $\indi{F_n} \convas 0$, as required.
    \subsection{Proof of~\eqref{eq:wealth-growth-rate-1}} Using the fact that the regret achieved by the ONS strategy with respect to the best constant bet $\lambda \in [-1/2, 1/2]$ is $\mc{O}(\log t)$; we first get the following lower bound, with the notation $v_t = g_t(X_t) - g_t(Y_t)$: 
    \begin{align}
        \frac{ \log \wealth_n}{n} &\geq 
        \sup_{\lambda \in [-1/2, 1/2]} \frac{1}{n} \sum_{t=1}^n \log(1 + \lambda v_t) - \mc{O}\lp \log n\rp \\ 
        &\geq \frac{ \sum_{t=1}^n v_t}{4} \lp \frac{ \sum_{t=1}^n v_t}{\sum_{t=1}^n v_t^2} \wedge 1 \rp - \mc{O}(\log n/n). 
    \end{align}
    The second inequality is a result of  choosing  a $\lambda \in [-1/2, 1/2]$ that optimizes the lower bound on $\sum_{t=1}^n \log(1+\lambda v_t)$, obtained by using the inequality $\log(1+x) \geq x - x^2$ for $x \geq -1/2$. 
    The result then follows by observing that under the assumption~\eqref{eq:limiting-average-payoff}, we have $ \lim_{n \to \infty} \frac{1}{n} \sum_{t=1}^n v_t \stackrel{(a.s.)}{=} \Delta$, and $\lim_{n \to \infty} \frac{1}{n} \sum_{t=1}^n v_t^2 \stackrel{(a.s.)}{=} \mathbb{E}\lb \big(g^*(X)-g^*(Y)\big)^2\rb$.

        \subsection{Proof of~\eqref{eq:main-test-stopping-time}} 
             Since $\tau$ is a non-negative integer-valued random variable, we have 
            \begin{align}
                \mathbb{E}[\tau] = \sum_{n =0}^{\infty} \mathbb{P}(\tau > n ) \leq \sum_{n= 0}^{\infty} \mathbb{P}\lp \log(\wealth_n) < \log(1/\alpha) \rp, 
            \end{align}
            where the inequality follows from the definition of our test $\tau = \inf \{n \geq 1: \wealth_n \geq 1/\alpha \}$. 
            We now use the regret guarantee of the ONS betting strategy to transform the above probabilities into a more convenient form. In particular we have 
            \begin{align}
                \log \wealth_n &= \sum_{t=1}^n \log \big( 1 + \lambda_t \lp g_t(X_t) - g_t(Y_t) \rp \big) \\
                & \geq \sup_{\lambda \in [-1/2, 1/2]} \sum_{t=1}^n \log \big( 1 + \lambda \lp g_t(X_t) - g_t(Y_t) \rp \big) - 12\log(n).
            \end{align}
            Here we used the fact that the regret of the ONS strategy in our problem is upper bounded by $12\log n$.

For all $n \geq 1$, introduce  the  event 
\begin{align}
   B_n = \left\{ |\sum_{t=1}^n v_t| \leq  \sum_{t=1}^n v_t^2 \right\},   \label{eq:bug-fix-0}
\end{align}
and  write $\{\tau > n\}$ as the union of two disjoint events $ \lp \{\tau>n\} \cap B_n \rp \cup \lp \{\tau>n\} \cap B_n^c \rp$. Let us consider the two cases separately: 
\begin{itemize}
    \item Under $B_n$, we know that $\lambda_0 = \frac{\sum v_t}{2\sum v_t^2} \in [-1/2, 1/2]$. Using the fact that $\log(1 + \lambda v_t) \geq \lambda v_t - \lambda^2 v_t^2$ for $\lambda v_t \geq -0.68$, and setting $\lambda = \lambda_0$, we get that 
    \begin{align}
        \{\tau > n \} \cap B_n &\subset \lbr \frac{1}{4} \lp \sum_{t=1}^n v_t \rp^2/ \lp \sum_{t=1}^n v_t^2\rp  < 12 \log (n/\alpha)\rbr \cap B_n \nonumber \\ 
        &= \lbr  \lp \frac{1}{n}\sum_{t=1}^n v_t \rp^2 <  48 \lp \frac{1}{n}\sum_{t=1}^n v_t^2\rp \frac{\log (n/\alpha)}{n}\rbr \cap B_n \coloneqq D_{n,1}. \label{eq:bug-fix-1}
    \end{align}

    \item We will further break down the event $B_n^c$ into $C_{n,1}$ and $C_{n,2}$ as follows: 
    \begin{align}
        B_n^c = \lbr \sum_t v_t >  \sum_t v_t^2 \rbr \cup \lbr \sum_t v_t < - \sum_t v_t^2 \rbr =: C_{n,1} \cup C_{n,2}.
    \end{align}
    Using the fact that $\log(1+ \lambda x) \geq \lambda x - \lambda^2 x^2$ for $\lambda x \geq -0.68$ again, we get 
    \begin{align}
        \{\tau > n \} \cap C_{n,1} &\subset \lbr \frac{1}{2} \sum_t v_t - \frac{1}{4} \sum_t v_t^2 < 12 \log(n/\alpha) \rbr \cap C_{n,1}  && \text{(set $\lambda=1/2$)} \nonumber \\
        & \subset \lbr \frac{1}{4} \sum_t v_t < 12 \log (n/\alpha) \rbr && \text{(property of $C_{n,1}$)}  \nonumber \\
        & = \lbr \frac{1}{n} \sum_{t=1}^n v_t < \frac{48 \log (n/\alpha)}{n}  \rbr \coloneqq D_{n,2}. \label{eq:bug-fix-2}
    \end{align}

    Similarly, we can consider the remaining portion of $\{\tau > n\}$ as 
    \begin{align}
        \{\tau > n \} \cap C_{n,2} & \subset C_{n,2} = \lbr \frac{1}{n} \sum_{t=1}^n v_t < -\frac{1}{n} \sum_{t=1}^n v_t^2 \rbr  \nonumber \\
        & \subset \lbr  \frac{1}{n} \sum_{t=1}^n v_t \leq 0 \rbr   && \text{\bigg( since $-\sum_{t=1}^n v_t^2 \leq 0$ \bigg)} \nonumber \\ 
        &\coloneqq D_{n,3}. \label{eq:bug-fix-3}
    \end{align}
\end{itemize}

The next step is to upper bound $\mathbb{P}(\tau>n)$ as follows: 
\begin{align}
    \mathbb{P}\lp \tau > n \rp &= \mathbb{P}(\{\tau > n\} \cap B_n) + \mathbb{P}(\{\tau > n\} \cap C_{n,1}) +  \mathbb{P}(\{\tau > n\} \cap C_{n,2}) \nonumber \\
    & \leq \mathbb{P}(D_{n,1}) + \mathbb{P}(D_{n,2}) + \mathbb{P}(D_{n,3}). \label{eq:bug-fix-4}
\end{align}
This implies that the expected stopping time admits the following upper bound: 
\begin{align}
    \mathbb{E}[\tau] \leq  1+ \sum_{n \geq 1} \bigg( \mathbb{P}(D_{n,1}) + \mathbb{P}(D_{n,2}) + \mathbb{P}(D_{n,3}) \bigg). \label{eq:bug-fix-5}
\end{align}
It remains to bound the three sums in the RHS above. To analyze these terms, we need the following two technical lemmas, whose proofs are in~\Cref{appendix:proof-auxiliary}. The first lemma obtains deviation bounds for the processes $\{g^*(X_t) - g^*(Y_t): t \geq 1\}$ and $\{v_t^2: t \geq 1\}$. 
    \begin{lemma}
        \label{lemma:alt-1}
        For all $n \geq 1$, define the event $G_n = G_{n,1} \cap G_{n,2}$, where 
        \begin{align}
            &G_{n,1} = \lbr \frac{1}{n} \sum_{t=1}^n g^*(X_t) - g^*(Y_t) \geq \Delta - \sigmatilde \sqrt{ \frac{4 \log n}{n}} - \frac{2 \log n}{3n} \rbr, \quad \text{and} \\
            &G_{n,2} = \lbr \frac{1}{n} \sum_{t=1}^n v_t^2 \leq \sigmatilde^2 + \gamma \sqrt{ \frac{4 \log n}{n}} + \frac{2 \log n}{3n} \rbr. 
        \end{align}
        Then, we have $\mathbb{P}\lp G_n \rp \geq 1 - 2/n^2$, which in turn, implies that $\sum_{n=1}^\infty \mathbb{P}\lp G_n^c \rp \leq \pi^2/3$. 
    \end{lemma}
The proof of this lemma is in~\Cref{proof:alt-1}. The next lemma records three inequalities between the terms $\sigmatilde^2, \sigma^2$ and $\Delta^2$ introduced at the beginning of~\Cref{appendix:proof-main-theorem}.  
         \begin{lemma}
             \label{lemma:alt-2}
             The following relations are true: 
             \begin{align}
                 \sigmatilde^2 \leq \sigma^2 + \Delta^2, \quad \sigmatilde \leq \sigma + \Delta, \quad \text{and} \quad \sqrt{\sigmatilde} \leq \sqrt{\sigma} + \sqrt{\Delta}. \label{eq:sigmatilde-inequalities} 
             \end{align}
         \end{lemma}
The proof of this lemma is in~\Cref{proof:alt-2}.  We now have the tools to obtain the required upper bounds on three sums in~\eqref{eq:bug-fix-5}. We begin with the terms $\sum_{n \geq 1} \mathbb{P}(D_{n,2})$  and $\sum_{n \geq 1} \mathbb{P}(D_{n,3})$ as they are easier to handle.

\paragraph{Bound on $\boldsymbol{\sum_{n\geq0} \mathbb{P}(D_{n,2})}$.} Recall the event $E_n = \{\mathcal{R}_n/n < r_n\}$. We can write $\mathbb{P}(D_{n,2}) = \mathbb{P}(D_{n,2} \cap E_{n}) + \mathbb{P}(D_{n,2} \cap E_{n^c}) \leq \mathbb{P}(D_{n,2} \cap E_n) + \mathbb{P}(E_n^c)$. We will carry forward the term $\sum_{n \geq 0} \mathbb{P}(E_n^c)$ into the final bound, so for now we only need to focus on the first term. Under the event $E_n$, we know that 
\begin{align}
    \sum_{t=1}^n v_t \geq \sup_{g \in \mathcal{G}}\sum_{t=1}^n 
  g(X_t) - g(Y_t) - n r_n \geq \sum_{t=1}^n g^*(X_t) - g^*(Y_t) - n r_n, 
\end{align}
where $g^*$ is the witness function associated with $P_X$ and $P_Y$. Let us introduce the notation $v^*_t = g^*(X_t) - g^*(Y_t)$. Then, we have 
\begin{align}
    \lbr \frac{1}{n} \sum_{t=1}^n v_t < \frac{48 \log (n/\alpha)}{n} \rbr \cap E_n \subset \lbr \frac{1}{n} \sum_{t=1}^n v^*_t < r_n + \frac{48 \log (n/\alpha)}{n}  \rbr .
\end{align}
Now, observe that 
\begin{align}
    \mathbb{P}\lp \lbr \frac{1}{n} \sum_{t=1}^n v^*_t < r_n + \frac{48 \log (n/\alpha)}{n}  \rbr \rp &\leq \mathbb{P} \lp \lbr \frac{1}{n} \sum_{t=1}^n v^*_t < r_n + \frac{48 \log (n/\alpha)}{n}  \rbr  \cap G_{n,1}\rp + \mathbb{P}(G_{n,1}^c) \\
    & \leq \mathbb{P} \lp \lbr \Delta < \sigmatilde\sqrt{\frac{4 \log n}{n}} + \frac{2 \log n}{3n} + r_n + \frac{48 \log (n/\alpha)}{n}  \rbr  \cap G_{n,1}\rp + \frac{1}{n^2}. 
\end{align}
Here $G_{n,1}$ is the high probability event in~Lemma~\ref{lemma:alt-1}.  
Now, we know from Lemma~\ref{lemma:alt-2} that $\sigmatilde \leq \sigma + \Delta$, and furthermore, assume that $n$ is large enough to ensure that $4 \log n /n \leq 1/4$~(a sufficient condition is if $n \geq 68$). Then, we get 
\begin{align}
    \lbr \Delta < \sigmatilde\sqrt{\frac{4 \log n}{n}} + \frac{2 \log n}{3n} + r_n + \frac{48 \log (n/\alpha)}{n}  \rbr &\subset \lbr \Delta < (\sigma + \Delta)\sqrt{\frac{4 \log n}{n}} + \frac{2 \log n}{3n} + r_n + \frac{48 \log (n/\alpha)}{n}  \rbr \\
    & \subset \lbr \frac{\Delta}{2} < \sigma\sqrt{\frac{4 \log n}{n}} + \frac{2 \log n}{3n} + r_n + \frac{48 \log (n/\alpha)}{n}  \rbr.
\end{align}
The second inclusion uses the fact that for $n \geq 68$, the term $\Delta \sqrt{4 \log n/n} \leq \Delta/2$. 
This allows us to conclude that 
\begin{align}
    &\sum_{n \geq 1} \mathbb{P}\lp D_{n,2} \rp  \leq n_{0,2} + \sum_{n \geq 1} \frac{1}{n^2} + \sum_{n \geq 1} \mathbb{P}(E_n^c), \nonumber \\
    \text{where} \quad 
    & n_{0,2} \coloneqq \inf \lbr n \geq 1: {\frac{\Delta}{2}} \geq \sigma\sqrt{\frac{4 \log n}{n}} + \frac{2 \log n}{3n} + r_n + \frac{48 \log (n/\alpha)}{n} \rbr. \label{eq:n02-def}
\end{align}

\paragraph{Bound on $\boldsymbol{\sum_{n\geq0} \mathbb{P}(D_{n,3})}$.} We again have $\mathbb{P}(D_{n,3}) = \mathbb{P}(D_{n,3} \cap E_{n}) + \mathbb{P}(D_{n,3} \cap E_{n^c}) \leq \mathbb{P}(D_{n,3} \cap E_n) + \mathbb{P}(E_n^c)$, and we only need to analyze the first term. Proceeding as before, we get 
\begin{align}
    \lbr \frac{1}{n} \sum_{t=1}^n v_t < 0 \rbr \cap E_n \subset \lbr \frac{1}{n} \sum_{t=1}^n v^*_t < r_n \rbr.  
\end{align}
Next, we have 
\begin{align}
    \mathbb{P}\lp \lbr \frac{1}{n} \sum_{t=1}^n v^*_t < r_n \rbr\rp & \leq \mathbb{P}\lp \lbr \frac{1}{n} \sum_{t=1}^n v^*_t < r_n \rbr \cap G_{n,1} \rp + \mathbb{P}(G_{n,1}^c) \\
    & \leq \mathbb{P}\lp \lbr \Delta < r_n + {\sigmatilde} \sqrt{\frac{4 \log n}{n}} + \frac{2 \log n}{3n} \rbr \cap G_{n,1} \rp + \frac{1}{n^2} \\
    & {\leq \mathbb{P}\lp \lbr \frac{\Delta}{2} < r_n + \sigma \sqrt{\frac{4 \log n}{n}} + \frac{2 \log n}{3n} \rbr \cap G_{n,1} \rp + \frac{1}{n^2}}. 
\end{align}
As before, the last inequality uses the observation that $\sigmatilde \leq \sigma + \Delta$ and that for $n \geq 68$, we have $\Delta \sqrt{4 \log n/n} \leq \Delta/2$. 
This implies that 
\begin{align}
    &\sum_{n \geq 1} \mathbb{P}\lp D_{n,3} \rp  \leq n_{0,3} + \sum_{n \geq 1} \frac{1}{n^2} + \sum_{n \geq 1} \mathbb{P}(E_n^c), \nonumber \\
    \text{where} \quad 
    & n_{0,3} \coloneqq \inf \lbr n \geq 1: {\frac{\Delta}{2}} \geq r_n + \sigma \sqrt{\frac{4 \log n}{n}} + \frac{2 \log n}{3n} \rbr .\label{eq:n03-def}   
\end{align}

\paragraph{Bound on $\boldsymbol{\sum_{n\geq0} \mathbb{P}(D_{n,1})}$.}  We begin by observing 
\begin{align}
        \lbr  \lp \frac{1}{n}\sum_{t=1}^n v_t \rp^2 <  48 \lp \frac{1}{n}\sum_{t=1}^n v_t^2\rp \frac{\log (n/\alpha)}{n}\rbr  &\subset  \lbr \left \lvert \frac{1}{n} \sum_{t=1}^n v_t \right\rvert < 7 \sqrt{\frac{1}{n}\sum_{t=1}^n v_t^2} \sqrt{\frac{ \log (n/\alpha)}{n}} \rbr  && (\text{since $\sqrt{48} \leq 7$}) \\
        & \subset \lbr \frac{1}{n} \sum_{t=1}^n v_t < 7 \sqrt{\frac{1}{n}\sum_{t=1}^n v_t^2} \sqrt{\frac{ \log (n/\alpha)}{n}} \rbr. 
\end{align}
The last inclusion uses the fact that $\{|A| \leq B\} = \{A \leq B\} \cap \{A \geq -B\} \subset \{A \leq B\}$. 
We again use the event $E_n$, and write 
\begin{align}
    \sup_{n \geq 1} \mathbb{P}(D_{n,1}) \leq \sum_{n \geq 0} \mathbb{P}\lp D_{n,1} \cap E_n \rp + \sum_{n \geq 1} \mathbb{P}(E_n^c). 
\end{align}
Furthermore, we also have 
\begin{align}
    \sum_{n \geq 1}\mathbb{P}(D_{n,1} \cap E_n) &\leq \sum_{n \geq 1}\mathbb{P}(D_{n,1} \cap E_n \cap G_{n,1} \cap G_{n,2}) + \sum_{n \geq 1} \mathbb{P}(G_{n,1}^c) +  \mathbb{P}(G_{n,2}^c) \nonumber \\
    & \leq \sum_{n \geq 1}\mathbb{P}(D_{n,1} \cap E_n \cap G_{n,1} \cap G_{n,2}) + \sum_{n \geq 1} \frac{2}{n^2} \nonumber \\
    & \leq \sum_{n \geq 1} \mathbb{P}\lp \lbr \frac{1}{n} \sum_{t=1}^n v^*_t < r_n + 7 \sqrt{\frac{1}{n}\sum_{t=1}^n v_t^2} \sqrt{\frac{ \log (n/\alpha)}{n}}   \rbr \rp + \frac{\pi^2}{3}.  \label{eq:bug-fix-7}
\end{align}
Under the event $G_{n,2}$~(introduced in Lemma~\ref{lemma:alt-1}), we have the following with probability at least $1-1/n^2$: 
\begin{align}
    \frac{1}{n} \sum_{t=1}^n v_t^2 \leq \widetilde{\sigma}^2 + \gamma \sqrt{\frac{4 \log n}{n}} + \frac{2 \log n}{3 n} ,
\end{align}
where $\gamma^2$ is defined as 
{
\begin{align}
    \gamma^2 = \sup_{g \in \mathcal{G}} \mathbb{V} \lp \big(g(X) - g(Y)\big)^2 \rp \leq \widetilde{\sigma}^2.   \label{eq:gamma-def}
\end{align}
The inequality above can be proved in the following steps: 
\begin{align}
    \gamma^2 &=  {\sup_{g \in \mathcal{G}} \mathbb{V}\lp \big( g(X) - g(Y) \big)^2 \rp} \leq {\sup_{g \in \mathcal{G}} \mathbb{E}\lb \lp g(X) - g(Y) \rp^4 \rb } \leq {\sup_{g \in \mathcal{G}} \mathbb{E}\lb \lp g(X) - g(Y) \rp^2 \rb } = \sigmatilde^2, 
\end{align}
where we have used the fact that $g(X)-g(Y) \in [-1, 1]$, and hence, $|g(X)-g(Y)|^4 \leq |g(X)-g(Y)|^2$.}

Now, using the fact that $\gamma \leq \widetilde{\sigma}$,  we have under $G_{n,2}$: 
\begin{align}
    \frac{1}{n} \sum_{t=1}^n  v_t^2 \leq \widetilde{\sigma}^2 + \sigmatilde \sqrt{\frac{4\log n}{n}} +  \frac{2 \log n}{3n}, \nonumber 
\end{align}
which implies that 
\begin{align}
    7\sqrt{\frac{1}{n} \sum_{t=1}^n v_t^2} \sqrt{\frac{\log (n/\alpha)}{n}} &\leq 7 \sqrt{\sigmatilde^2 + \sigmatilde \sqrt{\frac{4 \log n}{n}} + \frac{2 \log n}{3n} } \times \sqrt{\frac{\log (n/\alpha)}{n}} \\
    &\leq 7 \widetilde{\sigma} \sqrt{\frac{\log (n/\alpha)}{n}} + 7 \sqrt{2\sigmatilde} \lp \frac{\log (n/\alpha)}{n} \rp^{3/4} +7 \sqrt{\frac{2 \log n}{3n} \frac{\log (n/\alpha)}{n}} \nonumber\\
    &\leq 7 \widetilde{\sigma} \sqrt{\frac{2 \log (n/\alpha)}{n}} + 7 \sqrt{2\sigmatilde} \lp \frac{\log (n/\alpha)}{n} \rp^{3/4} + 7 \frac{\log (n/\alpha)}{n}   \nonumber\\ 
    &\leq 7 \sigma \sqrt{\frac{2 \log (n/\alpha)}{n}} + \frac{\Delta}{2} + 7 \sqrt{2\sigmatilde} \lp \frac{\log (n/\alpha)}{n} \rp^{3/4} + 7 \frac{\log (n/\alpha)}{n}   .
    \label{eq:bug-fix-8}
\end{align}
In the second to last inequality, we used the fact that $2/3 \leq 1$ and that $\log(n/\alpha) \geq \log n$, while in the last inequality, we used the fact that $\sigmatilde \leq \sigma + \Delta$, and assumed that $n$ is larger than 
\begin{align}
 n_{0, 4} \coloneqq \inf \{n \geq 1: 7 \sqrt{2 \log(n/\alpha)/n} \leq 1/2 \} = \inf \{n \geq 1: n/\log(n/\alpha) \geq 392\}.    \label{eq:n04-def}
\end{align}

Under the event $G_{n,1}$ we have
\begin{align}
    \frac{1}{n} \sum_{t=1}^n v^*_t \geq \Delta - \sigma \sqrt{\frac{4 \log n}{n}} - \frac{2 \log n}{3 n}.  \label{eq:bug-fix-9}
\end{align}

Plugging~\eqref{eq:bug-fix-8} and~\eqref{eq:bug-fix-9} into~\eqref{eq:bug-fix-7}, we get  the following (using the fact that $2\log n/(3n) \leq \log(n/\alpha)$), with $\widetilde{D}_{n,1} = D_{n,1} \cap E_{n} \cap G_{n,1} \cap G_{n,2}$: 
\begin{align}
    \sum_{n \geq 1} \mathbb{P}(\widetilde{D}_{n,1} ) &\leq \sum_{n \geq 1} \mathbb{P}\lp \lbr  \frac{\Delta}{2} < r_n + 7 \sigma \sqrt{\frac{2 \log (n/\alpha)}{n}} + 7\sqrt{2\sigmatilde}\lp \frac{\log (n/\alpha)}{n} \rp^{3/4} +  8 \frac{\log (n/\alpha)}{n}  + \sigma \sqrt{\frac{4 \log n}{n}}   \rbr \rp \nonumber \\
    & \leq \sum_{n \geq 1}\mathbb{P}\lp \lbr  \frac{\Delta}{2} < r_n + 9 \sigma \sqrt{\frac{2 \log (n/\alpha)}{n}}+ 7\sqrt{2\sigmatilde}\lp \frac{\log (n/\alpha)}{n} \rp^{3/4}  +    \frac{8\log (n/\alpha)}{n}  \rbr \rp,  \nonumber
\end{align}
where the second inequality uses the fact that $\sigma \sqrt{4 \log n /n} \leq 2 \sigma \sqrt{2 \log (n/\alpha)/n}$. This leads us to the required conclusion that 
\begin{align}
    &\sum_{n \geq 1} \mathbb{P}(D_{n,1}) \leq n_{0,1} + \frac{\pi^2}{3} + \sum_{n \geq 1} \mathbb{P}(E_n^c), \nonumber \\
    \text{where} \quad 
    & n_{0,1} \coloneqq \inf \lbr n \geq 1: \frac{\Delta}{2} \geq r_n + 9 \sigma \sqrt{\frac{2 \log (n/\alpha)}{n}} + 7\sqrt{2\sigmatilde}\lp \frac{\log (n/\alpha)}{n} \rp^{3/4}  +    \frac{8\log (n/\alpha)}{n}  \rbr.  \label{eq:n01-def}
\end{align}

\paragraph{Completing the proof.} Finally, we can upper bound the expected stopping time as 
\begin{align}
    \mathbb{E}[\tau] &\leq \sum_{n \geq 1} \mathbb{P}\lp D_{n,1}\rp + \sum_{n \geq 1} \mathbb{P}\lp D_{n,2}\rp + \sum_{n \geq 1} \mathbb{P}\lp D_{n,3}\rp \\
    & = \mathcal{O}\lp \max \{ n_{0,1},\; n_{0, 2},\; n_{0, 3}, n_{0, 4}\} + \sum_{n\geq 1} \mathbb{P}(E_n^c) \rp, 
\end{align}
where $n_{0,1}$, $n_{0,2}$, $n_{0,3}$, and $n_{0,4}$ were introduced in~\eqref{eq:n01-def},~\eqref{eq:n02-def},~\eqref{eq:n03-def}  and~\eqref{eq:n04-def} respectively. From their definitions, it follows that $n_{0,2}, n_{0,3},$ and $n_{0,4}$ are all $\mc{O}(n_0)$, where recall that 
\begin{align}
    n_0 \equiv n_0(\alpha, \Delta, \sigma) = \inf \lbr n \geq 1: r_n + \sigma \sqrt{\frac{\log(n/\alpha)}{n}} + \frac{\log(n/\alpha)}{n} \leq \Delta \rbr. 
\end{align}
To complete the proof we need to show the same for $n_{0,1}$. 
\begin{lemma}
    \label{lemma:n01} The term $n_{0,1}$ defined in~\eqref{eq:n01-def} satisfies the following: 
    \begin{align}
        n_{0,1} \leq \inf \lbr n \geq m_2\vee m_3: \frac{\Delta}{4} \geq r_n + 9 \sigma \sqrt{\frac{2 \log (n/\alpha)}{n}} +  \frac{8\log (n/\alpha)}{n}  \rbr = \mc{O}(n_0),  
    \end{align}
    where 
    \begin{align}
        &m_2 = \inf \lbr n \geq 1: \sigma\sqrt{\frac{\log (n/\alpha)}{n}} \leq \frac{\Delta}{(6272)^{1/3}}  \rbr = \mc{O}(n_0), \quad \text{and} \\
        &m_3 = \inf \lbr n \geq 1: \frac{\log (n/\alpha)}{n} \leq \frac{\Delta}{(6272)^{1/3}}  \rbr = \mc{O}(n_0).
    \end{align}
\end{lemma}
The proof of this statement is in~Section~\ref{proof:n01}.

        \subsection{Proof of~\eqref{eq:main-test-exponent}}
            To show this, we start with~\eqref{eq:proof-main-2}. Since, by assumption $r_n \to 0$ and $\mathbb{P}(\cap_{n=1}^{\infty} E_n)=1$, for every $\epsilon>0$, there exists a finite $N_{\epsilon}$, such that for all $n \geq N_{\epsilon}$, both of the following conditions hold simultaneously: 
            \begin{align}
                r_n &\leq \frac{\epsilon}{2}, \quad \text{and} %
                \sqrt{ \frac{ 8 \log(1/\alpha) + \log n}{n} }  \leq \frac{\epsilon}{2}. \label{eq:proof-main-3}
            \end{align}
            Then, plugging these into~\eqref{eq:proof-main-2} for any $n \geq N_{\epsilon}$, we get that 
            \begin{align}
                \{\tau > n\} \subset \{d_{\mc{G}} (\Phat_{X,n}, \Phat_{Y,n}) \leq \epsilon \} = \{ (\Phat_{X,n}, \Phat_{Y,n}) \in \probclasstwo \}. 
            \end{align}
            The  term $\probclasstwo$ used above was introduced in~\eqref{eq:two-sample-exponent}. 
            Finally, an application of Sanov's theorem implies the following result for any $\epsilon>0$
            \begin{align}
                \liminf_{n \to \infty}\; -\frac{1}{2n} \log \lp \mbb{P}(\tau > n ) \rp & \geq \inf_{(P', Q') \in \probclasstwo} \frac{1}{2} \lp \dkl(P', P_X) + \dkl(Q', P_Y) \rp 
            \end{align}
            as required.

    \section{Proof of~\Cref{prop:bounded-mean}}
    \label{proof:bounded-mean}
        It is easy to verify that under the null, the process $\{\wealth_t: t \geq 1\}$  is a nonnegative martingale with an initial value of $1$. Hence its type-I error is below the level $\alpha$, due to Ville's inequality.  

        Next, to obtain the upper bound on the expected stopping time, we need to get appropriate upper bounds on the terms $\nzero$, and $\sum_{n \geq 1} \mathbb{P}(E_n^c)$. The key to simultaneously bounding both these terms is to select the sequence $\{r_n: n \geq 1\}$ suitably to balance them. 
        Since $D = \sup_{x, x'\in \mc{U}} \|x-x'\|_2 = 1$,  the regret of the OGA strategy with adaptive step-sizes is $\mc{R}_n \leq (3/2) \sqrt{M_n} = (3/2) \sqrt{ \sum_{t=1}^n \|X_t - Y_t\|^2 }$. Using this we choose $r_n$ as follows: 
        \begin{align}
            r_n = \frac{3}{2n } \sqrt{n\sigmatilde^2 + \frac{n^2 \Delta^2}{225}}, \quad \text{and} \quad E_n = \lbr \frac{\mc{R}_n}{n} \leq r_n \rbr. 
        \end{align}
        Recall that~$\sigmatilde$ was defined in~\eqref{eq:sigmatilde}. This choice of $r_n$ allows us to keep $\nzero$ small enough, while also providing enough room to control the deviation of $\mc{R}_n/n$ beyond $r_n$, as we describe below.  

        \paragraph{Bound on $\boldsymbol{\nzero}$.} Plugging the above choice of $r_n$ in~\eqref{eq:n0-def}, we get 
        \begin{align}
            \nzero = \min \lbr n \geq 1: \frac{12}{n} \sqrt{n \sigmatilde^2 + \frac{n^2\Delta^2}{225}} + \frac{\log(n/\alpha)}{n} + \sigma\sqrt{ \frac{ \log (n/\alpha)}{n}} < \Delta\rbr, 
        \end{align}
        which implies that 
        \begin{align}
            \nzero \leq \min \lbr n \geq 1: \frac{12\sigmatilde}{\sqrt{n}}  + \frac{4}{5}\Delta + \frac{\log(n/\alpha)}{n} + \sigma\sqrt{ \frac{ \log (n/\alpha)}{n}} < \Delta\rbr. 
        \end{align}
        In the above display, we used the fact that $\sqrt{x + y} \leq \sqrt{x} + \sqrt{y}$. 
        By an application of~\Cref{lemma:logn-over-n-bound}~(stated at the end of this subsection), we can show that 
        \begin{align}
            \nzero &= \mc{O}\lp \frac{ \sigmatilde^2}{\Delta^2} + \frac{\sigma^2 \log(1/\Delta)}{\Delta^2} + \frac{\log(1/\alpha \Delta)}{\Delta}\rp \\
            &= \mc{O}\lp \frac{ \sigma^2 \log(1/\alpha\Delta)}{\Delta^2}  + \frac{\log(1/\alpha \Delta)}{\Delta} \rp, 
        \end{align}
        where we used the fact that $\frac{\sigmatilde^2}{\Delta^2} \leq \frac{\sigma^2 + \Delta^2}{\Delta^2} = 1 + \frac{\sigma^2}{\Delta^2}$. 

        \paragraph{Bound on $\boldsymbol{\sum_{n \geq 1} \mathbb{P}(E_n^c)}$.} We proceed as follows: 
        \begin{align}
            \mathbb{P}\lp E_n^c \rp &= \mathbb{P}\lp \frac{3}{2n} \sqrt{M_n} > \frac{3}{2n} \sqrt{ n \sigmatilde^2 + \frac{n^2 \Delta^2}{225}} \rp = \mathbb{P}\lp  M_n - \mathbb{E}[M_n] >   \lp n \sigmatilde^2 - \mathbb{E}[M_n] \rp + \frac{n^2 \Delta^2}{225} \rp  \\
            & \leq  \mathbb{P}\lp  M_n - \mathbb{E}[M_n] >  \frac{n^2 \Delta^2}{225} \rp, 
        \end{align}
        where the  inequality uses the fact that $\mathbb{E}[M_n] \leq n \sigmatilde^2$.  Now, by an application of Bernstein's inequality, we know that 
        \begin{align}
            \mathbb{P}\lp  M_n - \mathbb{E}[M_n] >  \gamma \sqrt{\frac{4 \log n}{n}} 
 + \frac{2 \log n}{n} \rp \leq \frac{1}{n^2},  \label{eq:proof-bounded-mean-10}
        \end{align}
        where $\gamma$ denotes the fourth moment term defined in~\eqref{eq:gamma-2}. 
        Now, introduce the term 
        \begin{align}
            n_1 = \min \lbr n \geq 1: \gamma \sqrt{\frac{4 \log n}{n}}  + \frac{2 \log n}{n}  < \frac{n \Delta^2}{225} \rbr, 
        \end{align}
        and note that for all $n \geq n_1$, we have $\mathbb{P}(E_n^c) \leq 1/n^2$, due to~\eqref{eq:proof-bounded-mean-10}. This implies that 
        \begin{align}
            \sum_{n \geq 1} \mathbb{P}\lp E_n^c \rp \leq n_1 + \sum_{n > n_1} \frac{1}{n^2} \leq n_1 + \frac{\pi^2}{6}. 
        \end{align}
        Finally, by two applications of~\Cref{lemma:logn-over-n-bound}, we can show that 
        \begin{align}
            n_1 = \mc{O}\lp \frac{ \log(1/\Delta)}{\Delta} + \frac{\gamma^2 \log(1/\Delta)}{\Delta^2} \rp = \mc{O}\lp \frac{ \log(1/\Delta)}{\Delta} + \frac{\sigma^2 \log(1/\Delta)}{\Delta^2} \rp.  
        \end{align}
        Combining with the bound on $\nzero$, we get the required result.

       We end this section by stating the following lemma, that was used to get the upper bounds on $\nzero$ and $n_1$. 
        \begin{lemma}
            \label{lemma:logn-over-n-bound}
            For $a>0$ and $0<b\leq 4$, introduce the term $N = \min \{n \geq 1: \log(b n)/n \leq a\}$. Then, we have the following: 
            \begin{align}
                N \leq 1 + \max \lbr 20,\; \frac{2 \log(2/a)}{a}\rbr. 
            \end{align}
        \end{lemma}
        The proof of this lemma is in~\Cref{appendix:proof-auxiliary}.

    \section{Proof of~\Cref{prop:kernel-mmd}}
    \label{proof:kernel-mmd}
        The proof of this result follows the same steps as the proof of~\Cref{prop:bounded-mean}.  
        The control of type-I error is due to the fact that the process $\{\wealth_t: t \geq 1\}$ is a nonnegative martingale with an initial value of $1$.
       
        Next, to obtain the upper bound on the expected stopping time, we follow the general strategy introduced in proving~\Cref{prop:bounded-mean}. In particular, we define 
        \begin{align}
            r_n = \frac{3}{2n}\sqrt{n \sigmatilde^2 + \frac{n^2 \Delta^2}{225}}, \quad \text{and} \quad 
            E_n = \lbr \frac{\mc{R}_n}{n} \leq r_n \rbr. 
        \end{align}
        Then, as in the case of the proof of~\Cref{prop:bounded-mean}, we can show that 
        \begin{align}
            &\nzero = \mc{O} \lp \frac{ \sigma^2 \log(1/\alpha \Delta)}{\Delta^2} + \frac{\log(1/\alpha \Delta)}{\Delta}\rp, \quad \text{and} \\ 
            &\sum_{n \geq 1} \mathbb{P}\lp E_n^c \rp = \mc{O} \lp \frac{\sigma^2 \log(1/\Delta)}{\Delta^2} + \frac{\log(1/\Delta)}{\Delta} \rp. 
        \end{align}
        Plugging these two expressions in~\eqref{eq:main-test-stopping-time}, we get the required bound on $\mathbb{E}[\tau]$. 
        
        To get the type-II error exponent, we need to find a sequence $\{r_n: n \geq 1\}$, satisfying the following properties: 
        \begin{itemize}
            \item $r_n \to 0$, and  
            \item $\mathbb{P}(E_n^c)=0$ for all $n \geq 1$. 
        \end{itemize}
        The above two conditions are satisfied by setting $r_n = (3/2)/\sqrt{2n}$, since we know that $\mc{R}_n \leq (3/2)\sqrt{M_n} \leq (3/2)\sqrt{2n}$; which in turn implies that $\mc{R}_n/n \leq (3/2)/\sqrt{2n} =r_n$ with probability $1$.  The result then follows by~\eqref{eq:main-test-exponent}. The specific form of the exponent $\beta$ is due to the fact that the kernel-MMD metric metrizes weak convergence, as noted by~\textcite{zhu2021asymptotically}. 
        
        This completes the proof of~\Cref{prop:kernel-mmd}. \hfill \qedsymbol
        
    \section{Proofs of~\Cref{prop:lower-bound-bounded-mean} and~\Cref{prop:lower-bound-kernel-mmd}}
    \label{proof:lower-bound}
        The first step in the proof of both these results is to obtain an information theoretic lower bound on the expected stopping time of any sequential power-one test $\tau'$ in terms of the KL divergence between the null and alternative distributions. We state and prove a general version of this result first, and then  use it to prove~\Cref{prop:lower-bound-bounded-mean} and~\Cref{prop:lower-bound-kernel-mmd} in~\Cref{proof:lower-bound-bounded-mean} and~\Cref{proof:lower-bound-kernel-mmd} respectively. 
        \begin{lemma}
            \label{lemma:lower-bound-general} Given a stream of observations $Z_1, Z_2, \ldots \stackrel{\iid}{\sim} P_Z$, consider the problem of testing the null $H_0: P_Z = Q_0$ against the alternative $H_1: P_Z = Q_1$, for $Q_0 \neq Q_1$. Let $\tau'$ denote any level-$\alpha$, power-one, sequential test for this problem. Then, we have 
            \begin{align}
                \mathbb{E}_{H_1}[\tau'] \geq \frac{\log(1/\alpha)}{\dkl(Q_1, Q_0)}. 
            \end{align}
        \end{lemma}
        \begin{remark}
            To apply this result to the bounded-mean testing and two-sample testing problems, we will construct distributions $P_X$ and $P_Y$, and set $Q_0 = P_X \times P_X$, and $Q_1 = P_X \times P_Y$. Then, for any level-$\alpha$ sequential test for these problems,~\Cref{lemma:lower-bound-general} immediately implies that $\mathbb{E}_{H_1}[\tau'] \geq \log(1/\alpha)/\dkl(P_Y, P_X)$. With an appropriate choice of the distributions $P_X$ and $P_Y$, we get lower bound expressions matching the upper bounds in~\Cref{prop:bounded-mean} and~\Cref{prop:kernel-mmd}.
        \end{remark}
        \begin{proof}
         The proof of this statement uses some standard properties of the KL-divergence that are often used in proving lower bounds for multi-armed bandit algorithms, as we detail below.  

        Since $\tau'$ is a power-one test with finite expectation under $H_1$, we can use Wald's Lemma to get 
        \begin{align}
            \dkl\lp Q_1^{\otimes \tau'}, Q_0^{\otimes \tau'} \rp &= \mathbb{E}_{H_1} \lb \log \lp \frac{dQ_1^{\otimes \tau'}}{dQ_0^{\otimes \tau'}} \rp \rb
             = \mathbb{E}_{H_1} \lb \sum_{t=1}^{\tau'} \log \lp \frac{dQ_1}{dQ_0} \rp \rb \\
             & = \mathbb{E}_{H_1}[\tau'] \dkl(Q_1, Q_0). \label{eq:lower-general-1}
        \end{align}
        Note that the last equality used the assumptions that $\tau'$ has finite expectation under the alternative, and that $0<\dkl(Q_1, Q_0)<\infty$. 

        Next, we consider the event $E = \{\tau' < \infty\}$. By definition of a sequential level-$\alpha$ test of power one, we have 
        \begin{align}
            q_1 \defined \mathbb{P}_{H_1}\lp E \rp = 1, \quad \text{and} \quad 
            q_0 \defined \mathbb{P}_{H_0} \lp E \rp \leq \alpha. 
        \end{align}

        Next, on evaluating KL divergence between two Bernoulli distributions with parameters $q_1$ and $q_0$, we get
        \begin{align}
            \dkl(q_1, q_0) &= q_1 \log \lp \frac{q_1}{q_0} \rp + (1-q_1) \log \lp \frac{1-q_1}{1-q_0} \rp \\
            & = \log \lp \frac{1}{q_0} \rp \geq \log \lp \frac{1}{\alpha} \rp. \label{eq:lower-general-2}
        \end{align}

        The final step is to relate~\eqref{eq:lower-general-1} and~\eqref{eq:lower-general-2} by using the data-processing inequality, as follows 
        \begin{align}
            \dkl(Q_1^{\otimes \tau'}, Q_0^{\otimes \tau'}) = \mathbb{E}_{H_1}[\tau'] \dkl(Q_1, Q_0) \geq \dkl(q_1, q_0),  
        \end{align}
        which implies the required statement
        \begin{align}
            \mathbb{E}_{H_1}[\tau'] \geq \frac{\log(1/\alpha)}{\dkl(Q_1, Q_0)}.
        \end{align}
        This completes the proof of~\Cref{lemma:lower-bound-general}. 
        \end{proof}

        \subsection{Proof of~\Cref{prop:lower-bound-bounded-mean}}        
        \label{proof:lower-bound-bounded-mean}
            The proof of this result follows directly by applying~\Cref{lemma:lower-bound-general} to the construction used by~\textcite[Lemma H.8]{lu2021variance} for deriving variance dependent lower bounds for best arm identification in multi-armed bandits. We present the details below for completeness. 
    
            In what follows, we show that there exist problem instances for which any sequential power-one test with finite mean must have $\mathbb{E}[\tau'] = \Omega (s^2/r^2)$. Then, to complete the proof, we show that for these instances we have $r = \Theta(\Delta)$ and $s^2 \geq \sigma^2$, for $\Delta$ and $\sigma$ defined in~\eqref{eq:Delta-sigma}. 
    
            We consider two cases, depending on the relative values of $r$ and $s^2$. 
            \paragraph{Case 1: $\boldsymbol{s^2 \geq 5r}$.}
            Introduce the following distributions for bounded random variables in $[0,1]$:
            \begin{align}
                P_X = \begin{cases}
                    0.5 + s, & w.p.\  0.5 - r/s \\
                    0.5 - s, & w.p. \ 0.5 + r/s
                \end{cases}, 
                \quad  \text{and} \quad 
                P_Y = \begin{cases}
                    0.5 + s, & w.p. \  0.5 \\
                    0.5 - s, & w.p. \ 0.5
                \end{cases}. 
            \end{align}
            Then, consider a bounded mean testing problem under the null with distribution $Q_0 = P_X \times P_X$, and another problem under the alternative with $Q_1 = P_X \times P_Y$. Then, by an application of~\Cref{lemma:lower-bound-general}, we get that for any $\tau'$ with finite expectation under $H_1$, we have 
            \begin{align}
                \mathbb{E}[\tau'] \geq \frac{\log(1/\alpha)}{\dkl(P_Y, P_X)} = \frac{2 \log(1/\alpha)}{\log(1- 4r^2/s^2)} = \Omega \lp  \frac{\log(1/\alpha)s^2}{r^2} \rp,  
            \end{align}
            where we used the bound $\log(1+x) \leq x$ and $1/(1-x) \geq x/2$ for $0<x<0.5$. The condition $s^2 \geq r$ was invoked to ensure that $r^2/s^2 \leq r/s^2$ is small enough.

            Note that in this case, the mean of $P_X$ is $\mu_X = 0.5 - {2r}$, which implies that $\Delta = |\mu_X - \mu_Y| = |0.5 {-2r} - 0.5| = {2r}$. To analyze $\sigma^2$, we note that by definition we have $\sigma^2 = \mathbb{V}(X) + \mathbb{V}(Y) = 2 s^2$. Together, these two results imply that $\mathbb{E}_{H_1}[\tau'] = \Omega\lp \log(1/\alpha) \sigma^2 / \Delta^2\rp$.

            \paragraph{Case 1: $\boldsymbol{s^2 < 5r}$.} In this case, we use a different distribution  {with $r \in (0, 1/8]$}: 
            \begin{align}
                P_X' = \begin{cases}
                    0.5 + s, & w.p. \ 0.5 - 2r \\
                    0.5 - s, & w.p. \ 0.5 - 2r \\
                    0, & w.p. \ 4r \\
                \end{cases}. 
            \end{align}
            Repeating the same argument with $Q_0 = P_X' \times P_X'$ and $Q_1 = P_X' \times P_Y$, we get that 
            \begin{align}
                \mathbb{E}[\tau'] \geq \frac{ \log(1/\alpha)}{\dkl(P_Y, P_X')} = {\frac{\log(1/\alpha)}{\log(1/(1-4r))} = \frac{\log(1/\alpha)}{8r}}, 
            \end{align}
            {where the last equality uses the fact that $4r \leq 1/2$ to get $\log(1/(1-4r)) \leq 4r/(1-4r) \leq 8r$.}
    
            As in the previous case, we can verify that $\mu_{X'} = (0.5-2r)$, which implies that $\Delta = |0.5 - 2r - 0.5| = 2r$. Similarly, we have $\sigma^2 = \mathbb{V}(X') + \mathbb{V}(Y) \leq {2s^2 +r}$, hence in this case we have $\mathbb{E}[\tau'] = \Omega \lp \log(1/\alpha) /\Delta\rp$. 

            {Combining the two cases, we can conclude that $\mathbb{E}[\tau'] = \Omega \lp \log(1/\alpha) \lp \frac{\sigma^2}{\Delta^2} + \frac{1}{\Delta} \rp \rp$.}
        
        \subsection{Proof of~\Cref{prop:lower-bound-kernel-mmd}}        
        \label{proof:lower-bound-kernel-mmd}
            To prove this statement, we set $\mc{X} = \mathbb{R}$, select the Gaussian kernel $k(x,y) = \exp(-(x-y)^2/b^2)$ for some fixed $b>0$, and consider the following two distributions: 
            \begin{align}
                P_X = N(0, s^2), \quad \text{and} \quad P_Y = N(r, s^2), 
            \end{align}
            for some $r>0$ and $s^2<1$. Using these two distributions,  we apply~\Cref{lemma:lower-bound-general}, with $Q_1 = P_X \times P_Y$ and $Q_0 = P_X \times P_X$, to get the following for any level-$\alpha$, power-one sequential test $\tau'$ with finite expectation: 
            \begin{align}
                \mathbb{E}[\tau'] \geq \frac{s^2\log(1/\alpha)}{r^2}. \label{eq:mmd-lower-1}
            \end{align}
            For the Gaussian kernel,~\textcite[Lemma 1]{reddi15high} showed that the MMD between distributions $P_X$ and $P_Y$ satisfies 
            \begin{align}
                \Delta^2 = \Omega \lp  \frac{2 r^2}{b^2} \rp. \label{eq:mmd-lower-2}
            \end{align}
    
            Now, observe that the term $\sigma^2$ can be written as 
            \begin{align}
                \sigma^2 = \sup_{g \in \mc{G}} \mathbb{V}_{P_X \times P_Y}\lp g(X) - g(Y) \rp   = \sup_{g \in \mc{G}}  \mathbb{V}_{P_X}\lp g(X) \rp  + \mathbb{V}_{P_Y}\lp g(Y) \rp. 
            \end{align}
            Recall that $\mc{G}$ in this case is the unit norm ball in the RKHS associated with the kernel $k$. Now, we upper bound $\mathbb{V}_{P_X}(g(X))$ by using the reproducing property of $\mc{G}$, to get 
            \begin{align}
                \mathbb{V}_{P_X}\lp g(X) \rp &= \mathbb{V}_{P_X}\lp \langle k(X, \cdot) - \mu_X, g \rangle \rp  \leq \mathbb{E}_{P_X}\lb \|g\|^2 \|k(X, \cdot) - \mu_X \|^2 \rb \\
                & \leq \mathbb{E}_{P_X}\lb \|k(X, \cdot) - \mu_X \|^2 \rb 
                {= \mathbb{E}_{P_X}[k(X, X)] - \|\mu_X\|^2 }\\
                & = 1 - \mathbb{E}_{P_X\times P_{X'}}[k(X, X')]. 
            \end{align}
            {In the last equality we have used the fact that $k(X, X)=1$, and that $\|\mu_X\|^2 = \langle \mathbb{E}_{P_X}[k(X, \cdot)], \mathbb{E}_{P_{X'}}[k(X', \cdot)] \rangle = \mathbb{E}_{P_X \times P_{X'}}[k(X, X')]$.}
            To get a more explicit upper bound, we use the fact that $1 - e^{-x} \leq x$ for all $x \geq 0$ to get 
            \begin{align}
                \mathbb{E}[1 - k(X, X')] = \mathbb{E}\lb 1 - \exp\lp - \frac{|X-X'|^2}{b^2} \rp \rb \leq \mathbb{E}\lb \frac{|X-X'|^2}{b^2}\rb = \frac{2 s^2}{b^2}. 
            \end{align}
            The last equality uses the fact that $X, X'$ are \iid $N(0, s^2)$. 
            Thus, we have proved that 
            \begin{align}
                \mathbb{V}_{P_X}\lp g(X) \rp &\leq  \frac{2s^2}{b^2}. 
            \end{align}
            Following the same steps to evaluate $\mathbb{V}_{P_Y}\lp g(Y) \rp$, we conclude that 
            \begin{align}
                \sigma^2 \leq \sup_{g \in \mc{G}} \lp \frac{2 s^2}{b^2} + \frac{2 s^2}{b^2}  \rp =  \frac{4 s^2}{b^2}.  \label{eq:mmd-lower-3}
            \end{align}
            Combining~\eqref{eq:mmd-lower-1}, \eqref{eq:mmd-lower-2}, and \eqref{eq:mmd-lower-3}, we get the required result that $\mathbb{E}[\tau'] = \Omega \lp \log(1/\alpha) \sigma^2/\Delta^2\rp$.

    \section{Proof of~\Cref{prop:time-varying}}    
    \label{proof:time-varying}
    
        The control of type-I error again follows  due to the fact that the wealth process is a test martingale under the null, and hence it crosses the level $1/\alpha$ with probability smaller than $\alpha$.  
        
        To show that the test has power one, our starting point is~\eqref{eq:proof-main-1}. Note that to prove that the test is consistent,  it suffices to show  
        \begin{align}
            \liminf_{n \to \infty} \lv \frac{1}{n} \sum_{t=1}^n g_t(X_t) - g_t(Y_t) \rv \stackrel{a.s.}{>} 0. 
        \end{align}
        
        We proceed as follows: 
        \begin{align}
            \liminf_{n \to \infty} \frac{1}{n} \sum_{t=1}^n g_t(X_t) - g_t(Y_t) & = \liminf_{n \to \infty} \lp -\frac{\mc{R}_n}{n} + \sup_{g \in \mc{G}} \frac{1}{n} \sum_{t=1}^n g(X_t) - g(Y_t) \rp  \\
            & \geq \liminf_{n \to \infty} - \frac{\mc{R}_n}{n} +\liminf_{n \to \infty} \sup_{ g \in \mc{G}} \frac{1}{n} \sum_{t=1}^n g(X_t) - g(Y_t) \\
            & \geq -\limsup_{n \to \infty} \frac{\mc{R}_n}{n} + \sup_{g \in \mc{G}} \liminf_{n \to \infty}  \frac{1}{n} \sum_{t=1}^n g(X_t) - g(Y_t) \\
            & \defined -\limsup_{n \to \infty} \frac{\mc{R}_n}{n} + \sup_{g \in \mc{G}} A(g).
        \end{align}
        By assumption on the prediction scheme, we know that $-\limsup_{n \to \infty} \mc{R}_n / n > - \dG(\boldsymbol{P}, \boldsymbol{Q})$. Hence, to complete the proof, it suffices to show that $\sup_{g \in \mc{G}} A(g) \geq \dG(\boldsymbol{P}, \boldsymbol{Q})$.

        To analyze the term $A(g)$ for some arbitrary $g \in \mc{G}$, we  first introduce some notation: 
        \begin{align}
            a_t(g) &\defined g(X_t) - g(Y_t), \quad b_t(g) = \mathbb{E}[a_t(g)|\mc{F}_{t-1}] \\
            A_n(g) &= \frac{1}{n} \sum_{t=1}^n a_t(g), \quad B_n(g) = \frac{1}{n}\sum_{t=1}^n b_t(g). 
        \end{align}
        Note that with this notation, we have $A(g) = \liminf_{n \to \infty} A_n(g)$. We can then proceed as follows: 
        \begin{align}
            A(g) &= \liminf_{n \to \infty} A_n(g) = \liminf_{n \to \infty} \big(B_n(g) + (A_n(g) - B_n(g)) \big) \\
            & \geq \liminf_{n \to \infty} B_n(g) - \limsup_{n \to \infty} |A_n(g) - B_n(g)| \\
            & = \liminf_{n \to \infty} B_n(g) = \liminf_{n \to \infty} \frac{1}{n} \sum_{t=1}^n \mathbb{E}[g(X_t) - g(Y_t)|\mc{F}_{t-1}].  \label{eq:proof-time-varying-0}
        \end{align} 
        In the above display,~\eqref{eq:proof-time-varying-0} follows due to the fact that  $\lim_{n \to \infty}|A_n(g) - B_n(g)| = 0$ almost surely, which we prove next. For any fixed $n$, define the event $E_n = \{|A_n(g) - B_n(g)|<\sqrt{8 \log(2n^2)/n} \}$, and note that, by Azuma's inequality, we have that $\mathbb{P}(E_n^c) \leq 1/n^2$. Hence, by the Borel-Cantelli lemma, we know that $\mathbb{P}\lp E_n^c\; i.o. \rp = 0$, which implies that the event $E_n$ occurs for all but finitely many $n$~(almost surely). Since $0 \leq |A_n(g) - B_n(g)|\leq \sqrt{8 \log(2n^2)/n}$ under the event $E_n$, the result follows as the upper bound converges to $0$ as $n\to \infty$.  
        
        Thus, we have proved that for any fixed $g \in \mc{G}$, we have 
        \begin{align}
            \liminf_{n \to \infty} \frac{1}{n} \sum_{t=1}^n g_t(X_t) - g_t(Y_t) & \geq \liminf_{n \to \infty} \frac{1}{n} \sum_{t=1}^n \mathbb{E}[g(X_t)-g(Y_t)|\mc{F}_{t-1}]. 
        \end{align}
        Since $g \in \mc{G}$ in the above inequality is arbitrary, and the left side is independent of $\mc{G}$, we can conclude the following by taking a supremum over all $g \in \mc{G}$: 
        \begin{align}
            \liminf_{n \to \infty} \frac{1}{n} \sum_{t=1}^n g_t(X_t) - g_t(Y_t) & ~\geqas~ \sup_{g \in \mc{G}} \liminf_{n \to \infty} \frac{1}{n} \sum_{t=1}^n \mathbb{E}[g(X_t)-g(Y_t)|\mc{F}_{t-1}] \\
            & ~=~ \dG\lp \boldsymbol{P}, \boldsymbol{Q} \rp. 
        \end{align}
        The required result then follows by combining the above with the observation that  for any sequence of observations $\{(X_t, Y_t): t \geq 1\}$, we have 
        \begin{align}
            \liminf_{n \to \infty} \lv \frac{1}{n} \sum_{t=1}^n g_t(X_t) - g_t(Y_t) \rv \geq             \liminf_{n \to \infty}  \frac{1}{n} \sum_{t=1}^n g_t(X_t) - g_t(Y_t). 
        \end{align}
        This completes the proof. \hfill \qedsymbol
               
    \section{Proof of~\Cref{theorem:abstract-power-result}}
    \label{proof:abstract-power-result}
        \paragraph{Consistency against a fixed alternative} 
            The plug-in or ERM prediction         strategy~$\Aerm$ selects the function $g_t$ at time $t \geq 2$ as follows: 
            \begin{align}
                \gtilde_t \in \argmax_{\gtilde \in \Gtilde} \frac{1}{t-1} \sum_{i=1}^{t-1} \gtilde(Z_i).  \label{eq:erm-prediction}
            \end{align}
            Recall that the notations $\gtilde$ and $\Gtilde$ were introduced in~\eqref{eq:G-tilde}.
            
            Consider any fixed alternative distribution $P \in \altclass$, and introduce $\gtilde^* \equiv \gtilde^*(P) = \argmax_{\gtilde \in \Gtilde} \mathbb{E}_P[\gtilde(Z)]$, to denote the witness function associated it. For simplicity, we assume that there exists a $\gtilde^*$ achieving the maximum. In case, this is not true, we can repeat the ensuing argument with some $\gtilde^*$ whose expected value is within $\delta$ of the supremum, and then take $\delta \to 0$. 
            Under the assumption that $\lim_{n \to \infty} C_n(\Gtilde, P) \to 0$, using~\textcite[Theorem~26.5]{shalev2014understanding},   we can identify a sequence of events $\{E_n: n \geq 1\}$ with $\mathbb{P}_P(E_n) \geq 1-1/n^2$ for all $n \geq 1$, defined as follows: 
            \begin{align}
                E_n \defined \left\{  \mathbb{E}[\gtilde^*(Z)] -\mathbb{E}[\gtilde_n(Z_t)|\mc{F}_{t-1}]   \leq \epsilon_n  \right\}, \quad \text{where} \quad {\epsilon_n = 2 C_n(\Gtilde, P) + 5\sqrt{\frac{2 \log n}{n}} \;\stackrel{n \to \infty}{\longrightarrow}\; 0}. \label{eq:proof-abstract-1}
            \end{align}
            Since $\sum_{n=1}^{\infty}1/n^2 < \infty$, we can conclude that $\mathbb{P}\lp\{ E_n^c \text{ infinitely often} \} \rp = 0$ by Borel-Cantelli Lemma.
            Next, we observe that 
            \begin{align}
                \frac{1}{n} \sum_{t=1}^n \gtilde_t(Z_t) =  \frac{1}{n} \sum_{t=1}^n \gtilde_t(Z_t) - \mathbb{E}[\gtilde_t(Z_t)|\mc{F}_{t-1}] +  \frac{1}{n} \sum_{t=1}^n \mathbb{E}[\gtilde_t(Z_t)|\mc{F}_{t-1}], 
            \end{align}
            which implies that 
            \begin{align}
                \liminf_{n \to \infty} \frac{1}{n} \sum_{t=1}^n \gtilde_t(Z_t) &=  \liminf_{n \to \infty} \frac{1}{n} \sum_{t=1}^n \gtilde_t(Z_t) - \mathbb{E}[\gtilde_t(Z_t)|\mc{F}_{t-1}] + \liminf_{n \to \infty}  \frac{1}{n} \sum_{t=1}^n \mathbb{E}[\gtilde_t(Z_t)|\mc{F}_{t-1}]  \\
                & \eqas 0 + \liminf_{n \to \infty}  \frac{1}{n} \sum_{t=1}^n \mathbb{E}[\gtilde_t(Z_t)|\mc{F}_{t-1}] \label{eq:proof-abstract-2} \\
                & \eqas  \liminf_{n \to \infty}  \frac{1}{n} \sum_{t=1}^n \mathbb{E}[\gtilde_t(Z_t)|\mc{F}_{t-1}] \indi{E_t} \label{eq:proof-abstract-3} \\ 
                & \geqas  \liminf_{n \to \infty}  \frac{1}{n} \sum_{t=1}^n \lp  \mathbb{E}[\gtilde^*(Z_t)] - \epsilon_t \rp\indi{E_t} \label{eq:proof-abstract-4} \\ 
                & \geqas  \liminf_{n \to \infty}  \frac{1}{n} \sum_{t=1}^n   \mathbb{E}[\gtilde^*(Z_t)] - \limsup_{n \to \infty} \frac{1}{n} \sum_{t=1}^n\epsilon_t  \label{eq:proof-abstract-5} \\  
                &=  \mathbb{E}[\gtilde^*(Z)] > 0.\label{eq:proof-abstract-6}
            \end{align}
            In the above display: \\
            \eqref{eq:proof-abstract-2} uses the strong law of large numbers~(SLLN) for bounded martingale difference sequences, \\
            \eqref{eq:proof-abstract-3} and~\eqref{eq:proof-abstract-5} use the fact that $\mathbb{P}\lp \{E_n^c \text{ i.o.}\}\rp =0$, \\
            \eqref{eq:proof-abstract-4} uses the definition of the event $E_n$, and \\
            \eqref{eq:proof-abstract-6} uses the fact that $\mathbb{E}_P[\gtilde^*(Z)]>0$ when $P \in \altclass$ as $\dG$ is \emph{characteristic} for the class of distributions defined in~\eqref{eq:general-pair-of-distributions}, and that $\lim_{n \to \infty} \epsilon_n = 0$ under the assumption that $\lim_{n \to \infty} C_n(\Gtilde, P) = 0$. 
            
            Finally, the result follows by noting that $\liminf_{n \to \infty}\lv \frac{1}{n} \sum_{t=1}^n \gtilde_t(Z_t) \rv \geq \liminf_{n \to \infty} \frac{1}{n} \sum_{t=1}^n \gtilde_t(Z_t)$. This completes the proof. \qed

        \paragraph{Minimum detectable separation}  We now consider the case when the stronger assumption on the complexity of $\Gtilde$ holds, that is, 
        \begin{align}
            C_n(\Gtilde) \defined \sup_{P \in \mc{P}(Z)} C_n(\Gtilde, P) \stackrel{n \to \infty}{\longrightarrow} 0. 
        \end{align}
        Now, for some $n \geq 1$, consider any $P \in \altclass(\Delta_n)$ and observe that the event that $\tau$ doesn't reject the null in the first $n$ rounds satisfies
        \begin{align}
            \{\tau > n\} &\subset  \lbr \log \wealth_n \leq \log (1/\alpha) \rbr \\
            & \subset \lbr \lv \frac{1}{n} \sum_{t=1}^n \gtilde_t(Z_t) \rv \leq \sqrt{\frac{8 \log(n/\alpha)}{n}} \rbr. \label{eq:proof-abstract-7}
        \end{align}
        The second inclusion in the above display uses the lower bound on the wealth process for ONS strategy as stated in~\eqref{appendix:background-ONS}. 
        Observe that 
        \begin{align}
            \lv \frac{1}{n} \sum_{t=1}^n 
            \gtilde_t(Z_t)\rv & = 
            \lv  \lp \frac{1}{n} \sum_{t=1}^n \gtilde_t(Z_t) - \mathbb{E}[\gtilde_t(Z_t)|\mc{F}_{t-1}] + \mathbb{E}[\gtilde_t(Z_t)|\mc{F}_{t-1}] - 
            \mathbb{E}[\gtilde^*(Z_t)] +
            \mathbb{E}[\gtilde^*(Z_t)]  \rp 
            \rv \\
            & \geq  \mathbb{E}[\gtilde^*(Z)] - 
            \lv   \frac{1}{n} \sum_{t=1}^n \gtilde_t(Z_t) - \mathbb{E}[\gtilde_t(Z_t)|\mc{F}_{t-1}] \rv   - \lv \frac{1}{n} \sum_{t=1}^n \mathbb{E}[\gtilde_t(Z_t)|\mc{F}_{t-1}] -  
            \mathbb{E}[\gtilde^*(Z_t)]  \rv \label{eq:proof-abstract-8}\\
            & \defined \mathbb{E}[\gtilde^*(Z)] - U_n - V_n. \label{eq:proof-abstract-9}
        \end{align}
        Combining~\eqref{eq:proof-abstract-7} and~\eqref{eq:proof-abstract-9}, we get that 
        \begin{align}
            \{\tau > n\} \subset \lbr \mathbb{E}[\gtilde^*(Z)] \leq \sqrt{\frac{8 \log(n/\alpha)}{n}} + U_n + V_n \rbr.\label{eq:proof-abstract-10}
        \end{align}
        We now introduce the following two events, which we will show later in~\Cref{lemma:abstract-1} occur with high probability. 
        \begin{align}
            E_{n,1} &\defined \{U_n \leq u_n\}, \quad u_n = \sqrt{ \frac{ 8 \log(4/\gamma)}{n}}, \quad \text{and} \label{eq:proof-abstract-11} \\
            E_{n,2} & \defined \{V_n \leq v_n \}, \quad v_n = \frac{1}{n}\lp 2 + \sum_{t=1}^{n-1} \lp 2C_t(\Gtilde) + 5 \sqrt{ \frac{2 \log(16n /\gamma)}{t}} \rp \rp. \label{eq:proof-abstract-12}
        \end{align}

        Next, starting with~\eqref{eq:proof-abstract-10}, we observe the following: 
        \begin{align}
            \{\tau > n \} \cap \lp E_{n,1} \cap E_{n,2} \rp &\subset \lbr \mathbb{E}[\gtilde^*(Z)] \leq \sqrt{\frac{8 \log(n/\alpha)}{n}} + U_n + V_n \rbr \cap \lp E_{n,1} \cap E_{n,2} \rp \\
            &\subset \lbr \mathbb{E}[\gtilde^*(Z)] \leq \sqrt{\frac{8 \log(n/\alpha)}{n}} + u_n + v_n \rbr \label{eq:proof-abstract-13} \\
            &= \lbr \mathbb{E}[\gtilde^*(Z)] \leq \Delta_n^* \rbr  \label{eq:proof-abstract-14}\\
            &\subset \lbr \Delta_n \leq \mathbb{E}[\gtilde^*(Z)] \leq \Delta_n^* \rbr  = \emptyset.  \label{eq:proof-abstract-15}
        \end{align}
        In the above display:\\
        \eqref{eq:proof-abstract-13} uses the fact that $U_n \leq u_n$ and $V_n \leq v_n$ under $(E_{n,1} \cap E_{n,2})$. \\
        \eqref{eq:proof-abstract-14} uses the fact that $\Delta_n^* = \sqrt{\frac{8 \log(n/\alpha)}{n}} + u_n + v_n$ by definition. \\
        \eqref{eq:proof-abstract-15} uses the fact that $\mathbb{E}_P[\gtilde^*(Z)]\geq \Delta_n$ as $P \in \altclass(\Delta_n)$, and the fact that $\Delta_n > \Delta_n^*$ by assumption. 
        
        Hence, from~\eqref{eq:proof-abstract-15}, we can conclude that 
        \begin{align}
            \mathbb{P}_P(\tau > n ) \leq \mathbb{P}_P\lp (E_{n,1} \cap E_{n,2})^c \rp \leq \mathbb{P}_P\lp E_{n,1}^c \rp + \mathbb{P}_P\lp E_{n,2}^c \rp. 
        \end{align}
        This implies that to complete the proof, it suffices to show that $\mathbb{P}_P(E_{n,i}^c) \leq \gamma/2$ for $i=1, 2$. We show this below in~\Cref{lemma:abstract-1}. 
        
        \begin{lemma}
            \label{lemma:abstract-1}
            For the events $E_{n,1}$ and $E_{n,2}$ defined in~\eqref{eq:proof-abstract-11} and~\eqref{eq:proof-abstract-12} respectively, we have $\mathbb{P}(E_{n,1}^c) \leq \frac{\gamma}{2}$ and $\mathbb{P}(E_{n,2}^c) \leq \frac{\gamma}{2}$. 
        \end{lemma}
        \begin{proof}
            We consider the two events separately: 
            \begin{itemize}
                \item \emph{Upper bound on $\mathbb{P}(E_{n,1}^c)$}. We note that the terms $\{\delta_t = \gtilde_t(Z_t) - \mathbb{E}[\gtilde(Z_t)|\mc{F}_{t-1}]: t \geq 1\}$ for a martingale difference sequence, taking values in the bounded interval $[-1,1]$. Hence, by an application of Azuma's inequality, we get the bound 
                \begin{align}
                    \mathbb{P}\lp \lv \frac{1}{n} \sum_{t=1}^n \delta_t \rv > a \rp \leq 2 \exp \lp -n \frac{a^2}{8} \rp. 
                \end{align}
                By setting $a = u_n = \sqrt{ 8 \log(4/\gamma)/n}$, we get the required result that $\mathbb{P}(E_{n,1}^c) \leq \gamma/2$. 
                
                \item \emph{Upper bound on $\mathbb{P}(E_{n,2}^c)$}.  Introduce the event $E_{n,2}^{(t)} \defined \{ |\mathbb{E}[\gtilde_t(Z_t)|\mc{F}_{t-1}]-\mathbb{E}[\gtilde^*(Z_t)]| \leq v_{n,t} \}$, where $v_{n,t} \defined 2C_{t-1}(\Gtilde) + 5\sqrt{2 \log(16n/\gamma)/(t-1)}$ for $t \geq 2$, and $v_{n,1}=2$. Note that $v_n = (\sum_{t=1}^n v_{n,t})/n$ and $E_{n,2}^c \subset \cup_{t=1}^n \lp E_{n,2}^{(t)}\rp^c$, and thus, to complete the proof, it suffices to show the following (since the result then follows by a union bound):  
                \begin{align}
                    \mathbb{P}\lp \lp  E_{n,2}^{(t)}\rp^c \rp  = \mathbb{P}\lp |\mathbb{E}[\gtilde_t(Z_t)|\mc{F}_{t-1}]-\mathbb{E}[\gtilde^*(Z_t)]| > 2C_{t-1}(\Gtilde) + 5\sqrt{ \frac{ 2\log(16n/\gamma)}{t-1}} \rp \leq \frac{\gamma}{2n}. \label{eq:proof-abstract-16}
                \end{align}
                Since $\gtilde_t$ corresponds to the ERM predictor based on $\{Z_1, \ldots, Z_{t-1}\}$, ~\eqref{eq:proof-abstract-16} follows by an application of \textcite[Theorem~26.5]{shalev2014understanding}. The remaining case of $t=1$ follows trivially, since $\gtilde_1$ and $\gtilde^*$ take values in $[-1,1]$, and thus $|\mathbb{E}[\gtilde_1(Z_1)|\mc{F}_0] - \mathbb{E}[\gtilde^*(Z_1)]|\leq 2 = v_{n,1}$ almost surely. 
            \end{itemize}
        \end{proof}

    \section{Proofs of Auxiliary Lemmas} 
    \label{appendix:proof-auxiliary} 

        \subsection{Proof of~\Cref{lemma:alt-1}} 
            \label{proof:alt-1}

                \noindent \textbf{Event $\boldsymbol{G_{n,1}}$:} Let $v^*_t$ denote the random variable $g^*(X_t)-g^*(Y_t)$. Then, we have $\mathbb{E}[v^*_t] = \Delta$ for all $t\geq 1$, and $\mathbb{V}(v^*_t) = \mathbb{E}[(v^*_t)^2] - \Delta^2 \leq \sigmatilde^2$. Hence, by a direct application of Bernstein's inequality, we get for any $\delta>0$:
                \begin{align}
                    \mathbb{P}\lp \frac{1}{n}\sum_{t=1}^n v_t^* - \Delta  < -\delta \rp \leq \exp \lp \frac{\frac{1}{2}n \delta^2}{ \sigmatilde^2 + \frac{1}{3}  \delta} \rp. 
                \end{align}
                On inverting this inequality, we get 
                \begin{align}
                    \mathbb{P}\lp \frac{1}{n} \sum_{t=1}^n v_t^* < \Delta - \sigmatilde \sqrt{ \frac{4 \log n}{n}} - \frac{2 \log n}{3 n} \rp \leq \frac{1}{n^2}. 
                \end{align}
                Hence, we have proved that $\mathbb{P}(G_{n,1}^c) \leq \frac{1}{n^2}$. 
                
                \noindent \textbf{Event $\boldsymbol{G_{n,2}}$:}
                The analysis of this event proceeds by considering the sequence of  bounded random variables, $v_t^2 - \mathbb{E}[v_t^2|\mc{F}_{t-1}]$. Introduce the shorthand $\beta_t$ for the random variable $\mathbb{E}[v_t^2|\mc{F}_{t-1}]$, and observe the following: 
                \begin{itemize}
                    \item For  all $t \geq 1$, we have $\beta_t \leq \sigmatilde^2$ almost surely, since $\beta_t = \mathbb{E}[ v_t^2 |\mc{F}_{t-1}] = \mathbb{E}[ (g_t(X_t) - g_t(Y_t))^2|\mc{F}_{t-1}] \leqas \sup_{g \in \mc{G}} \mathbb{E}[ (g(X) - g(Y))^2] \defined \sigmatilde^2$. Here, we have used the facts that $g_t$ is $\mathcal{F}_{t-1}$-measurable, and $(X_t, Y_t) \perp \mc{F}_{t-1}$.
                    \item For all $t \geq 1$, we have $\mathbb{E}[(v_t^2-\beta_t)^2|\mc{F}_{t-1}] \leq \gamma^2$ almost surely. This is because 
                    \begin{align}
                        \mathbb{E}\lb (v_t^2 - \beta_t)^2 | \mc{F}_{t-1} \rb &= \mathbb{E}\lb v_t^4 | \mc{F}_{t-1} \rb + \beta_t^2 - 2 \beta_t \mathbb{E}[v_t^2|\mc{F}_{t-1}] = \mathbb{E}\lb v_t^4|\mc{F}_{t-1} \rb - \beta_t^2  \\
                        & \leq \mathbb{E}\lb v_t^4 | \mc{F}_{t-1}\rb = \mathbb{E}\lb \lp g_t(X_t) - g_t(Y_t) \rp^4 \mid \mc{F}_{t-1} \rb \stackrel{a.s.}{\leq} \gamma^2. 
                    \end{align}
                    The last inequality above uses the fact that $g_t$ is $\mc{F}_{t-1}$-measurable, and that $(X_t, Y_t) \perp \mc{F}_{t-1}$. 
                \end{itemize}
        
                An application of  Freedman's inequality for the sequence of conditional zero-mean random variables, $(v_t^2 - \beta_t)_{t \geq 1}$, implies the following, with $V_n = \frac{1}{n} \sum_{t=1}^n v_t^2$: 
                \begin{align}
                     \mathbb{P}\lp V_n - \frac{1}{n} \sum_{t=1}^n \beta_t > \delta \rp  = \mathbb{P}\lp \sum_{t=1}^n v_t^2 - \beta_t > n \delta \rp \leq \exp \lp \frac{ \frac{1}{2} n^2 \delta^2}{  n\gamma^2 + \frac{1}{3}n\delta} \rp. \label{eq:alt-proof-17}
                \end{align}
               
                On inverting this concentration inequality, we get 
                \begin{align}
                    \mathbb{P}\lp V_n - \red{\sigmatilde^2} > \gamma \sqrt{ \frac{4 \log n}{n}} + \frac{2\log(n)}{3n} \rp   \leq 
                    \mathbb{P}\lp V_n - \frac{1}{n} \sum_{t=1}^n \beta_t > \gamma \sqrt{ \frac{4 \log n}{n}} + \frac{2 \log(n)}{3n} \rp  \leq  \frac{1}{n^2}. \label{eq:alt-proof-18}
                \end{align}
                The first inequality uses the fact that $\beta_t \leq \sigmatilde^2$ almost surely for all $t \geq 1$.

                Combining~\eqref{eq:alt-proof-17} and~\eqref{eq:alt-proof-18}, we get the required $\mathbb{P}(G_n^c) \leq \mathbb{P}(G_{n,1}^c) + \mathbb{P}(G_{n,2}^c) \leq \frac{2}{n^2}$.  \hfill\qedsymbol
        
        \subsection{Proof of~\Cref{lemma:alt-2}} 
        \label{proof:alt-2}
             Observe that 
             \begin{align}
                 \sigmatilde^2 &= \sup_{g \in \mathcal{G}} \mathbb{E}\lb (g(X) - g(Y))^2 \rb && \\
                 &= \sup_{g \in \mathcal{G}} \lp \mathbb{V}\big(g(X) - g(Y) \big) + \mathbb{E}\lb g(X) - g(Y) \rb^2 \rp && \\
                 & \leq \sup_{g \in \mathcal{G}} \mathbb{V}\lp g(X) - g(Y) \rp + \sup_{g \in \mathcal{G}} \mathbb{E}[g(X) - g(Y)]^2 && \\
                 & = \sigma^2 +  \sup_{g \in \mathcal{G}} |\mathbb{E}[g(X) - g(Y)]|^2 && \text{(Defn. of $\sigma^2$)}\\
                 & = \sigma^2 +  \lp \sup_{g \in \mathcal{G}} |\mathbb{E}[g(X) - g(Y)]| \rp^2 && \text{(monotonicity of $x \mapsto x^2$ for positive $x$)} \\
                 & = \sigma^2 +  \lp \sup_{g \in \mathcal{G}} \mathbb{E}[g(X) - g(Y)] \rp^2 && \text{($\mathcal{G}$ is closed under negation)} \\
                 & = \sigma^2 + \Delta^2. && \text{(Defn. of $\Delta$)} 
             \end{align}
             This proves the first claimed inequality.
             To prove the second inequality, we proceed as follows: 
             \begin{align}
                \sigmatilde^2 \leq \sigma^2 + \Delta^2 \leq \sigma^2 + \Delta^2 + 2 \sigma \Delta = (\sigma + \Delta)^2. 
             \end{align}
             On taking square roots and noting the quantities are positive, we get the required $\sigmatilde \leq \sigma + \Delta$.  By the same reasoning, we can also conclude that $\sqrt{\sigmatilde} \leq \sqrt{\sigma} + \sqrt{\Delta}$: 
             \begin{align}
                 (\sqrt{\sigmatilde})^2 \leq (\sqrt{\sigma})^2 + (\sqrt{\Delta})^2 \leq (\sqrt{\sigma})^2 + (\sqrt{\Delta})^2 + 2 \sqrt{\sigma \, \Delta} = (\sqrt{\sigma} + \sqrt{\Delta})^2,
             \end{align}
             which on taking square roots implies $\sqrt{\sigmatilde} \leq \sqrt{\sigma} + \sqrt{\Delta}$. \hfill \qedsymbol

        \subsection{Proof of~\Cref{lemma:n01}}
        \label{proof:n01}

Recall that $\sqrt{\sigmatilde} \leq \sqrt{\sigma} + \sqrt{\Delta}$ by~\eqref{eq:sigmatilde-inequalities} in~Lemma~\ref{lemma:alt-2}, which implies 
\[
\sqrt{2 \sigmatilde} \lp \frac{\log (n/\alpha)}{n} \rp^{3/4} \leq \sqrt{2 \sigma} \lp \frac{\log (n/\alpha)}{n} \rp^{3/4} + \sqrt{2 \Delta} \lp \frac{\log (n/\alpha)}{n} \rp^{3/4}. 
\]
This allows us to bound $n_{0,1}$ as 
\begin{align}
n_{0,1} \leq \inf {\bigg\{} n \geq 1:& \frac{\Delta}{2} \geq r_n + 9 \sigma \sqrt{\frac{2 \log (n/\alpha)}{n}} + 7\sqrt{2\sigma}\lp \frac{\log (n/\alpha)}{n} \rp^{3/4}   \\
& \;+ 7\sqrt{2\Delta} \lp \frac{\log (n/\alpha)}{n} \rp^{3/4}  +    \frac{8\log (n/\alpha)}{n}  {\bigg\}}.  
\end{align}

Let us define $m_1$ as 
\begin{align}
m_1 &= \inf \lbr n \geq 1:  7\sqrt{2\Delta} \lp \frac{\log (n/\alpha)}{n} \rp^{3/4}  \leq \frac{\Delta}{8} \rbr 
 = \inf \lbr n \geq 1: \frac{n}{\log (n/\alpha)} \geq \frac{ (6272)^{2/3}}{\Delta^{2/3}} \rbr, 
\end{align} 
and observe that for $n \geq m_1$, we have $7 \sqrt{2 \Delta}\lp \frac{\log (n/\alpha)}{n}\rp^{3/4} \leq \Delta/8$. Hence, we get 
\begin{align}
    n_{0,1} \leq \inf \lbr n \geq m_1: \frac{3\Delta}{8} \geq r_n + 9 \sigma \sqrt{\frac{2 \log (n/\alpha)}{n}} + 7\sqrt{2\sigma}\lp \frac{\log (n/\alpha)}{n} \rp^{3/4}  +    \frac{8\log (n/\alpha)}{n}  \rbr.  
\end{align}

To further analyze this term, we consider two cases: 

\paragraph{Case I: $\boldsymbol{\sigma^2 \geq \Delta}$.} Define $m_2$ as 
\begin{align}
    m_2 &\coloneqq \inf \lbr n \geq 1: 7 \sqrt{2\sigma} \lp \frac{\log(n/\alpha}{n} \rp^{3/4} \leq \frac{\Delta}{8} \rbr \\
    & = \inf \lbr n \geq 1: \frac{n}{\log (n/\alpha)} \geq \lp\frac{ 6272 \sigma}{\Delta^2}\rp^{2/3} \rbr  \\
    & \leq \inf \lbr n \geq 1: \frac{n}{\log (n/\alpha)} \geq 6272^{2/3} \frac{ \sigma^2}{\Delta^2} \rbr.
\end{align}
The inequality uses the fact that $(\sigma^2/\Delta^2) \geq (\sigma/\Delta^2)^{2/3}$ under the assumption that $\sigma^2 \geq \Delta$, or equivalently, that $\sigma^2/\Delta \geq 1$. To see this, we proceed as follows: 
\begin{align}
    \frac{\sigma^2}{\Delta^2} = \lp \frac{\sigma^3}{\Delta^3} \rp^{2/3} = \lp \frac{\sigma^2}{\Delta} \frac{\sigma}{\Delta^2} \rp^{2/3} \geq \lp \frac{\sigma}{\Delta^2} \rp^{2/3}. 
\end{align}

Now, for $n \geq m_2$, we have $7 \sqrt{2\sigma} \lp \frac{\log(n/\alpha}{n} \rp^{3/4} \leq \frac{\Delta}{8} $ by definition. Hence, we can further upper bound $n_{0,1}$ as 
\[
n_{0,1} \leq \inf \lbr n \geq m_1 \vee m_2: \frac{\Delta}{4} \geq r_n + 9 \sigma \sqrt{\frac{2 \log (n/\alpha)}{n}} +  \frac{8\log (n/\alpha)}{n}  \rbr.  
\]

\paragraph{Case II: $\boldsymbol{\sigma^2 < \Delta}$.} In this case, we have $\sqrt{\sigma} < \Delta^{1/4}$, which implies 
\begin{align}
    n_{0,1} \leq \inf \lbr n \geq m_1: \frac{3\Delta}{8} \geq r_n + 9 \sigma \sqrt{\frac{2 \log (n/\alpha)}{n}} + 7\sqrt{2} \Delta^{1/4}\lp \frac{\log (n/\alpha)}{n} \rp^{3/4}  +    \frac{8\log (n/\alpha)}{n}  \rbr.  
\end{align}
Define $m_3$ as 
\begin{align}
    m_3 &\coloneqq \inf \lbr n \geq 1: 7 \sqrt{2} \Delta^{1/4} \lp \frac{\log(n/\alpha}{n} \rp^{3/4} \leq \frac{\Delta}{8} \rbr = \inf \lbr n \geq 1: 7 \sqrt{2}  \lp \frac{\log(n/\alpha}{n} \rp^{3/4} \leq \frac{\Delta^{3/4}}{8} \rbr  \\
    & = \inf \lbr n \geq 1: \frac{n}{\log (n/\alpha)} \geq 6272^{2/3} \frac{1}{\Delta} \rbr.
\end{align}
Observe that $m_3 \geq m_1$ as $\Delta \leq 1$, and for $n \geq m_3$, we have $7 \sqrt{2\sigma} \lp \log(n/\alpha)/n \rp^{3/4} \leq \Delta/8$. Hence, we can again bound $n_{0,1}$ as 
\begin{align}
n_{0,1} &\leq \inf \lbr n \geq m_3\vee m_1: \frac{\Delta}{4} \geq r_n + 9 \sigma \sqrt{\frac{2 \log (n/\alpha)}{n}} +  \frac{8\log (n/\alpha)}{n}  \rbr  \\
        &= \inf \lbr n \geq m_3: \frac{\Delta}{4} \geq r_n + 9 \sigma \sqrt{\frac{2 \log (n/\alpha)}{n}} +  \frac{8\log (n/\alpha)}{n}  \rbr,  
\end{align}
where the equality uses the fact that $m_3 \geq m_1$.  \hfill \qedsymbol

        \subsection{Proof of~\Cref{lemma:logn-over-n-bound}}
            If $N<21$, then the statement of the theorem is trivially true. Hence we consider the $N\geq21$ case. 
            First, we note that for all $n \geq 20$, we have $\log(bn)/n \leq \log(4n)/n \leq 1$. If $N_1$ denotes the real-valued solution of the equation $\log(bn)/n = a$, then since  $N-1 \leq N_1 \leq N$, we get the following: 
            \begin{align}
                &a = \frac{\log b N_1}{N_1^{1/2}}  \times \frac{1}{N_1^{1/2}} \leq  \frac{1}{N_1^{1/2}},  \\
                 \text{which implies}  \quad & (bN_1)^{1/2}  \leq \frac{b^{1/2}}{a}\; \Rightarrow \log(bN_1) \leq 2 \log(b^{1/2}/a) \leq 2 \log(2/a), 
            \end{align}
            which implies that $N_1 \leq 2 \log(2/a)/a$. Since $N \leq N_1 + 1$, this completes the proof.\hfill\qedsymbol

    \section{Beyond IPMs}
        \label{appendix:beyound-IPMs}
        While our discussion in this paper focused on designing tests based on integral probability metrics~(IPMs), we note that similar arguments work with other distance measures that admit a variational representation, such as members of the $f$-divergence family. For instance, we can easily develop a sequential version of the popular $\chi^2$-test using our approach. 
        
        In addition, we  can also use our ideas to construct sequential tests using predictors that are learned incrementally using stochastic gradient descent~(or similar methods), such as neural networks. More specifically,  consider a class of functions $\mc{G}=\{g_{\theta}:\mc{X} \to [-1/2, 1/2], \text{ with } \theta \in \Theta\}$ for some parameter set $\Theta$, and define the distance measure  using a loss function $\ell$, as follows: 
        \begin{align}
            \dG(P_X, P_Y) \defined  \sup_{\theta \in \Theta} \mathbb{E}[\ell(\theta, X, Y)].  
        \end{align}
        Assume that $\dG(P_X, P_Y)=0$ under the null, and $\dG(P_X, P_Y)>0$ under the alternative. To construct a sequential  two-sample test using this $\dG$, we can use an SGD based prediction strategy, that sets $g_t = g_{\theta_t}$, and $\theta_t = \theta_{t-1} + \eta_t \ell'(\theta_{t-1}, X_t, Y_t)$ for some step size $\eta_t$. It is easy to verify that a sufficient condition for the resulting test to be consistent is that $\theta_t \convas \theta^*$, where $\theta^*$ is a random variable taking values in $\Theta$, such that $\mathbb{E}[\ell(\theta^*, X, Y)|\theta^*]>0$ almost surely. 
 
\end{appendix}
\end{document}